\documentclass[letterpaper,11pt,reqno]{amsart} 
\usepackage{sy}
\newcommand{\tR}{\mathsf{T}_R}
\newcommand{\tL}{\mathsf{T}_L}
\DeclareMathOperator{\typ}{typ}

\makeatletter
\renewcommand\paragraph{\@startsection{paragraph}{4}{\z@}%
  {3.25ex \@plus1ex \@minus.2ex}%
  {-1em}%
  {\normalfont\normalsize\bfseries}}
\makeatother
\begin{document}
\title[PA]{Cutoff for Simple Random Walk on\\
Preferential Attachment Graphs}

\author[A. Sly and S. Yang]
       {Allan Sly$^1$ and Seoyeon Yang$^1$}

\address{$^1$Department of Mathematics,
  Princeton University,
  Princeton, NJ 08540, USA}

\email{allansly@princeton.edu}
\email{syeon.y@princeton.edu}

%%%%%%%%%%%%%%%%%%%%%%%%%%%%%%%%%%%%%%%%%%%%%%%%%%%%
\begin{abstract}
We consider simple random walk on the sequential affine preferential attachment graphs with fixed
parameters $m\ge2$ and $\delta\ge0$. We prove the walk has total-variation cutoff, and for every fixed
$\varepsilon\in(0,1)$,
\[
t_{\mathrm{mix}}^{(n)}(\varepsilon)
=\frac{\log n}{\nu d}
 +O_{\mathbb P}\!\left(\sqrt{\log n}\right),
\]
where $\nu>0$ is the walk speed and $d>0$ is the entropy per level
of harmonic measure on the P\'olya--point tree, the local weak
limit of the graph. Thus cutoff occurs at the entropic time.

Entropic cutoff was previously established for a range of random graph models including configuration models under mild degree assumptions and  on Erd\H{o}s--R\'enyi random graphs for typical starting points.
We extend the harmonic-measure approach to preferential attachment,
accounting for its age-dependent structure and hubs of polynomial degree.
The local weak limit, the P\'olya--point tree, provides a more tractable
setting for analyzing the walk.
The proof combines regeneration estimates on the local limit tree, a quantitative graph--tree coupling on the mixing-time scale, and control of hub excursions.  We show that shortly before the cutoff time its distribution is close in total variation distance to a measure with small $L^\infty$ distance from stationarity.
When $m\ge2$, the spectral gap is bounded away from zero with
high probability which completes the upper-bound coupling argument.
% We prove worst-case total-variation cutoff for simple random walk
% on sequential affine preferential attachment graphs with fixed
% parameters $m\ge2$ and $\delta\ge0$. 
% %Each new vertex sends $m$ edges to older vertices, chosen sequentially with probabilities proportional to their current degrees plus $\delta$. 
% For every fixed
% $\varepsilon\in(0,1)$,
% \[
% t_{\mathrm{mix}}^{(n)}(\varepsilon)
% =\frac{\log n}{\nu d}
%  +O_{\mathbb P}\!\left(\sqrt{\log n}\right),
% \]
% where $\nu>0$ is the walk speed and $d>0$ is the entropy per level
% of harmonic measure on the P\'olya--point tree, the local weak
% limit of the graph. Thus cutoff occurs at the entropic time.

% Entropic cutoff was previously established for fixed-degree random
% regular graphs and for configuration models under degree-tail
% assumptions ensuring subpolynomial maximum degree. We adapt the
% harmonic-measure approach to preferential attachment despite its
% heavy-tailed degrees, polynomial-degree hubs, and age-dependent
% local structure. The proof combines regeneration estimates on the limiting tree, a quantitative graph--tree coupling on the mixing-time scale, and control of hub excursions. 
% We also establish, for each
% fixed $m\ge2$, an absolute spectral gap bounded away from zero with
% high probability, uniformly over $\delta\ge0$, and use it to complete
% the upper-bound argument.
\end{abstract}

\maketitle
% \setcounter{tocdepth}{1}
%\tableofcontents

%%%%%%%%%%%%%%%%%%%%%%%%%%%%%%%%%%%%%%%%%%%%%%%%%%%%
%%%%%%%%%%%%
%%%%%%%%%%%%
%%%%%%%%%%%%
%%%%%%%%%%%%
\section{Introduction}\label{sec:intro}

Preferential attachment is a widely studied model of growing networks with
power-law degree distributions~\cite{BarabasiAlbert,BollobasRiordanSpencerTusnady}. Vertices arrive successively and connect preferentially to vertices of high degree.
We study the mixing time of simple random walk on these graphs.
Our main result identifies the worst-case mixing time and proves
cutoff: the total-variation distance to equilibrium drops from near
one to near zero in a time interval negligible compared with the
mixing time. 
The cutoff time is of order \(\log n\), with the leading constant determined by the speed and harmonic-measure entropy of random walk on the local limiting tree.

We consider the sequential affine model, in which each new vertex
sends $m\ge2$ edges to older vertices, chosen sequentially with
probabilities proportional to their current degrees plus
$\delta\ge0$. 
The limiting probability of degree $k$ is asymptotic
to a positive constant times $k^{-(3+\delta/m)}$ as $k\to\infty$,
and fixed early vertices have degrees growing on the polynomial scale
$n^{m/(2m+\delta)}$ as the number $n$ of vertices
increases~\cite[Lemmas~5.1--5.2]{PA-local-limit}.
Thus the parameter range includes the infinite-variance case
$\delta=0$, while polynomial-degree hubs persist when $\delta>0$
and the limiting variance is finite. These hubs have much greater
stationary mass than a typical vertex, since the stationary
distribution is proportional to degree. Vertex ages also influence
the laws of neighboring vertices and the conditional distribution
of unrevealed edges.

\begin{figure}[ht]
    \centering
    \includegraphics[width=0.6\linewidth]{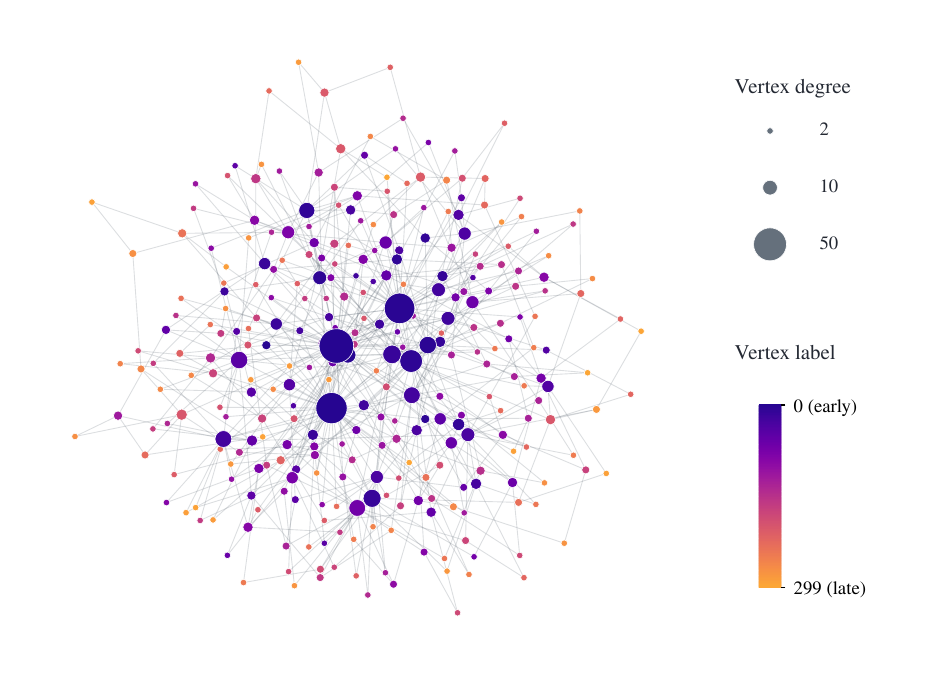}
    \caption{
    An illustrative preferential attachment variant.
Vertex area is proportional to degree (hence stationary mass);
color indicates birth order, from purple (early) to orange (late).}
    \label{fig:pa-degree-age}
\end{figure}

The local weak limit around a uniformly chosen vertex is the
P\'olya--point tree~\cite{PA-local-limit}. This is a random branching
tree whose offspring laws depend on auxiliary random variables which represent normalized birth times,
called ages, and on whether a vertex is older or younger than
its parent. Simple random walk on this tree is transient and {its loop erasure is an infinite ray.} The conditional law of that ray
given the tree is its \emph{harmonic measure}. 
This measure is far from uniform. The contrast with volume growth
is particularly pronounced in the critical case $\delta=0$:
neighborhoods in the P\'olya--point tree grow faster than
exponentially, and the diameter of the preferential attachment {is $O(\log n/\log\log n)$~\cite{BollobasRiordanDiameter}.}  
However, for every $\delta\ge0$, the harmonic measure is concentrated on a subtree with only exponential growth.
We let $d>0$ be the entropy per level of harmonic measure and so informally, a typical ray segment of length $k$, chosen according to harmonic measure,  has harmonic mass $\exp\{-dk+o(k)\}$.  Also, let $\nu>0$ denote
the speed of the walk. The walk reaches depth
about $\nu t$ at time $t$ so its entropy therefore
grows at rate $\nu d$. Setting $\exp(\nu dt)\approx n$ suggests the mixing scale.

Let $G_n$ be the preferential attachment graph defined formally in Definition~\ref{def:PA}. Let $P_n$ be its simple random walk transition
matrix and $\pi_n(v)=\deg(v)/(2m(n-1))$ its stationary distribution.
For an integer $t\ge0$ and $\varepsilon\in(0,1)$, write
\[
 d_n(t)=\max_{x\in V(G_n)}
       \|P_n^t(x,\cdot)-\pi_n\|_{\rm TV},
 \qquad
 t_{\rm mix}^{(n)}(\varepsilon)
       =\min\{t\ge0:d_n(t)\le\varepsilon\},
\]
and set $d_n(t)=d_n(\lfloor t\rfloor)$ for real $t\geq0$. We say that the walks exhibit cutoff at
$t_n^\star$ with window $w_n=o(t_n^\star)$ if, for every
$\varepsilon\in(0,1)$,
\begin{equation}\label{eq:cutoff-def}
\lim_{\gamma\to\infty}\liminf_{n\to\infty}
 \Pbb\bigl(d_n(t_n^\star-\gamma w_n)\ge1-\varepsilon\bigr)=1,
\qquad \lim_{\gamma\to\infty}\liminf_{n\to\infty}
 \Pbb\bigl(d_n(t_n^\star+\gamma w_n)\le\varepsilon\bigr)=1.
\end{equation}
Here $\Pbb$ denotes the graph law.

\begin{thm}[Cutoff for SRW on preferential attachment]\label{thm:main-cutoff}
Fix $m\ge2$ and $\delta\ge0$, and let $G_n$ be the preferential
attachment graph of Definition~\ref{def:PA}.
Let $\nu\in(0,1]$ and
$d\in(0,\infty)$ be the speed and the entropy per level of harmonic
measure of simple random walk on the P\'olya--point tree, defined in
\eqref{eq:def-speed} and \eqref{eq:def-entropy-rate}, and set
$t_n^\star:=\log n/(\nu d)$.
Simple random walk on $G_n$ exhibits
cutoff at $t_n^\star$ with a window of order at most
$\sqrt{\log n}$.
Equivalently,
for every fixed $\varepsilon\in(0,1)$,
\[
t_{\rm mix}^{(n)}(\varepsilon)=\frac{\log n}{\nu d}+O_{\Pbb}\bigl(\sqrt{\log n}\bigr).
\]
\end{thm}

The lower bound in Proposition~\ref{prop:cutoff-lower-bound}
already holds for a uniformly chosen starting vertex, whereas
the upper bound is uniform over all starting vertices, including
the hubs. Thus typical and worst starting vertices have the same
leading mixing time. The fluctuation estimates are used to obtain a cutoff window
of order at most $\sqrt{\log n}$.

\subsection{Related results}\label{subsec:related}
On random regular graphs of fixed degree at least three, simple
random walk exhibits cutoff when the distance traveled on the
limiting regular tree reaches the typical distance between graph
vertices~\cite{LubetzkySlyRegular}. Conditional on its distance
from the root, the tree walk is uniform on that level. For
nonregular trees, the distribution of mass within a level also
matters. Indeed, conditional on survival, harmonic measure on a
supercritical Galton--Watson tree with a nonconstant offspring
distribution of finite mean has dimension strictly smaller than
that of the tree boundary~\cite{LyonsPemantlePeres}.
The harmonic-measure dimension, together with the walk speed,
determines the cutoff location from a typical vertex of the
supercritical Erd\H os--R\'enyi giant and of configuration models
under suitable degree assumptions~\cite{RWRG}.
{
A similar harmonic-measure approach is used in~\cite{UnivCutoff} to describe the cutoff time by the speed and entropy of random walk for graphs with an added random matching, through an auxiliary quasi-tree.}

{
Entropic methods based on trajectory weights have also been
developed for non-backtracking walks on configuration
models~\cite{BenHamouSalezNBRW}, random walks on directed
configuration models~\cite{BCSDigraphs,CaiCaputoPerarnauQuattropani},
and more general sparse random Markov chains~\cite{BCSEntropic}.
On the forward tree approximations used in these arguments,
the probability of reaching a vertex is the product of the
original walk's one-step transition probabilities along its
ancestral path. In our undirected setting, returning excursions
make harmonic weights depend on escape probabilities in the
descendant subtrees.
}

Passing from a typical starting vertex to the worst one requires
additional control of exceptional neighborhoods.
Worst-start cutoff was established for configuration models with
minimum degree at least three, a finite second degree moment, and
a degree-tail condition giving a sufficiently small subpolynomial
maximum degree~\cite{Comparing}.
{
Their proof uses the degree assumptions to control small
neighborhoods around every vertex, then applies tree couplings
and concentration estimates to the walk beyond these neighborhoods.} 
In preferential attachment, polynomial-degree hubs violate the neighborhood-size bounds used in this argument. Moreover, uniform pairing of half-edges is replaced by
attachment probabilities that depend on ages and previous
revelations. The case $m=2$ also
falls outside the minimum-degree-three hypothesis and requires
using the structure imposed by attachment.

For preferential attachment, earlier work establishes several
geometric estimates relevant to mixing. Constant conductance was
proved for a standard degree-proportional
model~\cite{alp-zero-PA-Cheeger}. A logarithmic mixing
bound for the lazy walk is used in the cover-time
analysis~\cite{CooperFriezeCover}, while large-set expansion for the present
model was proved in~\cite{alp-nonzero-PA}. We establish an absolute
spectral gap with high probability, uniformly over $\delta\ge0$, {showing that there are no eigenvalues too close to -1.}
%The cutoff theorem gives the finer information needed to identify the leading mixing constant and bound the width of the transition.

Shortest-path asymptotics provide a complementary comparison.
For the loop-allowing model at $\delta=0$, the diameter is
asymptotic to $\log n/\log\log n$~\cite{BollobasRiordanDiameter}.
In the sequential model with $\delta>0$, typical distances are
logarithmic, with leading constant determined by the growth rate
of the local limit~\cite{vanderHofstadZhu2025}.
For every $\delta\ge0$, the mixing time in
Theorem~\ref{thm:main-cutoff} is logarithmic, with leading constant
determined by the walk speed and harmonic-measure entropy.

\subsection{Proof ideas}\label{subsec:strategy}

We follow the harmonic-measure approach of~\cite{RWRG},
beginning with quantitative estimates on the limiting tree.
For fixed $m\ge2$ and $\delta\ge0$, let $X$ be simple random
walk on the P\'olya--point tree $\Tcal$ with uniform root age,
started at its root $\rho$. Write $|X_t|$ for its distance from
the root at time $t$, and let $\xi=(\rho,\xi_1,\ldots)$ be its
limiting ray. For $v\in\Tcal$, define its \emph{harmonic information} by
\[
 W(v):=-\log\Pbb_\rho^{\Tcal}(v\in\xi),
 \qquad W_k:=W(\xi_k),
\]
where $\Pbb_\rho^{\Tcal}$ denotes the walk law conditional on
$\Tcal$. Thus $e^{-W_k}$ is the harmonic mass of rays sharing
the first $k$ edges of $\xi$.
Theorems~\ref{thm:distance-fluctuations}
and~\ref{thm:original-W-tightness} give
\[
 \max_{0\le s\le t}\bigl||X_s|-\nu s\bigr|
       =O_{\Pbb}(\sqrt t),
 \qquad
 W_k-dk=O_{\Pbb}(\sqrt k),
\]
under the joint law of the tree and the walk. Consequently,
$W_{|X_t|}=\nu dt+O_{\Pbb}(\sqrt t)$, providing the
fluctuation scale used in the mixing argument.

To handle the dependence on age, we construct a regeneration
structure with a fixed law for the age and type of the new root.
This yields the distance estimate. Harmonic weights still depend
on the tree beyond the next regeneration point. Approximating
escape probabilities at finite depth controls this dependence
and gives the covariance bounds needed for the entropy estimate.

%%%%%%%%%%%%%%%
To transfer estimates from the P\'olya--Point Tree to $G_n$ we need to construct a coupling between the graphs. We do not attempt to couple the entire ball of radius $\frac{1}{d}\log n$ in $G_n$ with the P\'olya--Point Tree as the number of vertices in the tree would be much larger than $n$.  Instead, we couple a truncated portion of the tree of depth $\frac{1-o(1)}{d}\log n$ that with high probability supports almost all of the harmonic measure.
We truncate the tree at points where it encounters either high degree vertices, very old vertices or along rays whose probability is too low. We show that few enough vertices are revealed along the way so that the P\'olya urn representation
of~\cite{PA-local-limit} continues to accurately describe the
subsequent attachment probabilities.
%The comparison is carried out on truncated neighborhoods that retain most of the walk mass while limiting the number of revealed vertices. The resulting continuation estimate, Proposition~\ref{prop:walk-coupling-log-time}, is designed for successive comparisons conditional on the preceding exploration.

Establishing the mixing bound from a typical starting point then follows from two consequences of the
coupling. For the lower bound, a walk from a uniformly chosen
vertex is localized on a set of $o(n)$ vertices with total degree $o(n)$, which has a vanishing stationary mass.
For the upper bound, we show that, near the cutoff time, the law of the walk is arbitrarily close in total variation to a probability measure with density relative to \(\pi_n\) at most \(\exp\{O(\sqrt{\log n})\}\).
The final step uses the uniform absolute spectral gap of
Proposition~\ref{prop:PA-uniform-two-sided-gap} which says that for each fixed $m\ge2$, there is $c=c(m)>0$ such that the spectral gap is at least $c(m)$.
Then, because of the spectral gap, after a further $O(\sqrt{\log n})$ steps, the $L^2$ distance from stationarity will be $o(1)$ which completes the proof.

Proving the result for arbitrary starting vertices requires a separate argument.
In~\cite{Comparing}, they extended the method of~\cite{RWRG} from typical starting points to prove cutoff in a configuration model with minimum degree 3 and degree distribution with exponential tails. They showed that after $O(\log \log n)$ steps, the walk is spread out enough so that with high probability it is at a vertex that serves as a typical starting point. In our setting, we allow a burn-in time of $O((\log\log n)^2)$, treating visits to high-degree vertices separately to control the amount of the graph revealed. This additional time is negligible compared with the $\sqrt{\log n}$ cutoff window.

When the walk hits a high degree vertex (greater than $\log^C n$) in the burn-in phase, it makes with high probability a bounded expected number of excursions before escaping down one of its subtrees.  Since the degree could be as large as $\sqrt{n}$, revealing all these subtrees would reveal too many vertices.  Instead, we sample $\lceil(\log n)^2\rceil$ incident slots and show that 
with probability $1-o(n^{-1})$, almost all subtrees will be typical and give the estimates we need.  Hence, with high probability it will escape down a typical subtree and we can apply the argument for mixing from a typical starting point.

Several questions remain. The lower bound from a uniformly chosen
vertex does not determine whether individual hubs mix earlier:
their many possible initial departures may accelerate initial mixing. The optimal window and any
limiting cutoff profile also remain to be determined.
The regime $-m<\delta<0$, corresponding to degree exponent between
two and three, is outside our scope. In particular, the regularity estimates in
Lemma~\ref{lem:regularity-PPT} are proved only for $\delta\ge0$.

Section~\ref{sec:PA-PPT} gives the model, urn estimates and spectral
bound. Section~\ref{sec:RegRWPA} develops the tree-walk estimates,
and Section~\ref{sec:coupling} constructs the comparison with the
finite graph. Section~\ref{sec:worst-start} treats arbitrary
starting vertices. Section~\ref{sec:cutoff-proof} proves the lower
bound and combines the estimates to establish
Theorem~\ref{thm:main-cutoff}.

\medskip

\section{Preferential Attachment and the P\'olya--Point Tree}\label{sec:PA-PPT}

\subsection{The finite preferential attachment graph}

Preferential attachment is a growth rule in which vertices with larger
current degree are more likely to receive new edges.

\begin{defn}[Sequential preferential attachment]\label{def:PA}
Fix \(m\geq 2\) and \(\delta\geq0\).

The graph \(G_2\) consists of vertices \(1,2\) joined by \(m\) parallel
edges.  For each \(k\geq 3\), form \(G_k\) from \(G_{k-1}\) by adding vertex
\(k\) and then adding \(m\) edges from \(k\) to older vertices, one at a
time.  For \(i=1,\ldots,m\), let \(\mathcal F_{k,i-1}\) be the history just
before the \(i\)-th edge is added, and let \(d_{k,i-1}(v)\) be the degree of
\(v<k\) at that time.  The older endpoint \(J_{k,i}\) is chosen according
to
\begin{equation}\label{eq:PA-rule-intro}
    \Pbb\bigl(J_{k,i}=v\mid\mathcal F_{k,i-1}\bigr)
    =\frac{d_{k,i-1}(v)+\delta}
         {2m(k-2)+\delta(k-1)+(i-1)},
    \qquad v<k.
\end{equation}
Endpoints chosen by different edges may coincide, so \(G_n\) is an
undirected loopless multigraph.  Degrees count edge multiplicity, and simple
random walk chooses an incident edge uniformly.
\end{defn}

Set
\(\alpha:=\delta/(2m+\delta)\in[0,1)\) and
\(u:=\alpha/(1-\alpha)=\delta/(2m)\).
Thus \(\alpha=0\) is pure degree-proportional attachment, whereas the
older endpoint tends to the uniform law as \(\alpha\uparrow1\).
The affine attractiveness is \(2mu\).  We also use the derived parameters
\[
    a:=m+2mu,
    \qquad
    \chi:=\frac{a}{m+a}=\frac{1+\alpha}{2}.
\]

\subsection{The P\'olya--point tree and the urn representation}

We now define the P\'olya--point tree (PPT), the local weak limit of \(G_n\).

\begin{defn}[P\'olya--point tree (PPT)]\label{def:PPT}
The P\'olya--point tree is a rooted random tree $(\Tcal,\rho)$ with ages
$A_v\in(0,1]$ and genealogical types $\typ(v)\in\{\tL,\tR\}$ for
non-root vertices. The root has no type and has age
$A_\rho\sim\Unif[0,1]$. Each vertex receives a mark with law
\[
    \Gamma_v\sim
    \begin{cases}
        \mathrm{Gamma}(a,1),   & v=\rho \text{ or } \typ(v)=\tR,\\
        \mathrm{Gamma}(a+1,1), & \typ(v)=\tL.
    \end{cases}
\]

Conditional on $(A_v,\Gamma_v)$, generate the following two families
independently.
\begin{enumerate}[label=\textup{(\roman*)},leftmargin=2.2em]
    \item \textbf{Left children.}
    There are $m$ left children if $v$ is the root or has type $\tL$,
    and $m-1$ if it has type $\tR$. Their ages are independent copies of
    $A_vV^{1/\chi}$, where $V\sim\Unif[0,1]$, and their type is $\tL$.

    \item \textbf{Right children.}
    The right-child ages form a Poisson point process on $[A_v,1]$ with intensity
    \[
        \rho_v(x)\dd x
        :=
        (1-\chi)\Gamma_v\frac{x^{-\chi}}{A_v^{1-\chi}}\dd x,
        \qquad x\in[A_v,1].
    \]
    The corresponding children have type $\tR$.
\end{enumerate}
At each step, new variables are sampled independently of the preceding
construction, subject to these conditional laws.
\end{defn}

Genealogical types and marks are retained when the tree is re-rooted;
the new root need not have the original root law.

If \(V_0\) is uniform on \([n]\), independently of \(G_n\), then
\((G_n,V_0)\) converges in the local weak sense to \((\Tcal,\rho)\)
\cite[Theorem~2.2]{PA-local-limit}. More precisely, for each fixed radius
and every $\varepsilon>0$, all sufficiently large $n$ admit a coupling
under which, with probability at least $1-\varepsilon$, the two rooted
balls are isomorphic and every matched pair $(k,v)$ satisfies
$|k/n-A_v|\leq\varepsilon$ \cite[Corollary~2.3]{PA-local-limit}.
Quantitative growing-radius comparisons are developed in
Section~\ref{subsec:ppt-pa-coupling}.

Let \(N_v\) be the number of right children of \(v\).  Conditional on
\((A_v,\Gamma_v)\),
\[
    N_v\sim\Poi(\Lambda_v),
    \qquad
    \Lambda_v:=\int_{A_v}^1\rho_v(x)\,\dd x
       =\Gamma_v\bigl(A_v^{-(1-\chi)}-1\bigr).
\]

\paragraph{P\'olya urn representation.}

We recall the P\'olya urn representation of the sequential preferential
attachment graph from \cite{PA-local-limit}.  Fix \(n\), set \(\psi_1=1\),
and let \(\{\psi_j\}_{2\le j\le n}\) be independent random variables with
\[
\psi_j\sim
\mathrm{Beta}\bigl(a,d_j\bigr),\qquad
d_j:=(2j-3)m+2mu(j-1),
\qquad
c_j:=a+d_j=2m(1+u)j-2m.
\]
For \(0\le k\le n\), define
\[
S_k:=\prod_{j=k+1}^n (1-\psi_j),
\]
with the conventions \(S_0=0\) and \(S_n=1\).  Equivalently, if
\[
    \varphi_k:=S_k-S_{k-1},
    \qquad
    I_k:=[S_{k-1},S_k),
\]
then \(\varphi_k=\psi_k S_k\) and the intervals \(I_1,\ldots,I_n\)
partition \([0,1)\).

Conditionally on \((\psi_j)_{1\le j\le n}\), let
\[
    \{U_{k,i}:2\le k\le n,\ 1\le i\le m\}
\]
be independent, with \(U_{k,i}\) uniform on \([0,S_{k-1}]\).  Draw the
slot \((k,i)\), the \(i\)-th edge of \(k\), to the unique \(j<k\) for which
\(U_{k,i}\in I_j\), ignoring the null event that \(U_{k,i}\) lies on an
interval boundary and retaining multiple edges.  Thus, conditionally on
\(\{\psi_j\}\),
\[
    p_{k\to j}:=\Pbb(J_{k,i}=j\mid\{\psi_\ell\})
    =\frac{\varphi_j}{S_{k-1}}
    =\psi_j\prod_{h=j+1}^{k-1}(1-\psi_h),
    \qquad j<k.
\]

The resulting multigraph has the same distribution as the sequential
preferential attachment graph \cite[Theorem 2.1]{PA-local-limit}.
The Beta parameters are the two attractiveness masses immediately after
vertex $j$ is added: $a$ at $j$ and $d_j$ on $[j-1]$. Subsequent
attachments restricted to $[j]$ reinforce these masses as a two-color
P\'olya urn.

The representation also explains Definition~\ref{def:PPT}. On the event
of Lemma~\ref{lem:S_k-estimate},
$S_{k-1}^{1/\chi}=(k/n)(1+O(\varepsilon+k^{-1}))$.
If $U_{k,i}\in I_j$, the two interval bounds similarly give
$U_{k,i}^{1/\chi}=(j/n)(1+O(\varepsilon+j^{-1}))$.
Writing $U_{k,i}=S_{k-1}V_{k,i}$, where $V_{k,i}$ is uniform, therefore
gives the limiting left-child age $A_vV_{k,i}^{1/\chi}$.
A uniform finite root has limiting age $\Unif[0,1]$.
Future vertices \(\ell>k\) attach to \(k\) when their uniform variables fall
in \(I_k\), whose length is \(\varphi_k=\psi_kS_k\). With $A=k/n$ and
$s=\ell/n$, the urn estimate and the approximation
$\psi_k\approx\Gamma/[2m(1+u)k]$ give
\[
    m p_{\ell\to k}
       \approx\frac1n(1-\chi)\Gamma A^{-(1-\chi)}s^{-\chi}.
\]
There are $n\,\dd s$ prospective younger vertices in $\dd s$, yielding
the right-child intensity in Definition~\ref{def:PPT} after Poisson
approximation. Conditioning a single attachment from a younger vertex
to end at $k$ weights the Beta density of $\psi_k$ by $\psi_k$, giving
$\mathrm{Beta}(a+1,d_k)$ and hence the limiting mark
$\mathrm{Gamma}(a+1,1)$. The outgoing attachments of a younger vertex
$\ell$ do not depend on $\psi_\ell$, so a right child
retains the base $\mathrm{Gamma}(a,1)$ law, as does the root.

\subsection{Quantitative urn and mark estimates}

The growing-ball coupling of Section~\ref{sec:coupling} uses three
quantitative inputs: Lemma~\ref{lem:S_k-estimate} controls urn endpoints
above a label cutoff; Lemma~\ref{lem:beta-gamma} compares Beta and Gamma
marks; and Lemma~\ref{lem:regularity-PPT} controls the ages, volume, and
marks in a PPT ball.

\begin{lem}[Uniform urn estimate]\label{lem:S_k-estimate}
There is a constant \(C=C(m,u)<\infty\) such that, for every
integer \(q\in\{3,\ldots,n\}\) and \(0<\varepsilon\leq1\),
\begin{equation}
    \Pbb\left(
       \max_{q-1\leq k\leq n}
       \left|\log S_k-\chi\log(k/n)\right|>\varepsilon
    \right)
    \leq \frac{C}{\varepsilon^2q}.
\label{eq:urn-failure}
\end{equation}
Consequently, on the complementary event,
\[
e^{-\varepsilon}\left(\frac{k}{n}\right)^\chi
    \leq S_k\leq
    e^{\varepsilon}\left(\frac{k}{n}\right)^\chi,
    \qquad q-1\leq k\leq n.
\]
\end{lem}

\begin{proof}
Let \(\Psi\) and \(\Psi_1\) denote the di- and trigamma functions.
For \(\psi_j\sim\mathrm{Beta}(a,d_j)\),
\[
    \Ebb\log(1-\psi_j)=\Psi(d_j)-\Psi(c_j),
    \qquad
    \Var\!\left(\log(1-\psi_j)\right)
       =\Psi_1(d_j)-\Psi_1(c_j).
\]
Their standard large-argument expansions give, uniformly in \(j\geq2\),
\[
    \left|\Ebb\log(1-\psi_j)+\frac{\chi}{j}\right|
       \leq \frac{C}{j^2},
    \qquad
    \Var\!\left(\log(1-\psi_j)\right)\leq\frac{C}{j^2}.
\]
The first bound and the harmonic-sum estimate show that the deterministic
bias is at most \(C/q\), uniformly for \(k\geq q-1\).  Kolmogorov's
maximal inequality, applied to the reversed centered sums, bounds the
probability that their maximum exceeds \(\varepsilon/2\) by
\(C/(\varepsilon^2q)\).  If \(\varepsilon<2C/q\), the same conclusion is
trivial after increasing \(C\).  This proves \eqref{eq:urn-failure}.
\end{proof}

For \(i>0\), let \(\mu_i\) denote the \(\mathrm{Gamma}(i,1)\) law and
let \(G_i\) be its distribution function.  Let \((Z_k)_{k\geq1}\) be
independent with common law \(\mu_a\); the variable \(Z_1\) is auxiliary.
For \(k\geq2\), define
the increasing quantile map \(f_k:[0,\infty)\to[0,1]\) by
\[
    \Pbb(\psi_k\leq f_k(z))=G_a(z),
    \qquad z\geq0.
\]
We henceforth realize \(\psi_k=f_k(Z_k)\), which has the required
\(\mathrm{Beta}(a,d_k)\) law.  The following lemma gives the quantitative
beta--Gamma comparison needed below.
For \(k\geq2\) and \(s>0\), define the good-mark set
\[
    \mathcal G_k^{\mathrm{BG}}(s)
       :=\left\{
          z>0:
          \left|\frac{f_k(z)}{z/[2m(1+u)k]}-1\right|\leq s
       \right\}.
\]

\begin{lem}[Uniform beta--Gamma approximation]\label{lem:beta-gamma}
There is a constant \(C=C(m,u)<\infty\) with the following property.
For any \(k_n^\circ\in\{3,\ldots,n\}\) satisfying
\(k_n^\circ\to\infty\), set
\begin{equation}
    \varepsilon_n^{\mathrm{BG}}
       :=C(k_n^\circ)^{-\frac{a}{2a+1}}.
\label{eq:beta-gamma-error}
\end{equation}
Then, for all sufficiently large \(n\), uniformly for
\(k_n^\circ\leq k\leq n\),
\begin{equation}
    \Pbb\biggl(
       \left|\frac{f_k(Z_k)}{Z_k/[2m(1+u)k]}-1\right|
          \leq\varepsilon_n^{\mathrm{BG}}
    \biggr)
    \geq1-\varepsilon_n^{\mathrm{BG}}.
\label{eq:beta-gamma-probability}
\end{equation}
The proof also gives the estimate
\begin{equation}
    \Ebb\left|c_kf_k(Z_k)-Z_k\right|\leq\frac{C}{k},
    \qquad k\geq3.
\label{eq:beta-gamma-L1}
\end{equation}
\end{lem}

\begin{proof}
Let \(p_k\)
and \(g_a\) denote the densities of \(c_k\mathrm{Beta}(a,d_k)\) and \(\mathrm{Gamma}(a,1)\) respectively.  For \(0<y<c_k\),
\[
    \log\frac{p_k(y)}{g_a(y)}
    =\log\frac{\Gamma(c_k)}{\Gamma(d_k)c_k^a}
      +y+(d_k-1)\log\left(1-\frac{y}{c_k}\right).
\]
The Gamma-ratio and Taylor expansions give, uniformly for
\(0<y\leq d_k^{1/4}\),
\[
    \left|\log\frac{p_k(y)}{g_a(y)}\right|
       \leq \frac{C(1+y+y^2)}{d_k}.
\]
Since the right-hand side is \(o(1)\), integration against \(yg_a(y)\)
bounds \(\int_0^{d_k^{1/4}}y|p_k-g_a|\,dy\) by \(C/d_k\).
Both laws have uniformly bounded sixth moments, using
\(\Ebb[(c_k\psi_k)^6]=c_k^6(a)_6/(c_k)_6\leq(a)_6\).
Their remaining first moments are therefore \(O(d_k^{-5/4})\), so
\[
    W_1\bigl(\mathcal L(c_k\psi_k),\mu_a\bigr)
       \leq\int_0^\infty y|p_k(y)-g_a(y)|\,dy
       \leq\frac{C}{d_k}\leq\frac{C}{k},
\]
which proves \eqref{eq:beta-gamma-L1}.

Since \(c_k/[2m(1+u)k]=1+O(k^{-1})\), whenever \(C/k\leq s\leq1\),
\(z\notin\mathcal G_k^{\mathrm{BG}}(s)\) implies
\(|c_kf_k(z)-z|\geq csz\) for some \(c=c(m,u)>0\).
For \(0<r\leq1\), Markov's inequality,
\eqref{eq:beta-gamma-L1}, and the Gamma small-ball estimate give
\begin{align*}
    \Pbb\bigl(Z_k\notin\mathcal G_k^{\mathrm{BG}}(s)\bigr)
       \leq\Pbb(Z_k<r)
          +\Pbb\bigl(|c_kf_k(Z_k)-Z_k|\geq csr\bigr)
       \leq Cr^a+\frac{C}{ksr}.
\end{align*}
Taking \(r=(ks)^{-1/(a+1)}\) gives the bound \(C(ks)^{-a/(a+1)}\).
For \(s=C_0(k_n^\circ)^{-a/(2a+1)}\), this is at most \(s\),
uniformly for \(k\geq k_n^\circ\), when \(C_0\) is sufficiently large.
Taking this constant in \eqref{eq:beta-gamma-error} proves
\eqref{eq:beta-gamma-probability}.
\end{proof}

%Exponential Markov's inequality and an integral comparison also give the following uniform tail estimate for the base Gamma variables. Using $\Ebb[e^{Z_k/2}]=2^a$, for every integer \(q\geq3\),
%\[
%    \Pbb\left(\exists k\geq q:
%       Z_k>(\log k)^2\right)
%    \leq 2^a\sum_{k=q}^\infty
%       \exp\left\{-\frac12(\log k)^2\right\}
 %   \leq C\exp\{-c(\log q)^2\}.
%\]

\subsection{Regularity of growing PPT balls}

For \(r\geq0\), let \(B_r(\rho)\) denote the radius-\(r\) ball in the
PPT, and write \(\Pbb_{a_0}\) for its law with deterministic root age
\(A_\rho=a_0\).

\begin{lem}[Regularity of the PPT]\label{lem:regularity-PPT}
Fix \(c_0,c_r,\zeta>0\), put
\[
    \ell_n:=\log\log n,\qquad r_n:=\lceil c_r\ell_n\rceil.
\]
There are constants \(C,c_{\mathrm{reg}},C_\Gamma>0\) and \(n_0\), depending only
on \(m,u,c_0,c_r,\zeta\), such that, for every \(n\geq n_0\) and
\(a_0\in[e^{-c_0\ell_n^2},1]\),
\begin{equation}\label{eq:deterministic-root-regularity}
    \Pbb_{a_0}\left(
       \begin{gathered}
          |B_{r_n}(\rho)|\leq a_0^{-(1-\chi)}e^{\zeta\ell_n^2},\\
          \min_{v\in B_{r_n}(\rho)}A_v\geq a_0e^{-\zeta\ell_n^2},\\
          \max_{v\in B_{r_n}(\rho)}\Gamma_v\leq C_\Gamma\ell_n^2
       \end{gathered}
    \right)
    \geq 1-Ce^{-c_{\mathrm{reg}}\ell_n^2}.
\end{equation}
Moreover, on this event, every \(v\in B_{r_n}(\rho)\) satisfies
\[
    \|\rho_v\|_\infty
       \leq (1-\chi)C_\Gamma\ell_n^2
          a_0^{-1}e^{\zeta\ell_n^2},
    \qquad
    \Lambda_v
       \leq C_\Gamma\ell_n^2
          a_0^{-(1-\chi)}
          e^{(1-\chi)\zeta\ell_n^2}.
\]
\end{lem}

\begin{proof}
Put \(\beta=1-\chi\leq\chi\) and
\(b=a_0e^{-\zeta\ell_n^2}\). The weight \(A^{-\beta}\) matches the
age dependence of the right-child intensity. Explore the PPT to depth
\(r_n\), retaining vertices of age below \(b\) but revealing none of
their offspring. Write \(\widehat B\) for the resulting
subtree and set
\[
    W:=\sum_{v\in\widehat B}(A_v\vee b)^{-\beta}.
\]
Each vertex contributes at least one, and a vertex crossing the age
cutoff contributes exactly \(b^{-\beta}\).

Consider a vertex of age \(A\geq b\), before its mark and offspring
are sampled. For \(V\sim\Unif[0,1]\), the expected weight of a left
child is
\[
    \Ebb\bigl[(AV^{1/\chi}\vee b)^{-\beta}\bigr]
    =A^{-\beta}\left(1+\beta\int_0^{\log(A/b)}
                    e^{(\beta-\chi)s}\,\dd s\right)
    \leq A^{-\beta}\bigl(1+\beta\log(A/b)\bigr).
\]
There are at most \(m\) left children. Since the mark has mean at most
\(a+1\) and right-child ages are at least \(A\), their expected total
weight is at most
\[
    (a+1)\beta A^{-\beta}
       \int_A^1 x^{-\chi-\beta}\,\dd x
    =(a+1)\beta A^{-\beta}\log(1/A).
\]
Thus the expected total child weight is at most
\(C(1+\log(1/b))\) times the parent's weight. Iterating this conditional
expectation bound over generations gives
\begin{equation}
\begin{split}
    \Ebb_{a_0}W
    \leq a_0^{-\beta}\sum_{j=0}^{r_n}
                  \bigl[C(1+\log(1/b))\bigr]^j
    \leq a_0^{-\beta}\exp\{O(\ell_n\log\ell_n)\},
\end{split}
\label{eq:regularity-weighted-moment}
\end{equation}
since \(\log(1/b)\leq(c_0+\zeta)\ell_n^2\) and \(r_n=O(\ell_n)\).
All bounds are uniform over the allowed root ages.

If the full ball contains a vertex of age below \(b\), the first such
vertex on its ancestral line belongs to \(\widehat B\). Hence
Markov's inequality and \eqref{eq:regularity-weighted-moment} give
\[
    \Pbb_{a_0}\left(\min_{v\in B_{r_n}(\rho)}A_v<b\right)
    \leq b^\beta\Ebb_{a_0}W
    \leq\exp\{-\beta\zeta\ell_n^2+O(\ell_n\log\ell_n)\}.
\]
On the complementary event, \(\widehat B=B_{r_n}(\rho)\) and
\(|\widehat B|\leq W\). The same estimate therefore bounds the
probability that the age bound holds but the volume bound fails by
\(\exp\{-\zeta\ell_n^2+O(\ell_n\log\ell_n)\}\).

Membership in \(\widehat B\) is determined before a vertex's own
mark is sampled. Conditional on its age and type, this mark has law
\(\mathrm{Gamma}(a,1)\) or \(\mathrm{Gamma}(a+1,1)\), so
\[
    \Pbb_{a_0}\left(
        \exists v\in\widehat B:\Gamma_v>z\right)
    \leq2^{a+1}e^{-z/2}\Ebb_{a_0}|\widehat B|
    \leq2^{a+1}e^{-z/2}\Ebb_{a_0}W.
\]
Since \(a_0^{-\beta}\leq e^{\beta c_0\ell_n^2}\), taking
\(z=C_\Gamma\ell_n^2\) with \(C_\Gamma>2\beta c_0+2\) bounds this
probability by \(e^{-\ell_n^2/2}\) for large \(n\).
As \(\ell_n\log\ell_n=o(\ell_n^2)\), these three estimates prove
\eqref{eq:deterministic-root-regularity}.

The remaining bounds follow from
\(\|\rho_v\|_\infty=(1-\chi)\Gamma_v/A_v\) and the formula for
\(\Lambda_v\).
\end{proof}

\subsection{A uniform spectral gap}
\label{sec:spectral}

Let $P_n$ be the transition matrix of simple random walk on $G_n$,
with eigenvalues $1=\lambda_1(P_n)\ge\cdots\ge\lambda_n(P_n)$.
Throughout this section, $m\ge2$ is fixed and all edge counts include
multiplicity. Let $\mathbb P_\alpha$ be the law of the model with
$\delta=2m\alpha/(1-\alpha)$, and let $\mathbb E_\alpha$ denote the
corresponding expectation.

\begin{prop}[Uniform two-sided spectral gap]
\label{prop:PA-uniform-two-sided-gap}
There is a constant $c=c(m)>0$ such that
\[
    \sup_{\alpha\in[0,1)}\mathbb P_\alpha\left(
        \max_{2\le i\le n}|\lambda_i(P_n)|>1-c
    \right)=o(1).
\]
\end{prop}

For a graph $G$ and $S\subseteq V(G)$, let $b_G(S)$ count the boundary edges
with exactly one endpoint in $S$, and let $e_G(S)$ count edges with
both endpoints in $S$. Write $d_G(v)$ for the degree of $v$ in $G$ and
$\operatorname{vol}_G(S)=\sum_{v\in S}d_G(v)$. Its conductance is
\[
    h(G)=\min_{\varnothing\ne S\subsetneq V(G)}
       \frac{b_G(S)}{
       \min\{\operatorname{vol}_G(S),\operatorname{vol}_G(S^c)\}}.
\]
Write $\tau(G)$ for the minimum number of monochromatic edges in a
two-coloring of $V(G)$.

\begin{lem}[Conductance bound]
\label{lem:PA-conductance}
There is a constant $c=c(m)>0$ such that
\[
    \sup_{\alpha\in[0,1)}\mathbb P_\alpha(h(G_n)<c)=o(1).
\]
\end{lem}

\begin{proof}
Choose $\eta=\eta(m)>0$ sufficiently small, with $\eta\le1/100$.
We first show that, uniformly in $\alpha$, with probability $1-o(1)$,
\begin{equation}
    b_{G_n}(S)\ge\eta|S|
    \qquad(1\le|S|\le n/2).
    \label{eq:PA-no-sparse-cuts}
\end{equation}
Connectivity gives this inequality when $|S|\le1/\eta$.
For a fixed set $S$ of size $1/\eta<k\le n/2$ and an integer
$1\le a\le\eta k$, we claim that
\begin{equation}
    \mathbb P_\alpha(b_{G_n}(S)=a)
       \le (mn)^m\binom{mn}{a}
          \frac{\binom{mk}{a}}{\binom{mn-a}{mk-a}}.
    \label{eq:fixed-cut-input}
\end{equation}

We adapt the crossing-pattern argument of
\cite[Lemma~3.2 and Appendix~A]{alp-nonzero-PA} to keep the bound
uniform in $\alpha$. Treat the seed edges as created by vertex $2$,
orient them from $2$ to $1$, and label all $m(n-1)$ edges in creation
order. Fix a set of $a$ labels prescribed
to cross $(S,S^c)$; incompatible seed prescriptions have probability
zero.

Consider a nonseed edge prescribed not to cross, created by a vertex
$r$ in a side $C\in\{S,S^c\}$. Let $i$ be its rank among the
noncrossing edges created by vertices in $C$, let $z,y$ count the preceding
noncrossing and crossing edges, and let $h$ count the preceding
edges from $r$. Let $h_C\le h$ count those preceding edges from $r$
whose target lies in $C$. Among the $i-1$ preceding noncrossing
edges created by vertices in $C$, the $i-1-h_C$ created by vertices other
than $r$ contribute $2$ each to the degree mass in $C$ (both endpoints
lie in $C$), while the $h_C$ created by $r$ contribute only their
target's degree, $1$ each;
the $y-(h-h_C)$ preceding crossing edges created by vertices other than $r$ contribute
$1$ each, from whichever endpoint lies in $C$. Summing and cancelling
$h_C$: on a history consistent with the prescription, the eligible
degree mass in $C$ is
$2(i-1-h_C)+h_C+[y-(h-h_C)]=2(i-1)+y-h$, and the same count without
restricting to $C$ gives $2(z+y)-h$. Also,
$m|C\cap[r-1]|\le i-1+y-h+m$ and
$m(r-1)=z+y-h+m$. Hence the probabilities of staying in $C$ under
degree-proportional and uniform attachment are bounded, respectively,
by
\[
    \frac{2(i-1)+y-h}{2(z+y)-h}
       \le\frac{i-1+y}{z+y},
    \qquad
    \frac{i-1+y+m-h}{z+y+m-h}.
\]
Since $z\ge i-1$ and $y\le a$, both are at most
$(i-1+a+m)/(z+a+m)$. Writing $D_{k,i}:=2m(k-2)+(i-1)$ for the
sum $\sum_{v<k}d_{k,i-1}(v)$, rule
\eqref{eq:PA-rule-intro} equals
$\beta_{k,i}\,d_{k,i-1}(v)/D_{k,i}+(1-\beta_{k,i})/(k-1)$, with
$\beta_{k,i}=(1-\alpha)D_{k,i}/[2m\alpha(k-1)+(1-\alpha)D_{k,i}]\in[0,1]$
depending on $(k,i)$ alone: a convex combination, with
history-independent weight, of the degree-proportional and uniform
rules above. Its conditional noncrossing probability is therefore at
most the same convex combination of the two bounds above, hence at
most
\begin{equation}
    \frac{i-1+a+m}{z+a+m}
       \le\frac{i+a+m}{z+a+m+1}.
    \label{eq:fixed-pattern-one-step}
\end{equation}

Put $M=m|S\setminus\{1\}|$, and let $s$ be the number of prescribed
crossing edges created by vertices in $S$. Bound each crossing probability by one
and multiply the right-hand bound of
\eqref{eq:fixed-pattern-one-step} over the remaining labels. The
noncrossing seed labels can be included: their ranks satisfy
$i=z+1$, so their factors equal one. The global noncrossing
ranks run from $1$ to $m(n-1)-a$, and the ranks on the two sides run
from $1$ to $M-s$ and from $1$ to $m(n-1)-M-a+s$. Thus the probability of the
prescribed pattern is at most
\begin{align*}
    \frac{(M+a-s+m)!(mn-M+s)!}{(mn)!(a+m)!}
    \le\frac{(M+m)!(mn-M)!}{(mn-a)!(a+m)!}
    \le(mn)^m\frac{(mk)![m(n-k)]!}{(mn-a)!a!}.
\end{align*}
For the first inequality, pair the $a$ extra numerator factors with
the last $a$ factors of $(mn)!$, using $a\le M$ and $M+m+a\le mn$.
For the second, use $M\in\{mk-m,mk\}$; at most $m$ extra numerator
factors remain, each bounded by $mn$. There are at most
$\binom{mn}{a}$ prescribed patterns, proving
\eqref{eq:fixed-cut-input}.

Vandermonde's identity gives
\[
    \binom{mn-a}{mk-a}
       \ge\binom nk\binom{(m-1)n-a}{(m-1)k-a}.
\]
The usual binomial bounds therefore imply, for $1\le a\le\eta k$,
\begin{align*}
    \binom nk\binom{mn}{a}
       \frac{\binom{mk}{a}}{\binom{mn-a}{mk-a}}
    \le\left[
       \left(\frac{\mathrm e mk}{a}\right)^{2a/k}
       \left(\frac kn\right)^{m-1-2a/k}
       \right]^k
    \le\left(\frac kn\right)^{(m-1)k/2}.
\end{align*}
The last inequality holds for sufficiently small $\eta$, uniformly
in $k/n\le1/2$. A union bound over $S$ and $a$ now gives
\[
    \sup_{\alpha\in[0,1)}\mathbb P_\alpha\left(
       \exists S:\ 1\le|S|\le n/2,\ b_{G_n}(S)<\eta|S|
    \right)
       \le C_m n^m\sum_{1/\eta<k\le n/2}
          k\left(\frac kn\right)^{(m-1)k/2}=o(1).
\]
Indeed, split the sum at $\sqrt n$. Below this point use
$k/n\le n^{-1/2}$; above it use $k/n\le1/2$. The resulting
geometric bounds tend to zero after multiplication by $n^m$ when
$\eta$ is sufficiently small.
This proves \eqref{eq:PA-no-sparse-cuts}.

Every vertex has at most $m$ edges directed toward older vertices,
so $e_{G_n}(U)\le m|U|$. For any nontrivial cut, take its side $U$
of smaller cardinality. On \eqref{eq:PA-no-sparse-cuts},
\[
    \frac{b_{G_n}(U)}{
       \min\{\operatorname{vol}_{G_n}(U),
              \operatorname{vol}_{G_n}(U^c)\}}
       \ge\frac{b_{G_n}(U)}{2m|U|+b_{G_n}(U)}
       \ge\frac{\eta}{2m+\eta}.
\]
Taking the minimum over cuts proves the lemma.
\end{proof}

\begin{lem}[Distance from bipartiteness]
\label{lem:PA-MaxCut-defect}
There is a constant $c>0$ such that, for all large enough $n$,
\[
    \sup_{\alpha\in[0,1)}\mathbb P_\alpha
       \left(\tau(G_n)\le\frac n{200}\right)\le e^{-cn}.
\]
\end{lem}

\begin{proof}
Every older vertex has degree at least $m$, so
\eqref{eq:PA-rule-intro} gives
\begin{equation}
    \mathbb P_\alpha(J_{k,i}=v\mid\mathcal F_{k,i-1})
       \ge\frac1{2(k-1)},\qquad v<k.
    \label{eq:spectral-attachment-lower-bound}
\end{equation}
One seed edge and the first edge from each vertex $k\ge3$ form
a spanning tree. Fix $\sigma_1=+1$. Every two-coloring
is then encoded by a binary sequence $(b_2,\ldots,b_n)$, where
$b_k=1$ precisely when the parent edge of $k$ is monochromatic.
A coloring with at most $n/200$ monochromatic edges has at most
$n/200$ ones in its code.

Fix such a code and expose the graph sequentially. Once the first
edge from $k$ reveals its parent, assign the color $\sigma_k$ prescribed
by $b_k$. Let $X_k$ indicate that the second edge from $k$ is
monochromatic, and put
\[
    r_k=\frac{|\{v<k:\sigma_v=\sigma_k\}|}{k-1}.
\]
The colors in this expression are known before the second edge is
exposed. Thus \eqref{eq:spectral-attachment-lower-bound} implies
$\mathbb P_\alpha(X_k=1\mid\mathcal F_{k,1})\ge r_k/2$.
Every two-coloring of $n$ vertices has at least $n(n-2)/4$
monochromatic vertex pairs. Since each denominator in $r_k$ is at
most $n-1$, omitting the pair $\{1,2\}$ gives, pathwise,
\[
    \sum_{k=3}^n r_k\ge\frac n4-O(1).
\]
The conditional Bernoulli exponential bound, iterated in the
exposure filtration, gives
\[
    \mathbb E_\alpha\exp\left\{
       -\sum_{k=3}^nX_k
       +\frac{1-e^{-1}}2\sum_{k=3}^n r_k
    \right\}\le1.
\]
Consequently, for this fixed code,
\[
    \mathbb P_\alpha\left(\sum_{k=3}^nX_k\le\frac n{200}\right)
       \le\exp\left\{\frac n{200}
                  -\frac{1-e^{-1}}8n+O(1)\right\}.
\]

If $\tau(G_n)\le n/200$, some code with at most $n/200$ ones
satisfies the event in the last display. For large $n$, the number
of these codes is at most
\[
    \sum_{0\le j\le n/200}\binom{n-1}{j}\le e^{n/25}.
\]
Since $1/200+1/25<(1-e^{-1})/8$, the union bound proves the result.
\end{proof}

\begin{proof}[Proof of Proposition~\ref{prop:PA-uniform-two-sided-gap}]
The usual Cheeger inequality \cite{ChungSpectralGraphTheory} gives
$1-\lambda_2(P_n)\ge h(G_n)^2/2$.
To control the other end of the spectrum, write $G=G_n$. We show
that for disjoint sets $L,R$ with $S=L\cup R\ne\varnothing$,
\begin{equation}
    \frac{2e_G(L)+2e_G(R)+b_G(S)}{\operatorname{vol}_G(S)}
       \ge h(G)\frac{\tau(G)}{|E(G)|}.
    \label{eq:spectral-cut-comparison}
\end{equation}
If $\operatorname{vol}_G(S)\le|E(G)|$, the left side is at least
$h(G)$, which suffices since $\tau(G)\le|E(G)|$.

Otherwise put $a=e_G(L)+e_G(R)$ and $b=b_G(S)$.
Extending the coloring of $S$ by independent fair colors on $S^c$
shows that
\[
    \tau(G)\le a+\frac{\operatorname{vol}_G(S^c)+b}{4},
    \qquad h(G)\operatorname{vol}_G(S^c)\le b.
\]
Since $h(G)\le1$, these inequalities give
$2h(G)\tau(G)\le2a+b$. Dividing by
$\operatorname{vol}_G(S)\le2|E(G)|$ proves
\eqref{eq:spectral-cut-comparison}.

The minimum of the left-hand side of
\eqref{eq:spectral-cut-comparison} over such pairs $L,R$ is
$1-\bar h(G)$, where $\bar h(G)$ is the dual Cheeger constant.
The dual Cheeger inequality
\cite[Theorem~3.2]{BauerJostDual}, applied with edge multiplicities
as weights, therefore yields
\[
    1+\lambda_n(P_n)
       \ge\frac{(1-\bar h(G_n))^2}{2}
       \ge\frac{h(G_n)^2}{2}
          \left(\frac{\tau(G_n)}{m(n-1)}\right)^2.
\]
Lemmas~\ref{lem:PA-conductance} and
\ref{lem:PA-MaxCut-defect} bound both $h(G_n)$ and
$\tau(G_n)/[m(n-1)]$ away from zero with probability $1-o(1)$,
uniformly in $\alpha$. The two spectral bounds prove the proposition.
\end{proof}

%=================================================================================================================
\medskip

\section{Random Walk on the P\'olya--Point Tree}\label{sec:RegRWPA}

Throughout this section, \(m\geq2\) and \(\alpha\in[0,1)\).
Let \(\Pbb\) denote the annealed law of the PPT of
Definition~\ref{def:PPT} and simple random walk \(X\) started at its
root \(\rho\). Write \(|v|=\dist_{\Tcal}(\rho,v)\); after re-rooting,
depth is measured from the new root. We prove concentration of the
distance and of harmonic information, and derive the estimates used
in the logarithmic-time exploration.

\subsection{Escape probabilities and harmonic measure}
\label{subsec:escape-harmonic}

For a fixed realization of \(\Tcal\), write \(\Pbb_x^{\Tcal}\) for
simple random walk started from \(x\), and put
\(\tau_z:=\inf\{n\geq0:X_n=z\}\). For a child \(y\) of \(z\), write
\(y^-=z\) and let \(h_z(y):=\Pbb_y^{\Tcal}(\tau_z=\infty)\), abbreviated
to \(h(y)\). We also use descendant trees with an artificial parent
attached to their root. The root degree includes this additional edge,
and no environment beyond the artificial parent is needed before it is
hit.

\begin{lem}[Uniform escape]\label{lem:escape}
Almost surely, for every parent--child edge \(z\to y\), including an
artificial parent edge,
\[
    h_z(y)\geq p_{\rm esc}:=\frac13.
\]
In particular, simple random walk on \(\Tcal\) is transient.
\end{lem}

\begin{proof}
Choose a left child \(w\) of \(y\). Since \(w\) has type \(\tL\),
recursively retaining two left children of each retained vertex gives an
infinite binary tree. Its effective resistance from \(w\) to infinity
is one. Hence the effective conductance from \(y\) to infinity within
this retained subtree, including the edge \(yw\), is \(1/2\). Rayleigh
monotonicity shows that the effective conductance \(C_y\) from \(y\)
to infinity through all its descendants satisfies \(C_y\geq1/2\).
Since the edge \(yz\) has unit conductance,
\[
    h_z(y)=\frac{C_y}{1+C_y}\geq\frac13.
\]
The argument is unchanged when \(z\) is an artificial parent. The
positive effective conductance to infinity also proves transience.
\end{proof}

Set \(\varrho:=1-p_{\rm esc}=2/3\). By the strong Markov property,
if \(x\) is a descendant of \(z\), then
\begin{equation}\label{eq:ancestor-return}
    \Pbb_x^{\Tcal}(\tau_z<\infty)
    \leq \varrho^{\dist(x,z)}.
\end{equation}

By transience and local finiteness, the walk converges to an end. Let
\(\xi=(\rho,\xi_1,\ldots)\) be the corresponding terminal ray. For a
child \(y\) of \(z\), let \(\theta_z(y)\) be the conditional
probability that this ray follows \(z\to y\), given that it contains
\(z\), and write \(\theta(y):=\theta_{y^-}(y)\).

\begin{lem}[Harmonic kernel]\label{lem:harmonic-kernel}
For every child \(y\) of \(z\),
\begin{equation}\label{eq:harmonic-kernel}
    \theta_z(y)=\frac{h_z(y)}{\sum_{w:w^-=z}h_z(w)}.
\end{equation}
Consequently,
\begin{equation}\label{eq:log-theta-degree}
    -\log\theta_z(y)
    \leq\log\deg_{\Tcal}(z)+\log p_{\rm esc}^{-1}
    =\log\deg_{\Tcal}(z)+\log3.
\end{equation}
\end{lem}

\begin{proof}
At each visit to \(z\), the probability of stepping to \(y\) and never
returning is \(h_z(y)/\deg_{\Tcal}(z)\). Summing over visits gives
the probability that the terminal ray follows \(z\to y\). The
expected number of visits is a common finite factor, so normalizing
over the children proves \eqref{eq:harmonic-kernel}. The bound follows
from \(h_z(y)\geq p_{\rm esc}\) and \(h_z(w)\leq1\).
\end{proof}

The harmonic information along the first \(k\) levels of the ray is
\[
    W_0:=0,\qquad
    W_k:=\sum_{j=1}^k-\log\theta(\xi_j),\qquad k\geq1.
\]
The same notation will be used after re-rooting at a regeneration
vertex, with the ray starting from that vertex.

We shall condition the walk to avoid an ancestor \(b\). For \(x\)
strictly below \(b\), put
\(H_b(x):=\Pbb_x^{\Tcal}(\tau_b=\infty)\) and
\(\Pbb_x^{\Tcal,b}:=\Pbb_x^{\Tcal}(\,\cdot\mid\tau_b=\infty)\).
Set \(H_b(b)=0\); when no ancestor is forbidden, write \(b=\partial\)
and \(H_\partial\equiv1\). Avoiding the parent of \(x\) prevents a
visit to \(b\), so \(H_b(x)\geq p_{\rm esc}\). The Markov property
gives the conditioned transition probabilities
\[
    p^b(x,y)=\frac{H_b(y)}{\deg_{\Tcal}(x)H_b(x)},\qquad y\sim x.
\]
For a child \(y\) of \(x\), every path to \(b\) passes through \(x\),
and hence
\begin{equation}\label{eq:H-monotone}
    H_b(y)=H_b(x)+(1-H_b(x))h_x(y)\geq H_b(x).
\end{equation}
Thus conditioning does not decrease the transition probability to any
prescribed child. Moreover,
\begin{equation}\label{eq:conditioned-nonreturn}
    \Pbb_y^{\Tcal,b}(\tau_x=\infty)
    =\frac{h_x(y)}{H_b(y)}\geq p_{\rm esc}.
\end{equation}

Conditioning the walk from a descendant root \(o\) never to hit its
artificial parent \(b\) preserves the harmonic kernel
\eqref{eq:harmonic-kernel}. Indeed, for each child \(y\) of \(x\),
\[
    p^b(x,y)\Pbb_y^{\Tcal,b}(\tau_x=\infty)
    =\frac{h_x(y)}{\deg_{\Tcal}(x)H_b(x)}.
\]
The conditioned walk is transient and stays in the descendant tree of
\(o\), so the same last-departure argument, normalized over the children
of \(x\), applies.

\subsection{Regeneration}
\label{subsec:regeneration}

We explore the tree along the walk. On the first visit to a vertex,
we reveal its mark and all its children, including their ages and types.
Let \(\Hcal_n\) denote the exploration history immediately
before any new data are revealed at time \(n\). Write \(\Tcal_v\) for
the complete marked descendant tree rooted at \(v\).

\begin{defn}[Regeneration times]\label{def:regeneration-times}
A step \(X_{n-1}=u\to X_n=y\) is a \emph{candidate} if \(u\) and \(y\)
are first visited at times \(n-1\) and \(n\), respectively, \(y\) is a
child of \(u\), and \(A_u<1/2<A_y\). It is \emph{accepted} if the walk
never subsequently returns to \(u\). Put \(r_0=0\), and let
\(r_1<r_2<\cdots\) be the accepted candidate times, with value
\(\infty\) if the next accepted candidate does not exist. When finite,
write \(v_i=X_{r_i}\) and \(b_i=X_{r_i-1}\).
\end{defn}

The age law at a candidate entrance is
\begin{equation}\label{eq:mu-R}
    \mu_R(\dd s)
    :=\frac{(1-\chi)s^{-\chi}}{1-2^{-(1-\chi)}}
      \mathbf 1_{(1/2,1)}(s)\,\dd s.
\end{equation}
Indeed, on the first visit to a vertex of age less than \(1/2\), the right-child
intensity on \((1/2,1)\) is a random constant times
\(s^{-\chi}\,\dd s\). Conditional on the number of children, the walk
chooses uniformly among the neighbors, so the degree bias does not
change this normalized age density. The selected child has type
\(\tR\), and, conditional on its age, its own mark and descendant tree
have the prescribed type-\(\tR\) PPT law.

\begin{defn}[Post-regeneration law]\label{def:Pstar}
Let \(\Pbb_{\mu_R}\) be the annealed law of a type-\(\tR\) PPT
rooted at \(o\), with \(A_o\sim\mu_R\), an artificial parent \(b_o\)
attached to \(o\), and simple random walk started from \(o\). Define
\[
    \Pbb_\star:=\Pbb_{\mu_R}(\,\cdot\mid\tau_{b_o}=\infty).
\]
\end{defn}

Lemma~\ref{lem:escape} gives
\(\Pbb_{\mu_R}(\tau_{b_o}=\infty)\geq p_{\rm esc}\).
Conditioning changes both the tree and the walk: the tree marginal is
biased by \(H_{b_o}(o)\), and, given the tree, the walk has law
\(\Pbb_o^{\Tcal,b_o}\). For every nonnegative functional \(Y\) of the
tree and walk,
\begin{equation}\label{eq:transfer}
\Ebb_\star Y\leq p_{\rm esc}^{-1}\Ebb_{\mu_R}Y.    
\end{equation}

\begin{lem}[Law at an accepted candidate]\label{lem:accepted-child}
Let \(\sigma\) be a finite exploration stopping time at which
\(u=X_\sigma\) is visited for the first time and \(A_u<1/2\).
The history \(\Hcal_\sigma\)
records the past walk and all information revealed through this
entrance, before \(\Gamma_u\) and the children of \(u\) are exposed.
Put \(Y=X_{\sigma+1}\).

Under either \(\Pbb\) or \(\Pbb_\star\), conditional on
\(\Hcal_\sigma\) and on \(u\to Y\) being an accepted candidate,
the marked descendant tree \(\Tcal_Y\), rooted at \(Y\), together
with the walk \((X_{\sigma+1+k})_{k\geq0}\), has law \(\Pbb_\star\).
Here \(u\) is replaced by an artificial parent carrying no age or mark.
\end{lem}

\begin{proof}
First work under the ordinary annealed law, either \(\Pbb\) or
\(\Pbb_{\mu_R}\). Conditional on \(\Hcal_\sigma\), expose \(\Gamma_u\)
and the numbers of children of \(u\) with ages below and above \(1/2\).
By the Poisson offspring rule and \eqref{eq:mu-R}, the descendant
trees rooted at children of age above \(1/2\), in an independent
random ordering, are i.i.d.\ with the tree marginal of
\(\Pbb_{\mu_R}\). The degree of \(u\) is fixed by the conditioning
data, and the walk chooses uniformly among its neighbors.
Conditional on the next step being a candidate, the selected tree
therefore has this same law. Since the law does not depend on the
exposed data, this remains true conditional only on \(\Hcal_\sigma\)
and the candidate event.

By the Markov property at \(\sigma+1\), the selected tree together
with the walk from \(Y\), stopped on hitting \(u\), has law
\(\Pbb_{\mu_R}\) stopped on hitting its artificial parent.
Acceptance is precisely the event that this parent is never hit.
Thus, writing \(B\) for the event that \(u\to Y\) is an accepted
candidate, for every bounded measurable functional \(F\) of the
selected tree and shifted walk,
\[
    \Ebb[F\mid\Hcal_\sigma,B]
    =
    \frac{\Ebb_{\mu_R}[F;\tau_{b_o}=\infty]}
         {\Pbb_{\mu_R}(\tau_{b_o}=\infty)}
    =\Ebb_\star[F].
\]

Finally, under \(\Pbb_\star\), the past has avoided the initial
artificial parent \(b_o\), and on \(B\) the future remains in
\(\Tcal_Y\). Hence, conditional on \(\Hcal_\sigma\) and \(B\), the
avoidance event defining \(\Pbb_\star\) is automatic. The conditional
law therefore agrees with that under \(\Pbb_{\mu_R}\), to which the
preceding argument applies.
\end{proof}

\subsection{Concentration of distance}
\label{subsec:distance}

\begin{prop}[Exponential regeneration tails]\label{prop:exp-moment-R}
There exists \(\varepsilon>0\) such that
\[
    \Ebb e^{\varepsilon r_1}+\Ebb_\star e^{\varepsilon r_1}<\infty.
\]
\end{prop}
\begin{proof}
We first give uniform bounds for a trial started at the first visit to
a vertex \(v\). If \(A_v>1/2\), a
specified left child has age below \(1/2\) with probability at least
\(2^{-\chi}\). The event that \(\Gamma_v\leq1\), there are no right
children, and this left child has age below \(1/2\) has uniformly
positive probability. On that event the degree is at most \(m+1\).
Thus the next step is to an unvisited child of age below \(1/2\) with a
uniformly positive probability.

If \(A_v<1/2\), restrict to \(\Gamma_v\geq1\). The number of right
children then dominates \(N\sim\Poi(2^{1-\chi}-1)\); each right child
has conditional probability at least \(1-2^{-(1-\chi)}\) of having
age above \(1/2\). Since the degree is at most \(m+1\) plus the
number of right children, the probability of choosing such a child is
at least
\[
    \Pbb(\Gamma_v\geq1)
    \bigl(1-2^{-(1-\chi)}\bigr)
    \Ebb\!\left[\frac{N}{m+1+N}\right]>0,
\]
uniformly in the age and type, including the original root. These two
bounds give a candidate within two steps of the first visit to any vertex, with
probability bounded below by a constant. To transfer this bound to
\(\Pbb_\star\), let \(\sigma\) be a stopping time at which \(X_\sigma\)
is visited for the first time,
and let \(E\) be the corresponding descendant trial event of length
\(\ell\leq2\).
On histories avoiding \(b_o\), conditioning on the tree and the walk
through \(\sigma+\ell\) gives
\[
    \Pbb_\star(E\mid\Hcal_\sigma)
    =\frac{\Ebb_{\mu_R}[\mathbf1_E H_{b_o}(X_{\sigma+\ell})
             \mid\Hcal_\sigma]}
           {\Ebb_{\mu_R}[H_{b_o}(X_\sigma)\mid\Hcal_\sigma]}
    \geq p_{\rm esc}\Pbb_{\mu_R}(E\mid\Hcal_\sigma).
\]
Indeed, on \(E\) the numerator's avoidance factor is at least
\(p_{\rm esc}\), and the denominator is at most one.
At the candidate entrance, the conditional probability of acceptance
is at least \(p_{\rm esc}\) by
\eqref{eq:conditioned-nonreturn}.

We next control the time lost between trials. Put
\[
    f(v):=|v|+\frac12\mathbf 1_{\{\typ(v)=\tL\}},
\]
where the indicator is interpreted as zero at the untyped original
root. In the following drift and hitting-time estimates, the tree is
fixed and all expectations and probabilities are quenched. If the
current vertex \(v=X_n\) has type \(\tL\), then under simple random walk
\[
 \Ebb[f(X_{n+1})-f(X_n)\mid X_n]
 \geq \frac{m-\frac32+\frac12N_v}{m+N_v+1}\geq\frac16.
\]
Indeed, a left-child step has increment \(1\), a right-child step has
increment \(1/2\), and the parent step has increment at least
\(-3/2\). If \(v\) has type \(\tR\), the corresponding
bound is
\[
 \Ebb[f(X_{n+1})-f(X_n)\mid X_n]
 \geq \frac{\frac32(m-1)+N_v-1}{m+N_v}\geq\frac14.
\]
At the original root every increment is positive. The same holds at
the post-regeneration root under the conditioned law, since its
artificial parent cannot be visited. Elsewhere under the conditioned
law, \eqref{eq:H-monotone} shows that the Doob transform transfers
transition mass from the parent to the children. Since every child
step has positive \(f\)-increment and the parent step has negative
increment, the same lower bounds remain valid. As
\(|f(X_{n+1})-f(X_n)|\leq3/2\), there are constants
\(\lambda_0,c_0>0\), uniform in the tree and the past, such that under
either quenched law
\begin{equation}\label{eq:radial-exponential-drift}
 \Ebb\!\left[
   e^{-\lambda_0(f(X_{n+1})-f(X_n))}
   \mid X_0,\ldots,X_n
 \right]\leq e^{-c_0}.
\end{equation}

Suppose that the walk starts from a child \(v\) of \(u\). Before
\(\tau_u\), the walk stays in the descendant tree of \(u\), where
\(f-f(u)>0\). Applying \eqref{eq:radial-exponential-drift} up to
\(\tau_u\) gives, for the walk started from any strict descendant
\(w\) of \(u\),
\[
 \Pbb_w(\tau_u<\infty)
 \leq e^{-\lambda_0(f(w)-f(u))}
\]
and
\[
 \Ebb_v\!\left[
   e^{-\lambda_0(f(X_n)-f(u))}
   \mathbf 1_{\{\tau_u>n\}}
 \right]
 \leq e^{-c_0n}e^{-\lambda_0(f(v)-f(u))}.
\]
The Markov property at time \(n\), with the first bound applied at
\(w=X_n\), therefore yields, uniformly in the tree and the past,
\begin{equation}\label{eq:late-return-tail}
    \Pbb_v(n<\tau_u<\infty)\leq e^{-c_0n}.
\end{equation}

If \(\sigma_k\) is the time to reach a level \(k\) above the current
height, then on \(\{\sigma_k>n\}\) one has
\(f(X_n)-f(X_0)\leq k\). Applying
\eqref{eq:radial-exponential-drift} before \(\sigma_k\) gives
\[
    \Pbb(\sigma_k>n)\leq\min\{1,e^{\lambda_0k-c_0n}\}.
\]
Consequently, for \(0<\varepsilon<c_0\), summation of the tail gives
\begin{equation}\label{eq:radial-hitting-mgf}
    \Ebb e^{\varepsilon\sigma_k}
    \leq C_\varepsilon e^{c_\varepsilon k},\qquad k\geq1,
    \qquad c_\varepsilon:=\frac{\lambda_0}{c_0}\varepsilon
    \downarrow0
    \quad\text{as }\varepsilon\downarrow0.
\end{equation}
All these bounds hold conditionally at stopping times.

Now test the one- or two-step age pattern at height records. If
the finite test fails to produce a candidate, resume at the first level
strictly above the maximum visited during the test. The required height
increase is at most five, so \eqref{eq:radial-hitting-mgf} gives a
uniform exponential moment for the time until the next test. The
conditional probability of a failed test is bounded away from one.
Consequently, for some small \(\varepsilon>0\), its duration \(D\)
satisfies \(\Ebb[e^{\varepsilon D}\mathbf 1_{\{\text{failure}\}}
\mid\text{past}]<1\), uniformly. Iterating this bound proves that the
time \(S_0\) of the first tested candidate has a uniform exponential
moment.

Let \(S_j\) be a tested candidate time, and let \(B_j\) be its
acceptance event. On \(B_j^c\), wait for the actual return to the old
endpoint \(u\). Its duration \(\tau_u\), measured from \(S_j\),
satisfies the defective exponential tail
\eqref{eq:late-return-tail}. Then wait to reach a level above the
maximum visited during this rejected excursion, and resume the same
candidate search. The required height increase is at most \(\tau_u+2\).
By \eqref{eq:radial-hitting-mgf}, the conditional exponential moment of
this catch-up time is at most
\(C_\varepsilon e^{c_\varepsilon(\tau_u+2)}\). The subsequent candidate
search has the uniform exponential moment just proved. Decreasing
\(\varepsilon\) so that \(\varepsilon+c_\varepsilon<c_0\) in
\eqref{eq:late-return-tail} therefore yields
\[
    \sup_j\Ebb\!\left[
       e^{\varepsilon(S_{j+1}-S_j)}\mathbf 1_{B_j^c}
       \,\middle|\,\Hcal_{S_j}\right]<\infty.
\]
At exponent zero these expectations are at most \(1-p_{\rm esc}\).
After decreasing \(\varepsilon\) further, uniform exponential
integrability makes them all at most a constant \(\kappa<1\). Iterating at
the successive candidate stopping times and summing over the number of
rejections gives an exponential moment for the first accepted time of
this restricted search. Since \(r_1\) is no later, the proposition
follows under both laws.
\end{proof}

Under \(\Pbb_\star\), set \(v_0=o\), \(b_0=b_o\), and define
\[
    R_i:=r_{i+1}-r_i,\qquad
    T_i:=\dist_{\Tcal}(v_i,v_{i+1}),\qquad i\geq0.
\]
Under the original law, use the same definitions for \(i\geq1\).
The regeneration vertices are successive descendants, so
\(1\leq T_i\leq R_i\). The \(i\)-th \emph{regeneration slab}
\(\mathcal C_i\) consists of the marked tree
\(\Tcal_{v_i}\setminus\Tcal_{v_{i+1}}\), rooted at \(v_i\),
together with the walk segment \((X_{r_i+t})_{0\leq t<R_i}\).
It includes every side branch attached to the retained tree.
We retain an unmarked incoming half-edge at \(v_i\) and a
distinguished unmarked outgoing half-edge at \(b_{i+1}\), recording
that the next step crosses the latter. The vertex \(v_{i+1}\), its
age and marks, and all its descendants are excluded. We forget
the original vertex addresses and child order, while retaining the
genealogical types. In particular, the outgoing half-edge carries no
age rank. Time is measured from
\(r_i\). Thus \(R_i\) and \(T_i\) are
determined by \(\mathcal C_i\).

Under the original law, the initial block \(\mathcal I_0\) consists of
\(\Tcal_\rho\setminus\Tcal_{v_1}\), rooted at \(\rho\), the path
\((X_t)_{0\leq t<r_1}\), and the unmarked outgoing half-edge at
\(b_1\), with the same convention on labels.

\begin{figure}[htbp]
    \centering
    \begin{tikzpicture}[x=0.6cm, y=0.6cm, thick, >=stealth]
        
        % -----------------------------
        % Coordinates
        % -----------------------------
        \coordinate (rho) at (0,-0.2);
        
        % Block 1 (Initial tree block before r_{i-1})
        \coordinate (b1_left) at (-1.5,-2.5);
        \coordinate (b1_right) at (1.5,-2.5);
        \coordinate (tau_im1_minus1) at (-0.5,-2.5); % Hollow dot
        
        % Regeneration node i-1
        \coordinate (tau_im1) at (-0.5,-3.5);       % Solid dot
        
        % Half-open slab C_{i-1}
        \coordinate (b2_left) at (-2,-6.5);
        \coordinate (b2_right) at (1,-6.5);
        \coordinate (tau_i_minus1) at (0,-6.5);     % Hollow dot
        
        % Regeneration node i
        \coordinate (tau_i) at (0,-7.5);            % Solid dot
        
        % Half-open slab C_i
        \coordinate (b3_left) at (-1,-10.5);
        \coordinate (b3_right) at (2,-10.5);
        \coordinate (tau_ip1_minus1) at (0.5,-10.5); % Hollow dot
        
        % Regeneration node i+1
        \coordinate (tau_ip1) at (0.5,-11.5);        % Solid dot
        
        % -----------------------------
        % Draw Subtrees (Triangles)
        % -----------------------------
        \draw[fill=gray!5] (rho) -- (b1_left) -- (b1_right) -- cycle;
        \draw[fill=gray!5] (tau_im1) -- (b2_left) -- (b2_right) -- cycle;
        \draw[fill=gray!5] (tau_i) -- (b3_left) -- (b3_right) -- cycle;
        
        % Labels inside subtrees
        \node at (-0.95, -5.85) {$\mathcal C_{i-1}$};
        \node at (0.3, -9.4) {$\mathcal C_i$};
        \node[above=4pt] at (rho) {\small $\rho$ (root)};
        
        % -----------------------------
        % Draw Regeneration Edges
        % -----------------------------
        \draw[ultra thick] (tau_im1_minus1) -- (tau_im1);
        \draw[ultra thick] (tau_i_minus1) -- (tau_i);
        \draw[ultra thick] (tau_ip1_minus1) -- (tau_ip1);
        
        % -----------------------------
        % Draw and Label Nodes
        % -----------------------------
        % Edge i-1
        \filldraw[fill=white, draw=black] (tau_im1_minus1) circle (2.5pt) node[left=24pt] {\small$X_{r_{i-1}-1}$};
        \filldraw[fill=black, draw=black] (tau_im1) circle (2.5pt) node[left=4pt] {\small $X_{r_{i-1}}$};
        
        % Edge i
        \filldraw[fill=white, draw=black] (tau_i_minus1) circle (2.5pt) node[left=40pt] {\small $X_{r_i-1}$};
        \filldraw[fill=black, draw=black] (tau_i) circle (2.5pt) node[left=1pt] {\small $X_{r_i}$};
        
        % Edge i+1
        \filldraw[fill=white, draw=black] (tau_ip1_minus1) circle (2.5pt); 
        \filldraw[fill=black, draw=black] (tau_ip1) circle (2.5pt) node[left=1pt] {\small$X_{r_{i+1}}$};
        
        % -----------------------------
        % Depth Arrows and Distance Markers
        % -----------------------------
        % Dashed projection lines
        \draw[dashed, gray] (tau_im1) -- (3.5, -3.5);
        \draw[dashed, gray] (tau_i) -- (5, -7.5);
        \draw[dashed, gray] (rho) -- (5, -0.2);
        
        % Arrows
        \draw[<->] (3.5, -0.2) -- (3.5, -3.5) node[midway, fill=white] {\small $|X_{r_{i-1}}|$};
        \draw[<->] (5, -0.2) -- (5, -7.5) node[midway, fill=white] {\small $|X_{r_i}|$};
        
        % LERW path notation
        \node[right] at (5.2, -7.5) {\small $X_{r_i} = \xi_{|X_{r_i}|}$};
        
        % -----------------------------
        % Optional: Indicate the LERW \xi passing through the spine
        % -----------------------------
        \draw[->, dotted, very thick, blue] (0, -0.5) .. controls (-0.5, -1.5) and (-0.2, -2) .. (tau_im1_minus1);
        \draw[->, dotted, very thick, blue] (tau_im1) .. controls (-0.8, -5) and (-0.2, -6) .. (tau_i_minus1);
        
    \end{tikzpicture}
    \caption{The half-open regeneration slabs separated by the accepted
    edges $(X_{r_i-1},X_{r_i})$.  The slabs are i.i.d.; the forward
    descendant trees rooted at $X_{r_i}$ are nested and are not
    independent.}
\end{figure}
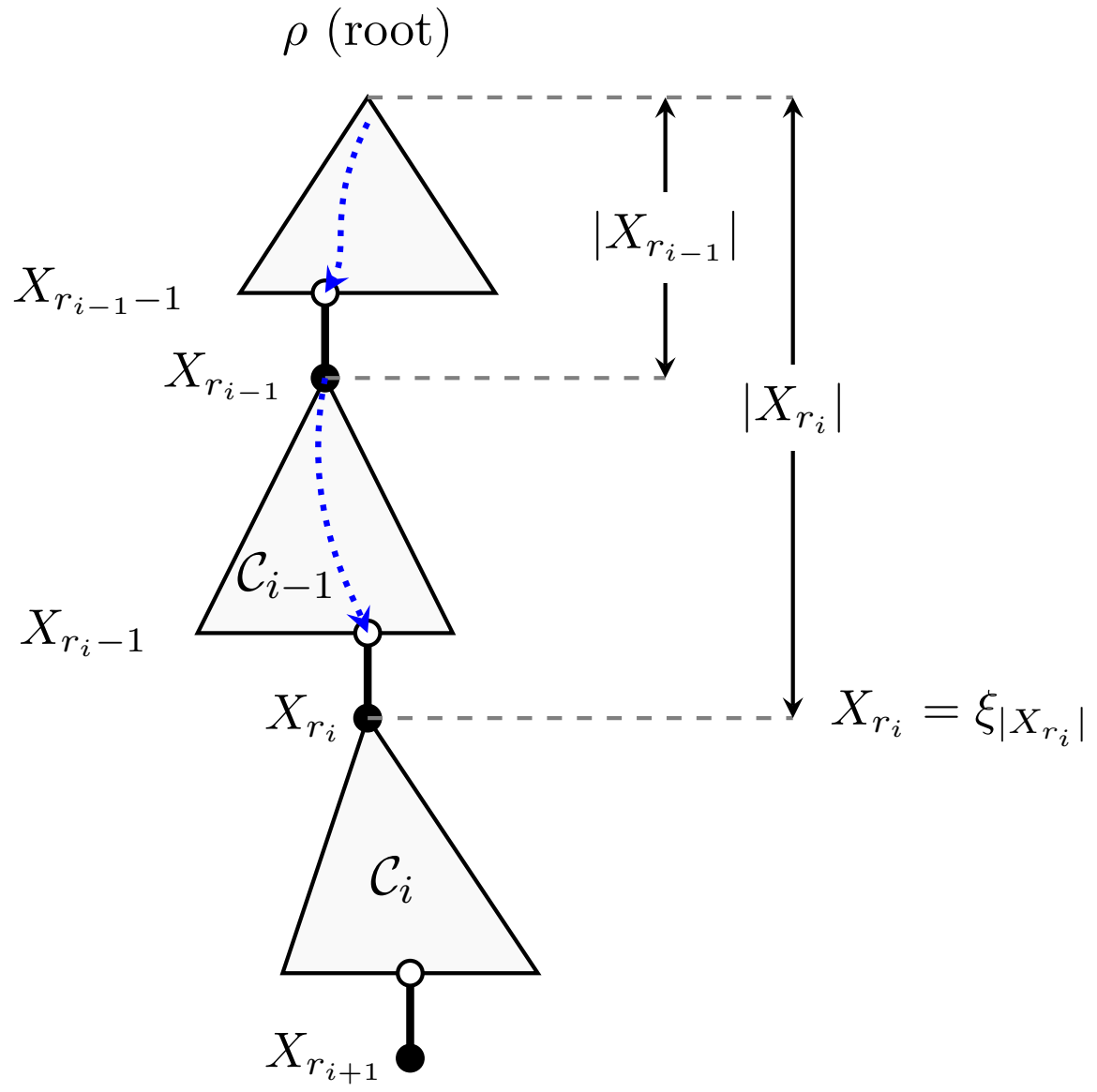

Let \(\Theta\) re-root the forward tree at \(v_1\), shift the walk by
\(r_1\), and retain \(b_1\) only as an artificial forbidden parent.
Every regeneration vertex lies on the terminal ray, since the walk never
returns across its incoming edge. Thus \(v_i=\xi_{|v_i|}\).

\begin{thm}[Regeneration product structure]\label{thm:regen-cycles}
Under \(\Pbb_\star\), the slabs \((\mathcal C_i)_{i\geq0}\) are
independent and identically distributed. For every \(k\geq1\), the
forward tree and shifted walk obtained by applying \(\Theta\) \(k\)
times have law \(\Pbb_\star\) and are independent of
\(\mathcal C_0,\ldots,\mathcal C_{k-1}\). In particular,
\[
    (R_i,T_i)_{i\geq0}\text{ are i.i.d.},\qquad
    \Ebb_\star e^{\varepsilon R_0}<\infty,\qquad 1\leq T_0\leq R_0
\]
for some \(\varepsilon>0\). Under the original law, the slabs
\((\mathcal C_i)_{i\geq1}\) are i.i.d.\ with this common law and are
independent of \(\mathcal I_0\).
\end{thm}

\begin{proof}
Under \(\Pbb_\star\), inspect the next step whenever the walk first
visits a vertex \(u\) with \(A_u<1/2\). After a candidate step \(u\to y\),
suspend the search until the walk returns to \(u\), or forever if it
does not. The inspected entrances are stopping times, and this search
finds the first accepted candidate: an accepted candidate inside a
suspended excursion would prevent the return, making \(u\to y\)
itself accepted.

At an inspected entrance, first work under the ordinary annealed law
\(\Pbb_{\mu_R}\), or \(\Pbb\) for the original root. Conditional on
the history, the mark of \(u\), and its offspring counts below and above
\(1/2\), the above-half child trees, in an independent random order,
are i.i.d.\ with the tree marginal of \(\Pbb_{\mu_R}\). A candidate
step selects one uniformly, so the selected tree retains this law
independently of all unselected child trees and the tree outside
\(\Tcal_u\). Here we use that the slab forgets child order and the
outgoing half-edge's age rank. The Markov property gives the same
independence jointly with the future walk stopped on returning to \(u\).
Conditioning on acceptance therefore gives an independent
\(\Pbb_\star\) forward pair, as in Lemma~\ref{lem:accepted-child}.
This accepted factorization also holds under \(\Pbb_\star\): on
histories avoiding the initial artificial parent, acceptance makes
its subsequent avoidance automatic.

At each inspected entrance, all earlier rejections are recorded in the
history. At the first successful cut, \(\mathcal C_0\) is determined by
the deleted exterior and the preceding path.
Proposition~\ref{prop:exp-moment-R} makes this
cut almost surely finite and gives the required exponential moment,
since \(R_0=r_1\). Thus, summing over the possible first successful
entrances, for bounded nonnegative \(F,G\),
\[
    \Ebb_\star[G(\mathcal C_0)\,F\circ\Theta]
    =\Ebb_\star[G(\mathcal C_0)]\Ebb_\star[F].
\]
Hence the shifted forward pair has law \(\Pbb_\star\) and is independent
of \(\mathcal C_0\); iteration proves the product structure. The claims
about \((R_i,T_i)\) follow since they are slab functionals and
\(1\leq T_0\leq R_0\). The same first-cut factorization under the
original law gives a \(\Pbb_\star\) forward pair independent of
\(\mathcal I_0\), proving the final assertion.
\end{proof}

\begin{thm}[Concentration of distance]\label{thm:distance-fluctuations}
Let
\begin{equation}\label{eq:def-speed}
    \nu:=\frac{\Ebb_\star T_0}{\Ebb_\star R_0}.
\end{equation}
Then \(0<\nu\leq1\), and 
\begin{equation}\label{eq:maximal-distance-tightness}
\lim_{\lambda\rightarrow\infty}\sup_{n\geq 1}\mathbb P\left(\max_{0\leq t\leq n}
        \left|\dist(\rho,X_t)-\nu\,t\right|>\lambda\sqrt{n}
    \right)
    =0
\end{equation}
\end{thm}
\begin{proof}
Theorem~\ref{thm:regen-cycles} gives \(1\leq T_0\leq R_0\), and
hence \(0<\nu\leq1\). Under the original law, the variables
\(T_i-\nu R_i\), \(i\geq1\), are independent and centered with finite
variance. If \(r_j\leq t<r_{j+1}\), then
\[
    |X_t|-\nu t
    =(|X_{r_1}|-\nu r_1)+\sum_{i=1}^{j-1}(T_i-\nu R_i)
      +(|X_t|-|X_{r_j}|-\nu(t-r_j)).
\]
The initial term is at most \(2r_1\) in absolute value, and the last
is at most \(2R_j\). Since \(r_j\geq j\), this gives, also for
\(t<r_1\),
\[
    \max_{t\leq n}\bigl||X_t|-\nu t\bigr|
    \leq2r_1+\max_{k\leq n}\Bigl|\sum_{i=1}^k(T_i-\nu R_i)\Bigr|
    +2\max_{1\leq i\leq n}R_i.
\]
Kolmogorov's maximal inequality and the exponential regeneration
tails imply, for \(\lambda\geq5\),
\[
    \Pbb\left(\max_{t\leq n}\bigl||X_t|-\nu t\bigr|>\lambda\sqrt n\right)
    \leq Ce^{-c\lambda\sqrt n}+\frac{C}{\lambda^2}+Cn e^{-c\lambda\sqrt n}.
\]
The right-hand side tends to zero uniformly in \(n\geq1\) as
\(\lambda\to\infty\).
\end{proof}

\begin{cor}[Logarithmic-time distance control]
\label{cor:logarithmic-time-distance}
For every \(c>0\),
\[
    \max_{t\leq\lfloor c\log n\rfloor}\bigl||X_t|-\nu t\bigr|
    =O_{\Pbb}(\sqrt{\log n}).
\]
Moreover, for every \(q>0\), there is \(C<\infty\) such that
\begin{equation}\label{eq:walk-near-ray}
    \Pbb\left(\max_{t\leq\lfloor c\log n\rfloor}
        \dist_{\Tcal}(X_t,\xi)>C\log\log n\right)
    \leq(\log n)^{-q}
\end{equation}
for all sufficiently large \(n\).
\end{cor}

\begin{proof}
The first assertion is \eqref{eq:maximal-distance-tightness}.
Every regeneration vertex lies on \(\xi\). Thus the distance from
the ray is at most \(r_1\) before the first regeneration, and at most
\(R_i\) during the \(i\)-th later interval. The exponential tails
and a union bound over at most \(\lfloor c\log n\rfloor+1\) intervals
prove \eqref{eq:walk-near-ray}.
\end{proof}

\medskip
\subsection{Concentration of entropy}
\label{subsec:entropy}

We first work under \(\Pbb_\star\), and define
\[
    Q_i:=W_{|v_{i+1}|}-W_{|v_i|},
    \qquad
    |v_j|=\sum_{i=0}^{j-1}T_i.
\]
By Theorem~\ref{thm:regen-cycles} and the invariance of
\eqref{eq:harmonic-kernel} under nonreturn conditioning,
\[
    Q_i=Q_0\circ\Theta^i,\qquad T_i=T_0\circ\Theta^i.
\]
Thus \((Q_i,T_i)_{i\geq0}\) is stationary. The variables \(Q_i\) need not be
independent, since their harmonic weights depend on the tree beyond
the next regeneration vertex. We control this dependence by computing
escape probabilities in finite descendant trees.

For a child \(y\) of \(x\), let \(h_x^{(\ell)}(y)\) be the probability,
starting from \(y\), of reaching descendant level \(\ell\) below \(y\)
before hitting \(x\). In particular, \(h_x^{(0)}(y)=1\). Define
\begin{equation}\label{eq:finite-depth-harmonic-kernel}
    \theta_x^{(\ell)}(y)
    :=\frac{h_x^{(\ell)}(y)}
       {\sum_{z:z^-=x}h_x^{(\ell)}(z)}.
\end{equation}
For any vertex \(v\), let \(W^{(\ell)}(v)\) be the sum of
\(-\log\theta_x^{(\ell)}(y)\) over the edges on the path from the
current root to \(v\), and put \(W_k^{(\ell)}=W^{(\ell)}(\xi_k)\).
Since each kernel sums to one,
\begin{equation}\label{eq:finite-depth-level-masses}
    \sum_{|v|=k}e^{-W^{(\ell)}(v)}=1,
    \qquad k\geq0.
\end{equation}

\begin{lem}[Finite-depth approximation]\label{lem:block-approx-error}
For every descendant edge and every \(\ell\geq0\),
\begin{equation}\label{eq:one-edge-harmonic-approximation}
    |\log\theta_x^{(\ell)}(y)-\log\theta_x(y)|
    \leq 2\varrho^\ell.
\end{equation}
Consequently,
\begin{equation}\label{eq:block-harmonic-approximation}
    |W_k-W_k^{(\ell)}|\leq 2k\varrho^\ell,
    \qquad
    |Q_i-Q_i^{(\ell)}|\leq 2\varrho^\ell T_i,
\end{equation}
where \(Q_i^{(\ell)}=W_{|v_{i+1}|}^{(\ell)}-W_{|v_i|}^{(\ell)}\).
\end{lem}

\begin{proof}
By the strong Markov property at descendant level \(\ell\) and
\eqref{eq:ancestor-return},
\[
    0\leq h_x^{(\ell)}(y)-h_x(y)\leq\varrho^{\ell+1}.
\]
Since \(h_x(y)\geq p_{\rm esc}=1/3\),
\[
    1\leq\frac{h_x^{(\ell)}(y)}{h_x(y)}
    \leq1+2\varrho^\ell.
\]
The ratio of the corresponding sums over the children of \(x\)
satisfies the same bounds.
Taking logarithms in \eqref{eq:finite-depth-harmonic-kernel} proves
\eqref{eq:one-edge-harmonic-approximation}. Summing over the relevant
paths proves \eqref{eq:block-harmonic-approximation}.
\end{proof}

The half-open slabs of Theorem~\ref{thm:regen-cycles} include their
complete side fringes. Thus
\[
    Q_i^{(\ell)}
    \in\sigma(\mathcal C_i,\ldots,\mathcal C_{i+\ell}).
\]
Indeed, only the terminal ray can leave a slab, and each regeneration
crossing advances at least one level. The boundary vertices in the
finite escape problem are absorbing, so their offspring need not be
included. In particular, \((Q_i^{(\ell)})_{i\geq0}\) is
\(\ell\)-dependent.

\begin{lem}[Second moment of harmonic block weights]\label{lem:Q0-L2}
One has \(\Ebb_\star Q_0^2<\infty\).
\end{lem}

\begin{proof}
At the first visit to a vertex \(v\), conditional on its age \(A_v=s\)
and type, its degree is at most \(m+1+N_v\), where \(N_v\), given \(\Gamma_v\), is
Poisson with mean \(\Lambda_v=\Gamma_v(s^{-(1-\chi)}-1)\).
The Poisson moment formula, the logarithmic moments of the Gamma
laws, and
\[
    \log(1+N_v)\leq\log(1+\Lambda_v)
       +\log(1+N_v/(1+\Lambda_v)),\qquad
    1+\Lambda_v\leq(1+\Gamma_v)s^{-(1-\chi)}
\]
give
\begin{equation}\label{eq:conditional-log-degree}
    \Ebb[\log^4(\deg(v)\vee2)\mid A_v=s,\typ(v)]
    \leq C(1+|\log s|^4).
\end{equation}

Consider first the walk under \(\Pbb_{\mu_R}\), stopped at its
artificial parent. By time \(n\), at most \(n+1\) vertices have been
visited, exposing at most \(m(n+1)\) left-child age ratios.
Their negative logarithms can be padded to a collection of
\(m(n+1)\) independent \(\operatorname{Exp}(\chi)\) variables.
Right edges increase age, so the negative logarithm of every visited
age is bounded by their sum plus \(\log2\). At a first visit, the mark and offspring of the
visited vertex have not been exposed. Applying
\eqref{eq:conditional-log-degree} at these successive visits, and
then \eqref{eq:transfer}, yields
\begin{equation}\label{eq:explored-degree-moment}
    \Ebb_\star\left[\max_{0\leq j\leq n}
        \log^4(\deg(X_j)\vee2)\right]\leq C(1+n)^5.
\end{equation}
The artificial parent is omitted from the maximum under the stopped
law; it is never visited under \(\Pbb_\star\).

Every parent of a ray edge contributing to \(Q_0\) is visited by
time \(R_0\). By \eqref{eq:log-theta-degree},
\[
    Q_0\leq R_0\left(\max_{j\leq R_0}\log(\deg(X_j)\vee2)+\log3\right).
\]
Cauchy--Schwarz, \eqref{eq:explored-degree-moment}, and the
exponential tail of \(R_0\) therefore give
\[
    \Ebb_\star Q_0^2
    \leq C\sum_{n\geq1}n^2(1+n)^{5/2}
                  \Pbb_\star(R_0=n)^{1/2}<\infty.
\]
\end{proof}

Define the entropy rate per ray level by
\begin{equation}\label{eq:def-entropy-rate}
    d:=\frac{\Ebb_\star Q_0}{\Ebb_\star T_0}\in(0,\infty).
\end{equation}
Finiteness follows from Lemma~\ref{lem:Q0-L2} and the exponential tail
of \(T_0\). For positivity, the terminal edge \(b_1\to v_1\) of the
first block is a candidate edge. Besides the right child \(v_1\), the
vertex \(b_1\) has a left child \(w\). Hence
\[
 \theta_{b_1}(v_1)
 \leq\frac{h_{b_1}(v_1)}{h_{b_1}(v_1)+h_{b_1}(w)}
 \leq\frac1{1+p_{\rm esc}}=\frac34.
\]
Thus \(Q_0\geq\log(4/3)\), which proves positivity.

\begin{prop}[Variance of harmonic block sums]\label{prop:cov-decay-B}
There are \(C,c>0\) such that
\[
    |\operatorname{Cov}_\star(Q_0,Q_j)|\leq Ce^{-cj},
    \qquad j\geq1.
\]
In particular,
\begin{equation}\label{eq:linear-block-variance}
    \operatorname{Var}_\star\left(\sum_{i=0}^{j-1}Q_i\right)
    \leq Cj.
\end{equation}
\end{prop}

\begin{proof}
Since \(T_0\leq R_0\) and Proposition~\ref{prop:exp-moment-R} gives
\(T_0\in L^2(\Pbb_\star)\), \eqref{eq:block-harmonic-approximation}
gives \(\|Q_i-Q_i^{(\ell)}\|_2\leq C\varrho^\ell\), uniformly in
\(i\). For every \(j\geq1\), \(Q_0^{(j-1)}\) depends only on
\(\mathcal C_0,\ldots,\mathcal C_{j-1}\), whereas \(Q_j\) is measurable
with respect to \(\sigma(\mathcal C_j,\mathcal C_{j+1},\ldots)\).
The latter follows by letting \(\ell\to\infty\) in the finite-depth
approximations \(Q_j^{(\ell)}\). The product structure of the slabs
therefore makes \(Q_0^{(j-1)}\) and \(Q_j\) independent. By stationarity
and Lemma~\ref{lem:Q0-L2},
\[
    |\operatorname{Cov}_\star(Q_0,Q_j)|
    =|\operatorname{Cov}_\star(Q_0-Q_0^{(j-1)},Q_j)|
    \leq\|Q_0-Q_0^{(j-1)}\|_2\|Q_j\|_2
    \leq C\varrho^{j-1}.
\]
Summing the covariances proves \eqref{eq:linear-block-variance}.
\end{proof}

\begin{thm}[Concentration of entropy]
\label{thm:original-W-tightness}
Under the original annealed law,
\begin{equation}\label{eq:entropy-tightness}
    \lim_{\lambda\to\infty}\sup_{k\geq1}
    \Pbb\bigl(|W_k-dk|>\lambda\sqrt{k}\bigr)=0.
\end{equation}
The same assertion holds under \(\Pbb_\star\), with depth measured
from its root.
\end{thm}

\begin{proof}
Work first under \(\Pbb_\star\). The regeneration depths and their
harmonic weights satisfy
\[
    \Ebb_\star|v_j|=j\,\Ebb_\star T_0,\qquad
    \Ebb_\star W_{|v_j|}=dj\,\Ebb_\star T_0,\qquad
    \operatorname{Var}_\star(|v_j|)
       +\operatorname{Var}_\star(W_{|v_j|})\leq Cj.
\]
Here the first variance bound uses the independent increments
\(T_i\), and the second is \eqref{eq:linear-block-variance}.
For \(\lambda\geq1\), choose the deterministic indices
\[
    j_-:=\max\{0,\lfloor k/\Ebb_\star T_0-\lambda\sqrt{k}\rfloor\},
    \qquad
    j_+:=\lceil k/\Ebb_\star T_0+\lambda\sqrt{k}\rceil.
\]
Chebyshev's inequality gives
\(\Pbb_\star(|v_{j_-}|>k\text{ or }|v_{j_+}|<k)\leq C/\lambda\).
It also bounds the probability that either
\(|W_{|v_{j_\pm}|}-dj_\pm\Ebb_\star T_0|>\lambda\sqrt{k}\) by
\(C/\lambda\). These bounds are uniform in \(k\geq1\), since
\(j_\pm\leq C(k+\lambda\sqrt{k})\); terms with \(j_-=0\) vanish.
On the complementary event, monotonicity gives
\[
    W_{|v_{j_-}|}\leq W_k\leq W_{|v_{j_+}|},
    \qquad |W_k-dk|\leq C(\lambda+1)\sqrt{k}.
\]
Letting \(\lambda\to\infty\) proves the assertion under \(\Pbb_\star\).

Under the original law, put \(L=|v_1|\) and let
\(\widehat W\) denote harmonic information measured from \(v_1\).
Theorem~\ref{thm:regen-cycles} and the invariance of
\eqref{eq:harmonic-kernel} under nonreturn conditioning imply that
\(\widehat W\) has the \(\Pbb_\star\)-law, is independent of \(L\),
and satisfies \(W_k=W_L+\widehat W_{k-L}\) when \(k\geq L\).
Monotonicity also covers \(k<L\), giving
\[
    |W_k-dk|\leq W_L+dL+
        \mathbf1_{\{L\leq k\}}|\widehat W_{k-L}-d(k-L)|.
\]
Thus the left-hand probability in \eqref{eq:entropy-tightness}, with
threshold \(2\lambda\sqrt{k}\), is at most
\[
    \Pbb(W_L+dL>\lambda)+
    \sup_{s\geq1}\Pbb_\star(|W_s-ds|>\lambda\sqrt{s}).
\]
The first term tends to zero because \(r_1<\infty\) almost surely and
the tree is locally finite. The second tends to zero by the
post-regeneration assertion.
\end{proof}
    
For integer \(\ell\geq3\log\log n/\log(3/2)\),
\(|W_k-W_k^{(\ell)}|=O((\log n)^{-2})\) uniformly for
\(k=O(\log n)\). Hence, if \(k_n=d^{-1}\log n+O(\sqrt{\log n})\),
\[
    W_{k_n}^{(\ell)}=dk_n+O_{\Pbb}(\sqrt{\log n}),
    \qquad
    e^{-W_{k_n}}=n^{-1}\exp\{O_{\Pbb}(\sqrt{\log n})\}.
\]
The rate \(d\) is per ray level. For the ray vertex
\(\xi_{|X_t|}\), distance concentration and monotonicity give
\(W_{|X_t|}=\nu dt+O_{\Pbb}(\sqrt t)\).
Indeed, on \(\{||X_t|-\nu t|\leq K\sqrt t\}\), sandwich this weight
between those at the deterministic levels
\(\max\{0,\lfloor\nu t-K\sqrt t\rfloor\}\) and
\(\lceil\nu t+K\sqrt t\rceil\), and apply
\eqref{eq:entropy-tightness} before letting \(K\to\infty\).

\subsection{Regularity and truncation}
\label{subsec:regularity}

Throughout this subsection, fix \(c>(\nu d)^{-1}\) and put
\(t_n:=\lfloor c\log n\rfloor\) and \(h_n:=\lceil c_h\ell_n\rceil\),
where \(\ell_n=\log\log n\) and the constant \(c_h>0\) is chosen in
Proposition~\ref{prop:no-truncation}.
Write \(B_r^+(v)\) for the radius-\(r\) ball in the descendant tree
\(\Tcal_v\). The following consequence of
Lemma~\ref{lem:regularity-PPT} controls the full neighborhood of the
walk.

\begin{prop}[Regularity of the walk neighborhood]
\label{thm:global-age-volume}
Fix \(A,C,\zeta,q>0\). Let \(X\) start at the root of a PPT
with prescribed age \(s\in[e^{-A\ell_n^2},1]\), mark
\(0<g\le C\ell_n^2\), and type, allowing the ordinary root.
Write \(\Pbb\) for its annealed law; if an entrance parent is attached,
condition never to hit it and omit it from all balls. Put
\[
    \mathcal W_n:=\bigcup_{0\leq t\leq t_n}B_{2h_n}(X_t).
\]
There are constants \(C_g\ge C\) and \(C_V<\infty\) such that,
uniformly over these root data, for all sufficiently large \(n\),
with probability at least \(1-(\log n)^{-q}\),
\[
    \min_{v\in\mathcal W_n}A_v\geq s e^{-\zeta\ell_n^2},
    \qquad
    \max_{v\in\mathcal W_n}\Gamma_v\leq C_g\ell_n^2,
    \qquad
    |\mathcal W_n|\leq e^{C_V\ell_n^2}.
\]
\end{prop}

\begin{proof}
Start the trials at first visits after the first step. The proof of
Proposition~\ref{prop:exp-moment-R} gives a uniform exponential tail
for \(r_1\), also under nonreturn conditioning, whose density is at
most \(p_{\rm esc}^{-1}\). The forward tree and shifted walk at each
subsequent regeneration time have law \(\Pbb_\star\).
Thus, for \(M=C_q\ell_n\) with \(C_q\) sufficiently large,
\[
    \Pbb\left(r_1\vee\max_{1\leq i\leq t_n}R_i>M\right)
    \leq\tfrac12(\log n)^{-q}.
\]
On the complementary event, each \(B_{2h_n}(X_t)\), \(t\leq t_n\), is
contained in a forward ball of radius \(4h_n+M\) about \(\rho\) or a
regeneration vertex. Indeed, choose the latest already reached
regeneration vertex \(v_j\) with \(|v_j|\leq |X_t|-2h_n\).
It is an ancestor of \(X_t\), and \(\dist(v_j,X_t)\leq 2h_n+M\):
either it is the current regeneration vertex, or the next
regeneration depth exceeds \(|X_t|-2h_n\), with gap at most \(M\).
Also \(\dist(v_j,X_t)\geq 2h_n\), so \(B_{2h_n}(X_t)\subseteq\Tcal_{v_j}\).
If no such vertex exists, \(|X_t|\leq 2h_n+r_1\leq 2h_n+M\), and the
ball is covered from \(\rho\). Therefore
\[
    \mathcal W_n\subseteq B_{4h_n+M}(\rho)
        \cup\bigcup_{i=1}^{t_n}B_{4h_n+M}^+(v_i).
\]

Apply Lemma~\ref{lem:regularity-PPT} to the initial ball with relative
age cutoff \(s e^{-\zeta\ell_n^2}\). Its proof applies to a prescribed
root type and mark: the first-generation bound gains only a factor
\(C'(1+g)\), absorbed by the slack in the exponential estimates.
The root itself satisfies the mark bound since \(C_g\ge C\).
Each descendant tree \(\Tcal_{v_i}\) has the tree marginal of
\(\Pbb_\star\), with
density at most \(p_{\rm esc}^{-1}\) relative to the tree marginal of
\(\Pbb_{\mu_R}\), and root age in \((1/2,1)\). Apply the same lemma there with
\(\zeta/2\); for large \(n\), its age bound is at least
\(\tfrac12e^{-\zeta\ell_n^2/2}\ge s e^{-\zeta\ell_n^2}\).
After increasing \(C_g,C_V\), the mark and volume bounds hold for
all these balls except with probability
\(O(\log n)e^{-c_{\mathrm{reg}}\ell_n^2}=o((\log n)^{-q})\).
Their total volume is still at most \(e^{C_V\ell_n^2}\).
These marginal estimates are used without conditioning on the
regeneration-length event.
\end{proof}

We now specify the truncation used in Section~\ref{sec:coupling-after-burn-in}.
Choose \(C_A>A\) and apply Proposition~\ref{thm:global-age-volume}
with \(\zeta=C_A-A\), so that \(s e^{-\zeta\ell_n^2}\ge e^{-C_A\ell_n^2}\).
With the corresponding mark and volume constants, set
\begin{equation}\label{eq:ppt-truncation-scales}
    a_n=e^{-C_A\ell_n^2},\qquad
    g_n=C_g\ell_n^2,\qquad
    V_n=e^{C_V\ell_n^2}
    ,\qquad
    \omega_n=\log n-\tfrac12\sqrt{\log n}.
\end{equation}
For a vertex \(v\), let \(v^{[h_n]}\) be its ancestor \(h_n\) edges
above it, or the root if \(|v|<h_n\). Construct the subtree in a fixed
breadth-first order, starting with the root in the queue. For each
vertex \(v\) taken from the queue, add its children to the queue if
\begin{equation}\label{eq:ppt-truncation-rule}
\begin{gathered}
    A_v\geq a_n,\qquad \Gamma_v\leq g_n,\qquad 
    W^{(h_n)}(v^{[h_n]})\leq \omega_n,\qquad
    |B_{h_n}^+(v^{[h_n]})|\leq V_n.
\end{gathered}
\end{equation}
If any of these conditions fails, declare \(v\) \emph{truncated}
and retain it as a boundary leaf. Denote by $\widetilde{\mathcal T}$
the resulting rooted subtree. All quantities refer to the full PPT;
the finite-depth tests can be evaluated by
completing the tree through the current level (and through level
\(h_n\) for the initial ball). At depths at least \(h_n\), the weight
test uses no offspring of \(v\). The age and mark bounds also imply
\(\|\rho_v\|_\infty\leq(1-\chi)g_n/a_n\), so a separate intensity
cutoff is unnecessary.
Write \(\mathcal B_{\rm tr}(v)\) for the event that \(v\) fails one of
the tests in \eqref{eq:ppt-truncation-rule}.

By \eqref{eq:finite-depth-level-masses}, at most
\(ne^{-\sqrt{\log n}/2}\) vertices \(v\) of a given level satisfy
\(W^{(h_n)}(v)\leq \omega_n\). Each untruncated vertex of
\(\widetilde{\mathcal T}\) at level \(k\geq h_n\) descends from such a
vertex at level \(k-h_n\), and at most \(V_n\) of these descendants
have the same ancestor.
The initial levels are bounded by the root-ball test. Thus the
number of untruncated vertices through depth \(t_n\) is at most
\begin{equation}\label{eq:truncated-exploration-size}
    (t_n+1)e^{\omega_n}V_n
    =n\exp\{-\tfrac12\sqrt{\log n}
               +O((\log\log n)^2)\},
    \qquad t_n=O(\log n).
\end{equation}
This bound excludes truncated boundary vertices; additional PA labels
are counted separately in the finite-graph comparison.

\begin{prop}[No truncation before the target time]
\label{prop:no-truncation}
Fix \(A,C>0\) and use any prescribed-root law of
Proposition~\ref{thm:global-age-volume}, with age
\(s\ge e^{-A\ell_n^2}\) and mark \(g\le C\ell_n^2\). Let
\[
    \ell_0=\left\lfloor
        \frac{\log n-2\gamma\sqrt{\log n}}{d}\right\rfloor,
    \qquad
    \ell_1=\left\lfloor
        \frac{\log n-\gamma\sqrt{\log n}}{d}\right\rfloor,
    \qquad
    t_0=\left\lfloor\frac{\ell_0+\ell_1}{2\nu}\right\rfloor .
\]
The constants \(c_h,C_A,C_g,C_V\), with \(C_A>A\) and \(C_g\ge C\),
can be chosen so that there is a function \(\eta(\gamma)\), tending
to zero as \(\gamma\to\infty\), for which,
for every fixed \(\gamma>0\), uniformly over the root data,
\[
    \Pbb\left(\bigcup_{0\leq t\leq t_0}
        \mathcal B_{\rm tr}(X_t)\right)
    \leq\eta(\gamma)+o(1).
\]
Moreover, for \(\tau:=\inf\{t:|X_t|=\ell_1\}\),
\begin{equation}\label{eq:ppt-stopped-no-truncation}
    \Pbb\left(\{\tau>\lceil t_n^\star\rceil\}
        \cup\bigcup_{0\leq t\leq\tau}\mathcal B_{\rm tr}(X_t)\right)
    \leq\eta(\gamma)+o(1).
\end{equation}
Here \(o(1)\) is taken as \(n\to\infty\) with \(\gamma\) fixed.
\end{prop}

\begin{proof}
The argument of Lemma~\ref{lem:Q0-L2}, with initial age \(s\) and
mark \(g\), gives
\[
    \Ebb W_{|v_1|}^2\le C'(1+|\log s|+\log(1+g))^2.
\]
The initial harmonic information is therefore
\(O_{\Pbb}(\ell_n^2)=o_{\Pbb}(\sqrt{\log n})\), uniformly.
Since \(r_1\) has a uniform exponential tail and the forward law after
the first regeneration is \(\Pbb_\star\), the distance and entropy
estimates used below hold uniformly at times and depths of order \(\log n\).

For each fixed \(\gamma\), one has \(t_0\leq t_n\) for all sufficiently
large \(n\). Choose \(c_h\) sufficiently large that
\eqref{eq:walk-near-ray}, with \(q=1\), gives
\[
    \Pbb\left(\max_{t\leq t_n}\dist(X_t,\xi)>h_n\right)=o(1)
\]
and \((\log n)\varrho^{h_n}=o(1)\). Apply
Proposition~\ref{thm:global-age-volume} with \(q=1\) and
\(\zeta=C_A-A\). Its bounds hold on \(\mathcal W_n\) with probability
\(1-o(1)\), and every age there is at least \(a_n\).

By \eqref{eq:maximal-distance-tightness}, there is
\(\eta_1(\gamma)\to0\) such that
\[
    \Pbb\left(\max_{t\leq t_0}|X_t|>\ell_1\right)
    \leq\eta_1(\gamma)+o(1).
\]
Indeed,
\[
    \ell_1-\nu t_0
    \geq\frac{\ell_1-\ell_0}{2}
    \geq \frac{\gamma}{2d}\sqrt{\log n}-1.
\]

Suppose the preceding regularity, proximity, and height events hold.
For \(t\leq t_0\), set \(v=X_t\) and \(u=v^{[h_n]}\). Then
\(B_{h_n}^+(u)\subseteq B_{2h_n}(v)\), so \(v\) passes the age,
mark, and ball tests. Moreover, \(\dist(v,\xi)\leq h_n\) implies
that \(u\in\xi\). Since \(|u|\leq\ell_1\), monotonicity and
Lemma~\ref{lem:block-approx-error} give
\[
    W^{(h_n)}(u)\leq W_{\ell_1}^{(h_n)}
        \leq W_{\ell_1}+2\ell_1 \varrho^{h_n}
        =W_{\ell_1}+o(1).
\]
Consequently,
\[
\begin{split}
    \Pbb\left(\bigcup_{0\leq t\leq t_0}
        \mathcal B_{\rm tr}(X_t)\right)
    &\leq \eta_1(\gamma)+o(1)
       +\Pbb\left(W_{\ell_1}^{(h_n)}>\log n-\tfrac12\sqrt{\log n}\right)\\
    &\leq \eta_1(\gamma)+o(1)
       +\Pbb\left(
          W_{\ell_1}-d\ell_1>
          \left(\gamma-\tfrac12\right)\sqrt{\log n}-o(1)
       \right).
\end{split}
\]
Since \(\ell_1\sim d^{-1}\log n\), Theorem~\ref{thm:original-W-tightness}
bounds the last probability by \(\eta_2(\gamma)+o(1)\), where
\(\eta_2(\gamma)\to0\).

For the stopped-time bound, distance concentration gives
\[
    \Pbb(\tau>\lceil t_n^\star\rceil)
    \leq\Pbb(|X_{\lceil t_n^\star\rceil}|<\ell_1)
    \leq\eta_3(\gamma)+o(1),
\]
where \(\eta_3(\gamma)\to0\), since
\(\nu\lceil t_n^\star\rceil-\ell_1
=\gamma\sqrt{\log n}/d+O(1)\).
On \(\{\tau\leq\lceil t_n^\star\rceil\}\), the regularity and proximity
events apply through \(\tau\), as \(\lceil t_n^\star\rceil\leq t_n\).
Since \(|X_t|\leq\ell_1\) for \(t\leq\tau\), the same comparison with
\(W_{\ell_1}^{(h_n)}\) bounds the probability of a truncation on this
event by \(\eta_2(\gamma)+o(1)\). Taking
\(\eta=\eta_1+\eta_2+\eta_3\) proves both bounds.
\end{proof}

For the ordinary PPT, \(A_\rho\) is uniform on \((0,1)\), so
\(\Pbb(A_\rho<e^{-A\ell_n^2})=e^{-A\ell_n^2}\); the Gamma tail also gives
\(\Pbb(\Gamma_\rho>C\ell_n^2)=o(1)\). Averaging over the remaining root
data yields the same no-truncation bound. With an entrance parent,
conditioning on reaching depth \(R\) before returning, in place of
never returning, adds at most \(\varrho^{R+1}\) to the error by
\eqref{eq:ancestor-return}.

\medskip

\section{Coupling preferential attachment to the P\'olya--point tree}
\label{sec:coupling}

We couple an adaptive sequence of radius-$O(\log\log n)$ PA explorations
to PPT explorations, retaining the information revealed by every trial.
We then construct coupled regions that the SRW does not leave until of order mixing time.

Throughout this section, $m\ge2$ and $\alpha\in[0,1)$ are fixed.
The symbols $v,w$ denote PPT vertices, $A_v$ their ages, $\Gamma_v$
their marks, and $K_v\in[n]$ their matched PA birth labels.

\subsection{Adaptive breadth-first coupling}
\label{subsec:ppt-pa-coupling}

Put $\ell_n=\log\log n$ and set
\begin{equation}
\begin{gathered}
 r_n=\lceil c_r\ell_n\rceil,\qquad
 b_n=e^{-c_0\ell_n^2},\qquad
 k_n^\circ=\lceil ne^{-c_1\ell_n^2}\rceil,\\
 B_n=\lceil e^{C_B\ell_n^2}\rceil,\qquad
 M_n=\lfloor ne^{-c_M\sqrt{\log n}}\rfloor,\qquad
 \bar g_n=C_G\log n.
\end{gathered}
\label{eq:coupling-scales}
\end{equation}
Each trial starts at a root of age at least $b_n$ and adds at most
$B_n$ PA labels to the exploration, within a total budget $M_n$.
All matched labels and targets queried during the trials are at least
$k_n^\circ$; the initial
record may contain complete stars at smaller labels.
The cap $g_n$ from \eqref{eq:ppt-truncation-scales} controls new PPT
marks, while $\bar g_n$ bounds $j\psi_j$ at labels $j\ge k_n^\circ$
whose marks have been exposed.

Fix $c_r,c_0,c_M,C_{\mathrm{test}}>0$, where
$e^{C_{\mathrm{test}}\ell_n^2}$ is the allowed number of local coupling trials.
Choose $\zeta>0$ sufficiently large and then take
\begin{equation}
 c_1>c_0+\zeta,\qquad
 C_B>(1-\chi)c_0+\zeta.
\label{eq:coupling-constant-order}
\end{equation}
We choose $C_B$ first and allow $c_1$ and $C_g$ to be increased below.
With the constant $C_g$ in \eqref{eq:ppt-truncation-scales}
sufficiently large, Lemma~\ref{lem:regularity-PPT} gives
the following bound. A radius-$(r_n+1)$ PPT ball with prescribed root
age $A\ge b_n$ has at most $B_n/2$ vertices, all with ages at least
$Ae^{-\zeta\ell_n^2}>8k_n^\circ/n$ and marks at most $g_n$, except
with probability $Ce^{-(C_{\mathrm{test}}+2)\ell_n^2}$. Increasing
$\zeta$ in the age and volume estimates, and then $C_g$ in the mark
estimate, gives the exponent $C_{\mathrm{test}}+2$.
The same bound holds for a forward ball with
prescribed root type and mark at most $g_n$: the first-generation
factor $C(1+g_n)$ is absorbed by the slack in
\eqref{eq:coupling-constant-order}.

Set $\varepsilon_n=n^{-1/3}$, and let
$\varepsilon_n^{\mathrm{BG}}$ be as in
\eqref{eq:beta-gamma-error}.

\paragraph{The adaptive exploration.}
The first coupling trial may start from a nonempty record of earlier
PA exploration.
Let $\mathcal A$ record all information revealed so far about the PA
graph and its urn marks. Keep this information even when a comparison
fails. Let $\mathcal H$ also contain the revealed PPT data and auxiliary choices.
To \emph{register} a PA vertex is to add its birth label to $E$;
put $t=|E|$. Register each root and both endpoints of every revealed
PA edge. Learning that a slot does not attach to a given target does
not by itself register its sender. Once registered, a label remains
in $E$, including after failed comparisons.

A registered vertex is \emph{explored} once its complete star, consisting
of all incident edges with their slot identities and multiplicities,
is recorded in $\mathcal A$. This depends only on the PA record and
does not require a successful PPT comparison. Before its complete star
is recorded, a vertex may already have recorded incidences or an exposed
mark. The set $E$ includes all registered labels, whether or not their
complete stars have been recorded.

While constructing this initial record, we may learn whether a vertex's
degree exceeds a given threshold before recording its complete star.
Auxiliary choices use independent randomness. Let $E_0$ consist of
the labels below $k_n^\circ$ that were queried as targets before the
first coupling trial, and put
\[
 \vartheta_\ell=\frac{\varphi_\ell}{S_{k_n^\circ-1}}
 =\psi_\ell\prod_{\ell<h<k_n^\circ}(1-\psi_h),\qquad \ell\in E_0.
\]
These quantities involve only urn variables below $k_n^\circ$.

A record $\mathcal H$ of this exploration is \emph{admissible} if it
satisfies the following conditions:
\begin{enumerate}[label=\textup{(\roman*)},leftmargin=2.2em]
\item Every target queried before the first coupling trial and every
vertex whose degree was tested has its complete star recorded. Thus all earlier
degree-test answers are determined by $\mathcal A$. The values
$\vartheta_\ell$, $\ell\in E_0$, are exposed and satisfy
\[
 \sum_{\ell\in E_0}\vartheta_\ell\le\frac{t\bar g_n}{k_n^\circ}.
\]
Once the trials begin, no target below $k_n^\circ$ is queried. Low endpoints
first disclosed thereafter are registered, but their complete stars
are not recorded and their marks remain unexposed.
\item Every queried target $j\ge k_n^\circ$ has its actual mark $Z_j$
exposed. Every exposed high-label mark satisfies $j\psi_j\le\bar g_n$.
\item The retained matching is injective, preserves types and PA slots,
and has PA labels at least $k_n^\circ$. Every retained pair, including
boundary pairs, satisfies
\[
 A_v\ge\frac{8k_n^\circ}{n},\qquad
 Z_{K_v}=\Gamma_v\le g_n,\qquad
 Z_{K_v}\in\mathcal G_{K_v}^{\mathrm{BG}}
                    (\varepsilon_n^{\mathrm{BG}}).
\]
\end{enumerate}
Admissibility is determined by the revealed record.
The root and label-budget hypotheses are imposed
separately in the lemmas below. Since $E_0$ stays fixed and $t$ can only
increase, its aggregate bound persists throughout the coupling trials.

During trials, expose $Z_j$ before querying a target $j\ge k_n^\circ$.
Complete any mark coupling before checking $j\psi_j\le\bar g_n$.
In each star, compare all left positions, then all right positions,
and finally all new child marks. Do not query a new child as a target
or interleave another star before completing these comparisons. Stop
a trial at its first failed comparison or regularity check, or before
exceeding its label budget.

Write $\Delta_v=|K_v/(nA_v)-1|$ for the label discrepancy, and set
\begin{equation}
 \Lambda_n=g_n(8k_n^\circ/n)^{-(1-\chi)},\qquad
 H_n=\begin{cases}
 1,&\chi>1/2,\\
 1+\log(n/k_n^\circ),&\chi=1/2.
 \end{cases}
\label{eq:coupling-envelopes}
\end{equation}
The effect of the accumulated PA record is controlled by
\[
 \eta_n(t)=\frac{(t+1)\bar g_n}{k_n^\circ}
          +\frac{nt\bar g_n}{(k_n^\circ)^2}.
\]
Uniformly for $t\le M_n$, one has $\eta_n(t)=o(1)$.

\begin{lem}[Adaptive BFS coupling]\label{lem:ext-coup}
Fix $C_0<\infty$. From an admissible record $\mathcal H$ with
$|E|\le M_n-B_n$, choose a matched root $v_0$ whose complete PA star
has not yet been recorded, before revealing its remaining offspring.
Its only recorded PA incidence is its
entrance, or none for an initial root. Assume that its coupled mark is
good, its prescribed PPT type is consistent with the entrance, and
\[
 A_{v_0}\ge b_n,\qquad
 \Delta_{v_0}\le C_0\sqrt{\eta_n(M_n)}.
\]
There is a coupling to the forward radius-$r_n$ PPT exploration with
the prescribed conditional PA and PPT laws, adding at most $B_n$
labels. Uniformly over these records, including after earlier failed
comparisons, its conditional probability of success is at least
$1-e^{-(C_{\mathrm{test}}+1)\ell_n^2}$ for all sufficiently large $n$.
On success:
\begin{enumerate}[label=\textup{(\arabic*)},leftmargin=2.2em]
\item The extended record is admissible. The matching is injective,
preserves types and slots, and matches every edge incident to a vertex
at depth less than $r_n$.
\item Every vertex $w$ in the extension satisfies
$A_w\ge A_{v_0}e^{-\zeta\ell_n^2}>8k_n^\circ/n$ and $K_w\ge k_n^\circ$.
\item Every newly matched vertex, including boundary vertices,
satisfies $Z_{K_w}=\Gamma_w\le g_n$ and
$Z_{K_w}\in\mathcal G_{K_w}^{\mathrm{BG}}(\varepsilon_n^{\mathrm{BG}})$.
\item The label discrepancy satisfies
\begin{equation}
 \max_w\Delta_w\le\Delta_{v_0}+Cr_n\eta_n(M_n).
\label{eq:extension-position-error}
\end{equation}
\end{enumerate}
All PA data revealed during a failed trial remain in the record.
\end{lem}

\begin{proof}
We follow the breadth-first exploration of
\cite[Lemma~3.5 and Proposition~3.6]{PA-local-limit}.
Explore vertices of depth less than $r_n$ under the preceding
comparison and stopping rules. Retain all generated PPT data. Defer completion
of any reference subtree that may be used in a later comparison until
the adaptive sequence ends.

Each PA answer is sampled conditional on $\mathcal A$ and retained.
Each PPT comparison uses only the preceding record, these answers,
and independent auxiliary randomness. Hence, also after failed comparisons,
\[
 \mathcal L(\text{unrevealed PA variables}\mid\mathcal H)
 =\mathcal L(\text{unrevealed PA variables}\mid\mathcal A).
\]
For the one-trial estimates, write $v_0$ for the current root and
$M$ for the bound on $|E|$ during this trial.

\paragraph{The accumulated PA information.}
The estimates below use the current value $t=|E|$ and include
the initial record and all earlier trials.
For $\ell<k_n^\circ\le k$,
$p_{k\to\ell}=\vartheta_\ell\prod_{k_n^\circ\le h<k}(1-\psi_h)$.
For each slot whose endpoint has not been revealed, let $F$ be the
set of target labels excluded by previous queries of that slot. Hence,
for every urn compatible with an admissible record and every sender
$k\ge k_n^\circ$,
\begin{equation}
 \sum_{\ell\in F}p_{k\to\ell}
 \le\sum_{\substack{\ell\in F\\\ell\ge k_n^\circ}}\psi_\ell
       +\sum_{\ell\in E_0}\vartheta_\ell
 \le\frac{2t\bar g_n}{k_n^\circ}=o(1).
\label{eq:left-forbidden-bound}
\end{equation}
Fix an unexposed coordinate $x=\psi_j$, $j\ge k_n^\circ$, and
condition on all other urn variables. Recall that the prior law of
$\psi_j$ is $\operatorname{Beta}(a,d_j)$, where
$d_j=(2j-3)m+2mu(j-1)$. Among slots with revealed endpoints, let $s_j$
count those ending at $j$, and let $N_j$ count those with target label
$\ell$ and sender label $k$ satisfying $\ell<j<k$.
Every slot with a revealed endpoint has its sender in $E$, so
$\sum_j s_j\le mt$ and $N_j\le mt$; endpoints below $k_n^\circ$
are included in these counts.

For each slot $e$ with sender $k>j$ whose endpoint has not been revealed,
write $F_e$ for its excluded target set. Since every queried high target has its mark
exposed, $j\notin F_e$. Define
\[
\begin{aligned}
 \alpha_e&:=\sum_{\substack{\ell\in F_e\\j<\ell<k}}
      \psi_\ell\prod_{\ell<h<k}(1-\psi_h),\\
 \beta_e&:=\sum_{\substack{\ell\in F_e\\\ell<j}}
      \psi_\ell\prod_{\substack{\ell<h<k\\h\ne j}}(1-\psi_h).
\end{aligned}
\]
These coefficients do not depend on $x$: $\alpha_e$ is the excluded
mass above $j$, and $(1-x)\beta_e$ is the excluded mass below $j$.
The avoidance probability factors as
\[
 1-\alpha_e-(1-x)\beta_e
 =(1-\alpha_e-\beta_e)(1+b_ex),\qquad
 b_e:=\frac{\beta_e}{1-\alpha_e-\beta_e}.
\]
The values $\vartheta_\ell$ do not depend on $x$, and removing the
factor $1-x$ from either a low- or high-target probability still leaves
$\alpha_e+\beta_e\le2t\bar g_n/k_n^\circ$. Hence
$b_e\le Ct\bar g_n/k_n^\circ$.
Slots with sender at most $j$ contribute factors independent of $x$;
products and sums over $e$ below range over slots with sender greater
than $j$ whose endpoints have not been revealed. The likelihood of the PA record gives
\[
 \mathcal L(\psi_j\mid\mathcal A,(\psi_i)_{i\ne j})(dx)
 =\frac{h_j(x)}{\mathbb E_Qh_j(X)}\,Q(dx),\qquad
 Q=\operatorname{Beta}(a+s_j,d_j),\qquad
 h_j(x)=(1-x)^{N_j}\prod_e(1+b_ex).
\]
There are at most $mn$ such factors. Put
$\lambda_j=\sum_e b_e$ and $L=t+nt\bar g_n/k_n^\circ$;
then $N_j+\lambda_j\le CL$ and $L/k_n^\circ=o(1)$.

Likelihood-ratio comparison bounds this posterior stochastically
between $\operatorname{Beta}(a+s_j,d_j+N_j)$ and
$\operatorname{Beta}(a+s_j,d_j-\lambda_j)$.
For the upper bound, the logarithmic derivative of the density ratio is
\[
 -\frac{N_j+\lambda_j}{1-x}
 +\sum_e\frac{b_e}{1+b_ex}\le0;
\]
the lower bound follows because $\prod_e(1+b_ex)$ is increasing.
The beta moment formula therefore gives, for $r=1,2$,
\[
 \left|\log\frac{
   \Ebb[(1-\psi_j)^r\mid\mathcal A,(\psi_i)_{i\ne j}]
 }{(d_j)_r/(c_j)_r}\right|
 \le C\left(\frac{s_j}{j}+\frac{(1+s_j)L}{j^2}\right),
\]
where $(x)_r$ is the rising factorial. For an exposed coordinate,
the corresponding logarithmic error is at most $C\bar g_n/j$.
Summing over coordinates gives $O(\eta_n(t))$, since
$\sum_j s_j\le mt$ and $\sum_{j\ge k_n^\circ}j^{-2}\le C/k_n^\circ$.
Successive conditional expectations consequently yield
\[
 \Ebb\left[\prod_j(1-\psi_j)^{r_j}\,\middle|\,\mathcal A\right]
 =\prod_j\frac{(d_j)_{r_j}}{(c_j)_{r_j}}
   \exp\{O(\eta_n(t))\},\qquad r_j\in\{1,2\},
\]
for any subset of coordinates above $k_n^\circ$, without conditional
independence. In particular, for
$T_{u,v}=\prod_{u<h<v}(1-\psi_h)$ and $k_n^\circ\le u<v\le n$,
\begin{equation}
 \Ebb[T_{u,v}\mid\mathcal A]=(u/v)^\chi(1+O(\eta_n(t))),\qquad
 \Ebb\left[
   \left(\frac{T_{u,v}}{(u/v)^\chi}-1\right)^2
   \,\middle|\,\mathcal A\right]\le C\eta_n(t).
\label{eq:conditional-urn-products}
\end{equation}
Here the unconditioned first moment is
$(u/v)^\chi(1+O(u^{-1}))$, and the ratio of its second moment
to the square of its first is $1+O(u^{-1})$.

\paragraph{Left children.}
At the current pair $(k,v)$, suppose $\Delta_v\le1/4$. Consider a
remaining left slot $(k,i)$ with endpoint $J:=J_{k,i}$ and excluded target set $F$.
Before revealing $J$, put $F_{\mathcal A}(j)=\Pbb(J\le j\mid\mathcal A)$.
For $k_n^\circ\le j<k$, conditional on the urn,
\[
 \Pbb(J\le j\mid\psi,\mathcal A)
 =\frac{S_j}{S_{k-1}}\frac{1-a_j}{1-q_F},\qquad
 a_j=\sum_{\ell\in F,\,\ell\le j}p_{j+1\to\ell},\quad
 q_F=\sum_{\ell\in F}p_{k\to\ell}.
\]
Both excluded masses are bounded by \eqref{eq:left-forbidden-bound}.
Integrating over the urn and using \eqref{eq:conditional-urn-products}
gives the uniform relative estimate
\[
 F_{\mathcal A}(j)=(j/k)^\chi(1+O(\eta_n(t))).
\]
Reveal and retain the actual endpoint $J$, then draw
\[
 V\mid\mathcal H,J\sim
 \operatorname{Unif}[F_{\mathcal A}(J-1),F_{\mathcal A}(J)),
 \qquad A_w=A_vV^{1/\chi}.
\]
The interval lengths are exactly the conditional endpoint
probabilities, so $V$ is uniform given the preceding $\mathcal H$.
This kernel uses no hidden urn information beyond the retained $J$.
Register each new endpoint immediately and repeat for every
remaining left slot, omitting the entrance slot for a type-$\tR$ vertex.

On the reference age requirement, $y=kV^{1/\chi}\ge6k_n^\circ$.
The CDF estimate at $2k_n^\circ$ forces $J>2k_n^\circ$; applying it
at $J-1$ and $J$ gives
$(J-1)(1-C\eta_n(t))\le y\le J(1+C\eta_n(t))$. Consequently
\begin{equation}
 \left|\frac{J}{kV^{1/\chi}}-1\right|\le C\eta_n(t),\qquad
 \Delta_w\le\Delta_v+C\eta_n(t).
\label{eq:left-error-recursion}
\end{equation}
This position comparison holds for every realization of the urn.
An endpoint whose label was in $E$ before its disclosure causes failure.
Its high-label probability
is at most $Ct\bar g_n/k_n^\circ$: exposed coordinates satisfy their
caps, while unexposed coordinates have posterior means at most
$C(1+s_j)/j$ and $\sum_j s_j\le mt$. A low endpoint requires failure
of the reference age requirement.

\paragraph{Right children.}
Let $\mathcal D_i=[(i-1)/n,i/n)$ and put
\[
 \lambda_i=\int_{\mathcal D_i\cap[A_v,1]}\rho_v(s)\,ds,
 \qquad \sum_i\lambda_i
 =\Gamma_v\bigl(A_v^{-(1-\chi)}-1\bigr)\le\Lambda_n.
\]
Condition on the current record $\mathcal H$, fixing the parent data.
The current mark is coupled and good. Since
$p_{i\to k}=\psi_kT_{k,i}$, \eqref{eq:conditional-urn-products},
the beta--Gamma estimate, and Cauchy--Schwarz give
\begin{equation}
 \sum_{i>k}\Ebb[|mp_{i\to k}-\lambda_i|\mid\mathcal H]
       +\sum_{i\le k}\lambda_i
 \le C\Lambda_n\left(
   \Delta_v+\varepsilon_n^{\mathrm{BG}}
   +\sqrt{\eta_n(t)}+\frac1{k_n^\circ}\right).
\label{eq:right-intensity-L1}
\end{equation}
The Riemann-sum error is at most $Cg_n/k_n^\circ$;
replacing $k/n$ by $A_v$ in the intensity costs
$C\Lambda_n\Delta_v$, including the interval between the two ages.

For a slot with sender $i>k$ whose endpoint has not been revealed,
let $F$ be its excluded target set and suppose $k\notin F$.
Conditional on the urn and $\mathcal A$, the probability
that it attaches to $k$ is $\widetilde p=p_{i\to k}/(1-q)$, where
$q=\sum_{\ell\in F}p_{i\to\ell}<1/2$ by
\eqref{eq:left-forbidden-bound}. The mixed product moments above give
\[
 \Ebb[p_{i\to k}p_{i\to\ell}\mid\mathcal A]
 \le Cg_n\bar g_n k^{\chi-1}\ell^{\chi-1}i^{-2\chi},
 \qquad \ell\in F,\quad \ell\ge k_n^\circ.
\]
For the low targets, the same moments give
\[
 \sum_{i>k}\Ebb\left[
  p_{i\to k}\sum_{\ell\in E_0}p_{i\to\ell}
  \,\middle|\,\mathcal A\right]
 \le Cg_n\left(\sum_{\ell\in E_0}\vartheta_\ell\right)
       k^{\chi-1}(k_n^\circ)^\chi\sum_{i>k}i^{-2\chi}
 \le\frac{Ct g_n\bar g_nH_n}{k_n^\circ}.
\]
Consequently the total correction from exclusions is at most
\begin{equation}
 \sum_{i>k}\sum_{j=1}^m
 \Ebb[|\widetilde p_{i,j}-p_{i\to k}|\mid\mathcal A]
 \le\frac{Ct g_n\bar g_nH_n}{k_n^\circ}.
\label{eq:right-conditioning-error}
\end{equation}
Here and below only slots with unrevealed endpoints and $k\notin F$
enter these sums.
For $\chi>1/2$, the sum of $i^{-2\chi}$ above
$\max\{k,\ell\}$ is at most $C\max\{k,\ell\}^{1-2\chi}$;
for $\chi=1/2$, it contributes $H_n$.
The same calculation bounds the conditional expected sum of squared
attachment probabilities, together with $\sum_i\lambda_i^2$, by
$Cg_n^2H_n/k_n^\circ$.
Edges to $k$ from senders with labels in $E$, apart from the entrance, and
reference points in cells indexed by $E$ cost at most $Ct g_n/k_n^\circ$.

Scan the cells in increasing order. In each cell, retain the indicators
that eligible PA slots attach to $k$, and couple their sum to a
$\operatorname{Poisson}(\lambda_i)$ count, updating exclusions and
registering their senders whenever an attachment is found. Declare
failure if the counts differ or exceed one, if an edge to $k$ other
than the entrance has sender in $E$, or if a
reference point has $i\le k$ or $i\in E$; membership in $E$ is tested
before the query.

Average the urn-conditional coupling of the indicators and reference
count over the current urn posterior, then sample the reference count
conditional on the retained indicators using independent auxiliary randomness.
This preserves the posterior identity and the prescribed Poisson
law given the preceding full record. Le Cam's bound and the preceding
estimates bound the count and collision errors throughout the block by
\begin{equation}
 C\left[
  \Lambda_n\left(\Delta_v+\varepsilon_n^{\mathrm{BG}}
          +\sqrt{\eta_n(M)}+\frac1{k_n^\circ}\right)
  +\frac{Mg_n\bar g_nH_n}{k_n^\circ}
 \right].
\label{eq:right-block-error}
\end{equation}
The estimates hold at each intervening record with $t\le M$.
Indeed, the conditional error in replacing $mp_{i\to k}$ by
$(1-\chi)\Gamma_vk^{\chi-1}i^{-\chi}$ is at most
$C\Gamma_vk^{\chi-1}i^{-\chi}
(\varepsilon_n^{\mathrm{BG}}+\sqrt{\eta_n(M)}+1/k_n^\circ)$.
These deterministic intensity envelopes sum to at most $C\Lambda_n$.
For exclusions, one may bound
$\sum_{\ell\in F,\,\ell\ge k_n^\circ}\ell^{\chi-1}
\le M(k_n^\circ)^{\chi-1}$ and use the fixed low-target cap
before summing over $i$; labels newly added to $E$ belong to cells
already scanned. Thus the same bound applies as the record grows.

For each reference point, draw its age in $\mathcal D_i\cap[A_v,1]$ with
density $\rho_v/\lambda_i$, still without exposing child marks.
In a cell with exactly one PA slot attaching to $k$ and one reference
point, set $K_w=i$ and retain that slot as the entrance.
The reference positions have exactly
the prescribed Poisson-process law, and matched pairs satisfy
\begin{equation}
 \Delta_w\le\frac1{nA_w}\le\frac{C}{k_n^\circ}.
\label{eq:right-age-label-error}
\end{equation}
After successful position comparisons, all outgoing slots of $k$ have
revealed endpoints, and every slot with a younger sender either has a
revealed endpoint or has been tested for attachment to $k$.
Thus the complete incident-edge multiset at $k$ is matched.
The vertex $k$ is now explored, while the complete stars of its newly
registered children have not yet been recorded.

\paragraph{New marks.}
Condition on the completed left and right geometry and all retained
PA answers of the star. On a successful position comparison, a newly
registered child with label $j$ has $s_j=1$ if it is a left child and
$s_j=0$ if it is a right child. Write $\iota$ for this value; an initial root without a
recorded incidence also has $\iota=0$. For
$Q=\operatorname{Beta}(a+\iota,d_j)$,
\[
 \Ebb_Q|h_j(X)-1|
 \le N_j\Ebb_QX+\Ebb_Q(e^{\lambda_jX}-1)\le CL/j,
 \qquad \Ebb_Qh_j(X)\ge1-CL/j.
\]
Here the beta moments give
$\Ebb_Qe^{\lambda_jX}\le
(1-\lambda_j/(a+\iota+d_j))^{-(a+\iota)}$.
Thus the posterior is within $CL/j$ of $Q$ in total variation.
Under $f_j^{-1}$, $Q$ becomes $\mu_a$ when $\iota=0$; when $\iota=1$,
its density relative to $\mu_a$ is $c_jf_j(z)/a$, whereas that of
$\mu_{a+1}$ is $z/a$. Equation~\eqref{eq:beta-gamma-L1} bounds this
additional error by $C/j$. These bounds are uniform in the other urn
variables and hence remain valid after integrating them out.

Before exposing each $Z_j$, let $\mathcal H^-$ be the full record,
including the completed positions and all earlier mark comparisons.
Maximally couple $\mathcal L(Z_j\mid\mathcal H^-)$ to $\mu_{a+\iota}$.
Implement this coupling by revealing and retaining the actual $Z_j$,
then drawing $\Gamma_w$ from its conditional coupling kernel using
independent auxiliary randomness. Thus
$\mathcal L(\Gamma_w\mid\mathcal H^-)=\mu_{a+\iota}$. Successive reference
marks therefore have the prescribed product law conditional on the positions.
After failure, keep generated values and apply these same target laws
to remaining marks under the completion convention above. For
$\mathcal G=\mathcal G_j^{\mathrm{BG}}(\varepsilon_n^{\mathrm{BG}})$,
$\mu_{a+1}(\mathcal G^c\cap(0,g_n])\le
(g_n/a)\varepsilon_n^{\mathrm{BG}}$, so
\begin{equation}
 \Pbb\left(
   \bigl(\{Z_j\ne\Gamma_w\}\cup\{Z_j\notin\mathcal G\}\bigr)
   \cap\{\Gamma_w\le g_n\}
   \,\middle|\,\mathcal H^-\right)
 \le C\left[
   g_n\varepsilon_n^{\mathrm{BG}}
   +\frac{t+1}{k_n^\circ}
   +\frac{nt\bar g_n}{(k_n^\circ)^2}
 \right].
 \label{eq:joint-mark-error}
\end{equation}
On mark agreement and goodness, the definition of
$\mathcal G_j^{\mathrm{BG}}$ gives
$j\psi_j\le Cg_n<\bar g_n$. Thus the new target cap costs no further
exceptional event.

\paragraph{Errors within a trial.}
Every run of left edges starts at $v_0$ or a right child.
Equations~\eqref{eq:left-error-recursion} and
\eqref{eq:right-age-label-error} give
$\max_w\Delta_w\le\Delta_{v_0}+Cr_n\eta_n(M)<1/4$,
proving \eqref{eq:extension-position-error} and justifying the
comparisons. Stop before a query could exceed the label budget.
On matching records, an overflow is a PPT volume failure. The PPT
age, mark, and volume requirements cost at most
$Ce^{-(C_{\mathrm{test}}+2)\ell_n^2}$. There are at most $mB_n$ left
slots, $B_n$ right blocks, and $B_n$ new marks. Summing their bounds
therefore gives, for the trial failure $F$,
\begin{equation}
 \Pbb(F\mid\mathcal H)
 \le Ce^{-(C_{\mathrm{test}}+2)\ell_n^2}
 +CB_n\left[
   \Lambda_n\left(\Delta_{v_0}+r_n\eta_n(M)
       +\varepsilon_n^{\mathrm{BG}}+\sqrt{\eta_n(M)}\right)
   +\frac{Mg_n\bar g_nH_n}{k_n^\circ}
 \right].
\label{eq:extension-error}
\end{equation}
For $M=M_n$,
\[
 \sqrt{\eta_n(M_n)}
 =e^{-\frac12c_M\sqrt{\log n}+O(\ell_n^2)},\qquad
 B_n\Lambda_n\varepsilon_n^{\mathrm{BG}}=n^{-a/(2a+1)+o(1)}.
\]
Since $r_n\eta_n(M_n)=o(\sqrt{\eta_n(M_n)})$, the root hypothesis
makes the right-hand side at most
$e^{-(C_{\mathrm{test}}+1)\ell_n^2}$ for large $n$.
\end{proof}

\paragraph{Repeated applications.}
Apply Lemma~\ref{lem:ext-coup} at at most
$N\le e^{C_{\mathrm{test}}\ell_n^2}$ adaptively chosen roots, stopping
when fewer than $B_n$ labels remain in the budget. Initiation is
measurable from the preceding record $\mathcal H_i$, which retains
all earlier failures. With $F_i=\varnothing$ for an uninitiated trial,
\begin{equation}
 \Pbb\left(\bigcup_{i\le N}F_i\right)
 \le\sum_{i\le N}\Ebb\left[
   \mathbf1_{\{i\text{ initiated}\}}\Pbb(F_i\mid\mathcal H_i)\right]
 \le Ne^{-(C_{\mathrm{test}}+1)\ell_n^2}
 \le e^{-\ell_n^2}=o(1).
\label{eq:repeated-extension-error}
\end{equation}
To initialize from an empty record, fix
$k_0\in\{\lceil nb_n\rceil,\ldots,n\}$, set $K_{v_0}=k_0$,
$A_{v_0}=k_0/n$, and $\Gamma_{v_0}=Z_{k_0}$.
This mark has law $\mu_a$ and $\Delta_{v_0}=0$; its tail and
beta--Gamma goodness errors add $o(1)$. Later roots retain their
previously assigned data. Along at most $Nr_n$ left steps,
\begin{equation}
 \Delta_v\le C\bigl(\varepsilon_n+Nr_n\eta_n(M_n)\bigr)
       =o(\sqrt{\eta_n(M_n)}),
\label{eq:coupling-depth-error}
\end{equation}
where the $\varepsilon_n$ term also allows an initial discrepancy
$O(\varepsilon_n)$. Thus the root hypothesis persists throughout the
sequence. No urn exceptional probability is used in these bounds.

\subsection{Coupling the walks of mixing time order}
\label{sec:coupling-after-burn-in}

Fix $C_0<\infty$. Let $\mathcal H$ be admissible with $|E|\le M_n/2$,
and let $V_0$, matched to a PPT root $\rho$, satisfy the root hypotheses
of Lemma~\ref{lem:ext-coup}.

Let $b$ be the parent at the entrance, and put
$\sigma_i=\inf\{t:|X_t|=r_n+i\}$ for $i\ge0$ and
$\tau_b=\inf\{t:X_t=b\}$,
with $\tau_b=\infty$ if there is no entrance. Use $\ell_0,\ell_1$
from Proposition~\ref{prop:no-truncation}, and write
$\tau=\inf\{t:|X_t|=\ell_1\}$.
Time and depth are measured from $V_0$ and $\rho$ throughout;
when an entrance is present, root degrees include its edge to $b$.
Choose the constants so that $h_n=r_n$,
$a_n\ge b_n$, $V_n\le B_n$, $c_M<1/2$, and
$(\log n)\varrho^{h_n}=o(1)$.

We construct a subgraph $\widetilde\Gamma\subseteq G_n$ and an injective
map $\Phi:\widetilde\Gamma\to\widetilde{\mathcal T}$. First compare the
radius-$(r_n+1)$ ball by the proof of Lemma~\ref{lem:ext-coup}, setting
$\widetilde\Gamma=\varnothing$ on failure. On success, apply the PPT
truncation tests to the initial ball. Retain in $\widetilde\Gamma$
the vertices of depth at most $r_n$ whose ancestral lines, including
the vertices themselves, pass these tests. Record the initial matching
in $\Phi'$, with $\Phi'(K_v)=v$.
Queue the depth-$(r_n+1)$ pairs whose parents are retained, in a fixed
breadth-first order, retaining their ages, marks, and entrance slots.

For a queued pair $(k,z)$, declare the pair truncated if
$\mathcal B_{\rm tr}(z)$ occurs. Otherwise, if an additional incidence at $k$
has already been recorded, declare the pair \emph{mismatched}.
Otherwise compare the remaining PA incident slots with the PPT
offspring of $z$ as in Lemma~\ref{lem:ext-coup}, stopping the comparison
once $B_n$ labels have been added to $E$ in this comparison.
The coupling gives a type-preserving bijection and every new pair $(j,y)$ must
satisfy $A_y\ge8k_n^\circ/n$, $Z_j=\Gamma_y\le g_n$, and
$Z_j\in\mathcal G_j^{\mathrm{BG}}(\varepsilon_n^{\mathrm{BG}})$.
On success, add $k$ to $\widetilde\Gamma$, extend $\Phi'$ by the new
pairs, and queue its children;
on failure, declare the pair mismatched. If an exposed value satisfies
$j\psi_j>\bar g_n$, stop the entire construction and set
$\widetilde\Gamma=\varnothing$.
Do not continue PA--PPT comparisons below a truncated or mismatched pair,
or compare stars at depth $\ell_1+1$.
Thus complete stars are compared through depth $\ell_1$ along the remaining branches.
For the finite-depth tests in \eqref{eq:ppt-truncation-rule}, use the
already generated data through the current level in active branches;
complete missing data only in truncated or mismatched branches, without
PA queries. No future offspring of a queued pair are thereby revealed.
Once this continuation ends, its reference tree is not used for further
PA comparisons and may be completed with its conditional law.

Truncated and mismatched pairs
remain in $\Phi'$,
but only explored vertices whose complete stars were successfully compared enter
$\widetilde\Gamma$. Keep the edges between these retained vertices
and put $\Phi=\Phi'|_{\widetilde\Gamma}$. The matching is injective, the matched labels are at least
$k_n^\circ$, and
\[
 \deg_{G_n}(x)=\deg_{\mathcal T}(\Phi(x)),
 \qquad x\in\widetilde\Gamma.
\]

Defining
\[
 \Pi_t:=\{X_s\in\Phi(\widetilde\Gamma)\text{ for all }s\le t\},
\]
we have coupling of SRW
$\mathscr X$ on $G_n$ from $V_0$ with SRW $X$ on the PPT from $\rho$
until time $t$ on $\Pi_t$.

\begin{prop}
\label{prop:walk-coupling-log-time}
Fix $\gamma>0$. Uniformly over the admissible histories and root data,
including records of earlier failed comparisons, the construction has
the prescribed conditional PA and PPT laws and satisfies
\begin{equation}\label{eq:walk-target-coupling}
\begin{split}
 \Pbb\left(\Pi_\tau^c
       \,\middle|\,\sigma_0<\tau_b,\mathcal H\right)
 \le \Pbb\left(\bigcup_{s\le\tau}\mathcal B_{\rm tr}(X_s)
       \,\middle|\,\sigma_0<\tau_b,\mathcal H\right)
 +o(1).
\end{split}
\end{equation}
\end{prop}

\begin{proof}
Lemma~\ref{lem:ext-coup} bounds failure of the initial comparison
by $o(1)$, uniformly in $\mathcal H$.
We first check the bound on $|E|$.
By \eqref{eq:truncated-exploration-size}, at most
$(t_n+1)e^{\omega_n}V_n$ stars are compared. Each comparison adds at most
$B_n$ labels to $E$, even on failure. Thus the number of labels added
to $E$ after the initial ball is at most
\begin{equation}
 CB_n(t_n+1)e^{\omega_n}V_n
 =n\exp\left\{-\tfrac12\sqrt{\log n}+O(\ell_n^2)\right\}=o(M_n),
\label{eq:walk-registration-budget}
\end{equation}
since $c_M<1/2$. Thus the initial record leaves enough room for the
construction. Before the first violation of $j\psi_j\le\bar g_n$, each newly exposed
mark has only its entrance recorded. The posterior stochastic bound
$\operatorname{Beta}(a+\iota,d_j-\lambda_j)$, $\iota\in\{0,1\}$,
and its exponential moments give, uniformly in the current record,
\[
 \Pbb(j\psi_j>\bar g_n\mid\mathcal H_{\rm current})
 \le Ce^{-c\bar g_n}.
\]
At most $M_n$ marks are exposed. Increasing $C_G$ therefore bounds
the probability of any such global stop by $CM_ne^{-c\bar g_n}=O(n^{-8})$.

Apart from initialization failure, a cap stop, and an entrance return, the first
exit by time $\tau$ visits a truncated or mismatched vertex. We bound
the latter possibility by exposing the walk as in \cite{RWRG}, while
keeping the fixed breadth-first order and coupling kernels.
Before exploring a later level $r_n+i$, $i\ge1$, reveal the walk
through $\sigma_i$, stopping earlier at an entrance return or a
previously detected defect.
This path uses only complete stars at smaller depths and independent
walk randomness, so it gives no further information about unqueried
PA variables or unrevealed PPT offspring.

For $\sigma_i\le\tau$, let $u_i$ be the ancestor of $X_{\sigma_i}$
at distance $r_n-1$. Let $F_i$ be the reached matched candidates in
$\partial B_{r_n-1}^+(u_i)$, including candidates subsequently truncated
or mismatched, and let $L_i$ be the children reached when the stars
of vertices in $F_i$ are successfully compared. At $\sigma_0$, take $L_0$
to be the initially queued children of $F_0$; their parents' star
comparisons are covered by the initialization estimate.
For $i\ge1$, $F_i$ is selected before its comparisons,
and each new child is added to $L_i$ upon completion of its parent's
star comparison.
Stop the analysis at a visit to a truncated vertex, a return to a previously
selected $u_j$ after $\sigma_j$, or the first charged failure.
The volume bound at $X_{\sigma_i}$ gives $|F_i|\le V_n\le B_n$;
at most $B_n$ labels are added per star comparison, so $|L_i|\le CB_n^2$.
The truncation checks use PPT data only through the current level, so the
offspring of candidates about to be compared remain unrevealed.

Unless the walk returns to $u_i$ after $\sigma_i$, its first visit
to depth $r_n+i+1$ lies below $u_i$, and $u_{i+1}$ is a child of $u_i$.
Consequently $F_{i+1}\subseteq L_i$ until a truncated or mismatched vertex
is visited. By \eqref{eq:ancestor-return}, all such ancestor returns
have total probability at most
$(\ell_1+1)\varrho^{r_n-1}=o(1)$.
Thus it suffices to control comparisons in the selected lists and
additional incidences recorded at their children before their complete
stars are recorded.

Track each vertex in $L_i$ from $\sigma_0$ if $i=0$, and from completion
of its parent's star comparison otherwise, until its own comparison begins.
Its mark is already exposed and capped when tracking begins. Within the
parent's comparison, a label already
in $E$ before the query or two or more attachments from the same sender
to the current target already cause failure; the
whole left--right--mark comparison finishes before another star's
queries interleave.
Subsequent contacts are charged when first
disclosed, using membership in the current list. For a slot $(u,s)$ whose
endpoint has not been revealed and an exposed target $j$, the conditional query law gives
\[
 \Pbb(J_{u,s}=j\mid\psi,\text{current record})
 \le \frac{\bar g_n/k_n^\circ}{1-2M_n\bar g_n/k_n^\circ}
 \le \frac{2\bar g_n}{k_n^\circ}.
\]
The endpoints of at most $mM_n$ outgoing slots are revealed while
exploring their senders; their contacts with $L_i$ cost at most
$C|L_i|M_n\bar g_n/k_n^\circ$. Conversely, each listed sender has
at most $m$ slots, each tested at most once against any given target
in $E$, giving the same bound for contacts from listed senders to other
targets. Hence the total contact cost per list is at most
$CB_n^2M_n\bar g_n/k_n^\circ$. This includes queries in failed
off-list comparisons. Before the first charged failure, every
selected root therefore has only its entrance recorded, as required
by Lemma~\ref{lem:ext-coup}.

Every star compared after initialization has root age at least $a_n$
and mark at most $g_n$. Its remaining offspring contain a left child below $8k_n^\circ/n$
with probability at most
$Ce^{-\chi(c_1-C_A)\ell_n^2}$; right children cannot cross this
cutoff. The probability of a child mark exceeding $g_n$ is at most
$C(m+g_ne^{(1-\chi)C_A\ell_n^2})e^{-g_n/2}$.
Choose $c_1$ and then $C_g$, after $C_B$, so that both bounds are
$o(B_n^{-2}e^{-2\ell_n^2})$.
The number of right children is Poisson with mean at most
$g_ne^{(1-\chi)C_A\ell_n^2}=o(B_n)$, so its probability of exceeding
$B_n/2$ is smaller still. Thus the summed age, mark, and label-budget
overflow errors per frontier are at most $Ce^{-2\ell_n^2}$.

The left-edge recursion and the right-child bound give
\[
 \Delta_v\le\Delta_\rho+C(\ell_1+1)\eta_n(M_n)
          \le(C_0+1)\sqrt{\eta_n(M_n)}
\]
for all sufficiently large $n$, since $\ell_1=O(\log n)$.
Each generation has at most $B_n$ selected stars and $CB_n^2$
slot, block, and mark comparisons. Summing the conditional estimates
from Lemma~\ref{lem:ext-coup} and adding the contact bound gives error
at most $Cq_n$ per frontier, where
\[
 q_n:=e^{-2\ell_n^2}
 +B_n^2\left[
   \Lambda_n\bigl(\varepsilon_n^{\mathrm{BG}}+\sqrt{\eta_n(M_n)}\bigr)
   +\frac{M_ng_n\bar g_nH_n}{k_n^\circ}
 \right].
\]
All charges use the record before the relevant query or comparison;
future comparison success and future nonreturn are not added to
the conditioning. There are at most
$\ell_1+1$ selected generations. The terms containing $M_n$ are at most
$\exp\{-\tfrac12c_M\sqrt{\log n}+O(\ell_n^2)\}$, and the remaining
comparison term is $n^{-a/(2a+1)+o(1)}$. Hence $(\ell_1+1)q_n=o(1)$, and
summing over the selected frontiers gives
\begin{equation}\label{eq:walk-level-check}
\Pbb(\text{a cap stop or a selected comparison or contact failure}
     \mid\mathcal H)
\le C(\ell_1+1)q_n+O(n^{-8})=o(1).
\end{equation}

An entrance return after $\sigma_0$ costs at most $\varrho^{r_n}$.
Since $\Pbb(\sigma_0<\tau_b\mid\mathcal H)\ge p_{\rm esc}$,
conditioning on $\sigma_0<\tau_b$ multiplies the preceding error
bounds by at most $p_{\rm esc}^{-1}$. Thus
$\Pi_{\tau}^c$ requires a visit to a truncated vertex apart from
an error $o(1)$, uniformly in $\mathcal H$.
\end{proof}

We now define a continuation from one entrance. Fix an integer
testing depth $R\le r_n$ with $R\to\infty$, keeping the explored
region $\widetilde\Gamma$ constructed above. Run $\mathscr X$ until
its first entrance return, first other exit from $\widetilde\Gamma$, or
first visit to a retained vertex at depth $R$. These outcomes are called
\emph{return}, \emph{defect}, and \emph{penetration}, respectively.
A defect or penetration is a \emph{resolution}. Write $\xi$ for
the duration of this test. If $\widetilde\Gamma=\varnothing$,
declare an immediate defect and put $\xi=0$. After penetration
within $R^2$ steps, continue until the first exit from
$\widetilde\Gamma$, the first visit to matched depth $\ell_1$, or
time $\lceil t_n^\star\rceil$ measured from $V_0$. On the coupled
event, the visit to depth $\ell_1$ occurs at time $\tau$.

Let $\tau_{\rm exit}$ be the first exit from $\widetilde\Gamma$, including an
entrance return, and for a retained vertex $y$ put
\[
 h_{\widetilde\Gamma}(y)
   =P_{V_0}^{G_n}(\tau_y<\tau_{\rm exit}).
\]
This probability is determined by the recorded incident slots. Fix
$C_{\rm cap}>2\gamma$ and put
$\kappa_n=n^{-1}e^{C_{\rm cap}\sqrt{\log n}}$.
The continuation is \emph{accepted} if the test penetrates within
$R^2$ steps, the walk reaches depth $\ell_1$ before exiting and
by time $\lceil t_n^\star\rceil$, and its depth-$\ell_0$ ancestor $y$
satisfies $h_{\widetilde\Gamma}(y)\le\kappa_n/p_{\rm esc}$.
All these conditions are measurable from the finite exploration
and the walk at its stopping time.

For a fixed record, let $r$ be the resolution probability and put
\[
 q(z)=P_{V_0}^{G_n}(\text{continuation accepted at }z).
\]
Thus $r-q(V(G_n))$ is the probability of a defect or a penetration
with a rejected continuation.

\begin{lem}[Continuation from one entrance]
\label{lem:local-continuation}
For every realization of the construction,
\[
 r\ge p_{\rm esc},\qquad E_{V_0}^{G_n}\xi\le C_mR,\qquad
 q(z)\le\frac{\kappa_n}{p_{\rm esc}^{\,2}}r
 \quad(z\in V(G_n)).
\]
Fix $A,C>0$ and $\varepsilon\in(0,1)$. In addition to the preceding
coupling hypotheses, suppose that, conditional on $\mathcal H$, the
unexposed forward PPT has its prescribed-root law, with
$A_\rho\ge e^{-A\ell_n^2}$ and $\Gamma_\rho\le C\ell_n^2$.
Choose $\gamma$, then $C_{\rm cap}$ sufficiently large, and the
truncation and coupling constants with $A<C_A\le c_0$ and $C_g\ge C$.
Uniformly over these histories and root data, for all sufficiently
large $n$,
\begin{equation}\label{eq:local-continuation-loss}
 \Ebb\left[r-q(V(G_n))
       \,\middle|\,\mathcal H\right]\le\varepsilon^2p_{\rm esc}.
\end{equation}
Consequently, the event
$\mathcal E=\{q(V(G_n))\ge(1-\varepsilon)r\}$ is measurable from
the finite exploration record and satisfies
\[
 \Pbb(\mathcal E^c\mid\mathcal H)
 \le\varepsilon.
\]
\end{lem}

\begin{proof}
An immediate defect has $r=1$, $\xi=0$, and $q=0$. Otherwise compare
the test with the completed PPT walk. Never hitting the entrance
forces a resolution, so Lemma~\ref{lem:escape} gives
$r\ge p_{\rm esc}$. The type-corrected depth used in
\eqref{eq:radial-exponential-drift} has drift at least $1/6$ before
the test stops and is at most $R+1/2$ at that time. Optional
stopping gives $E_{V_0}^{G_n}\xi\le C_mR$.

If $q(z)>0$, the endpoint has a unique depth-$\ell_0$ ancestor $y$
in the matched tree. Every accepted path ending at $z$ hits $y$
before exiting. Hence
\[
 q(z)\le h_{\widetilde\Gamma}(y)
       \le\frac{\kappa_n}{p_{\rm esc}}
       \le\frac{\kappa_n}{p_{\rm esc}^{\,2}}r.
\]
For the probability estimate, let $\mu$ be the law of the limiting
ray in the completed forward PPT, conditioned never to hit the
entrance parent. By the coupling and uniform escape,
\[
 h_{\widetilde\Gamma}(y)
 \le P_\rho^{\mathcal T}(\tau_{\Phi(y)}<\tau_b)
 \le\frac{\mu(\Phi(y))}{p_{\rm esc}}.
\]
Indeed, after hitting $\Phi(y)$, the probability of never returning
to its parent is at least $p_{\rm esc}$. Thus the condition
$\mu(\Phi(y))\le\kappa_n$ suffices for the hitting-probability
condition in the acceptance rule.

Let $E_s$ be the event that the completed PPT walk reaches depth
$s$ before the entrance parent. Since $R\le r_n$, the strong
Markov property and \eqref{eq:ancestor-return} give
\[
 E_{r_n}\subseteq E_R,\qquad
 \Pbb(E_R\setminus E_{r_n}\mid E_R,\mathcal H)
 \le\varrho^{R+1}=o(1).
\]
Thus the coupling estimate conditioned on $E_{r_n}$ also holds
conditioned on $E_R$, with an additional error $\varrho^{R+1}$.
The duration bound discards $\{\xi>R^2\}$ at cost $C_m/R$.
For a defect before penetration, the completed PPT walk reaches
depth $r_n$ from the defect location before the entrance with
probability at least $p_{\rm esc}$. Such a defect forces
$\Pi_\tau^c$, so its probability is at most $p_{\rm esc}^{-1}$
times the coupling-failure probability conditioned on $E_{r_n}$.
For an immediate defect, use the PPT root as the defect location.

Proposition~\ref{prop:walk-coupling-log-time} bounds the coupling
loss by the truncation loss through $\tau$ and
$o(1)$. The stopped-time estimate \eqref{eq:ppt-stopped-no-truncation}
bounds this loss and the deadline failure
$\{\tau>\lceil t_n^\star\rceil\}$ by $\eta(\gamma)+o(1)$,
uniformly over the prescribed root data.
Apply Proposition~\ref{thm:global-age-volume} with relative age
parameter $C_A-A>0$. Its age bound gives
$A_\rho e^{-(C_A-A)\ell_n^2}\ge a_n\ge b_n$.

At time $\tau$, the depth-$\ell_0$ ancestor of $X_\tau$ can differ
from the corresponding vertex of the limiting ray only if the walk
subsequently hits that ancestor's parent. By \eqref{eq:ancestor-return},
this has probability at most
$C\varrho^{\ell_1-\ell_0+1}=o(1)$, also under nonreturn conditioning.
Outside this event, its harmonic mass is $e^{-W_{\ell_0}}$, and
\[
 d\ell_0-\log\kappa_n^{-1}
   =(C_{\rm cap}-2\gamma)\sqrt{\log n}+O(1).
\]
Theorem~\ref{thm:original-W-tightness} makes failure of the
harmonic-mass bound arbitrarily small by increasing
$C_{\rm cap}-2\gamma$. The proof of Proposition~\ref{prop:no-truncation}
gives the same estimate for the prescribed root data: the initial
harmonic contribution is $O_{\Pbb}(\ell_n^2)=o_{\Pbb}(\sqrt{\log n})$.
Replacing nonreturn conditioning by $E_R$ changes the error by
at most $\varrho^{R+1}$.

Together these estimates bound the conditional expected rejected
mass in \eqref{eq:local-continuation-loss}. Choose $\gamma$, then
$C_{\rm cap}$, and finally $n$ to make their sum at most
$\varepsilon^2p_{\rm esc}$.
The event $\mathcal E$ is measurable
because the killed transition probabilities and acceptance rule
are determined by the finite record. Its probability bound follows
from Markov's inequality and $r\ge p_{\rm esc}$.
\end{proof}

The conditional loss bounds hold at every admissible retained record.
Section~\ref{sec:worst-start} verifies the conditions on the initial
record and the budget bounds on a common event, and applies these
conditional loss bounds at the initiated couplings.
The common exceptional probability is excluded once, after combining
the estimates over starting vertices.

\section{Upper bound on the worst-case mixing time}\label{sec:worst-start}

We extend the good-start analysis of
Proposition~\ref{prop:walk-coupling-log-time} to arbitrary initial states.
For each exploration below, we use the registered label set $E$
from Section~\ref{sec:coupling}, relative to the current record.

\subsection{Good roots and common events}
\label{subsec:stopped-geometry}

Throughout this section, write $\beta=1-\chi$, fix $A>0$ and a
sufficiently large $c_R$, and set
\[
 L=\log\log n,\qquad R=\lceil c_RL\rceil,\qquad
 D=\lceil(\log n)^3\rceil,\qquad a=\lfloor ne^{-AL^2}\rfloor+1.
\]
Labels below $a$ are called old; the remaining labels are called young.
The set of hubs is $\mathcal C_n=\{v<a:\deg_{G_n}(v)>D\}$.
For $r\ge0$, form the stopped neighborhood $B_r^{\mathrm{stop}}(x)$
by breadth-first search from $x$, retaining every edge incident to
a nonhub vertex at depth less than $r$. The search stops at hubs
and at depth $r$, retaining these vertices.
If $x\in\mathcal C_n$, the neighborhood is $\{x\}$.
Call $x$ a \emph{good root} if $B_R^{\mathrm{stop}}(x)$ is a tree.
In particular, every hub is good.

SRW reaches a good root quickly because each larger stopped
neighborhood has at most one cycle. For a connected multigraph $H$, write
$\operatorname{tx}(H)=|E(H)|-|V(H)|+1$ for its tree excess, counting
edges with multiplicity.

\begin{lem}\label{lem:stopped-tangle-free}\label{lem:cycle-clearing}
With probability $1-o(1)$,
\begin{equation}\label{eq:stopped-geometry}
 \{1,2\}\subseteq\mathcal C_n,\qquad
     \operatorname{tx}\bigl(B_{6R}^{\mathrm{stop}}(x)\bigr)\le1
     \quad\text{for every }x\in[n].
\end{equation}
For every realization $G$ of $G_n$ satisfying \eqref{eq:stopped-geometry},
\[
 \sup_{x\in V(G)}P_x^G\left(
   X_{\tau_{\mathcal C_n}\wedge5R}\text{ is not a good root}
 \right)\le e^{-cR},
\]
where $\tau_{\mathcal C_n}$ is the first hitting time of $\mathcal C_n$
and $c>0$ depends only on the model parameters.
\end{lem}

\begin{proof}
We first prove \eqref{eq:stopped-geometry}. Put
$k=\lfloor n^{\beta/2}\rfloor$, so $k<a$ for large $n$.
The attractiveness of a vertex is its current degree plus $\delta$.
Before the $j$-th edge of vertex $t$ is attached, the sum of these
weights over vertices $v<t$ is
\[
 T=2m(t-2)+\delta(t-1)+j-1.
\]
A vertex of attractiveness $d$ is selected with probability $d/T$.
The reciprocals
$1/T$ sum to $\beta\log(t/i)+O(1)$ over slots between birth times
$i$ and $t$, with seed times interpreted as $2$.
For a vertex of attractiveness $d\ge m+\delta\ge2$, the weight
$1/(d-1)$ has conditional multiplier $1-1/T$. Iterating from birth
and applying Markov's inequality gives
\[
 \mathbb E\frac1{\deg_{G_n}(i)+\delta-1}\le C(i/n)^\beta,
 \qquad
 \mathbb P\bigl(\deg_{G_n}(i)\le D\bigr)\le CD(i/n)^\beta.
\]
Consequently,
\[
 \mathbb P\bigl([k]\nsubseteq\mathcal C_n\bigr)
 \le CDn^{-\beta}k^{1+\beta}
 \le CDn^{-\beta\chi/2}=o(1).
\]

We use the edge-set estimate of
Garavaglia--Stegehuis~\cite[Lemma~1]{SubgraphPA} in a form uniform
over graphs of maximum degree at most four. For such a directed
multigraph $F$, with distinct birth labels $i_v$, edges directed
from younger to older vertices, and no deterministic seed edges,
\begin{equation}\label{eq:embedding-bound}
 \mathbb P(F\subseteq G_n)
 \le C^{|E(F)|}\prod_{u\to v\in E(F)}i_u^{-\chi}i_v^{-\beta}
 \le C^{|E(F)|}\prod_{v\in V(F)}i_v^{-\deg_F(v)/2}.
\end{equation}
The last inequality uses $\chi\ge1/2$ and $i_u>i_v$.
Since the cited statement fixes the number of edges, we verify the
dependence of the constant. Fix an injective assignment of the
prescribed outgoing edges to birth slots. After each vertex's birth,
give its attractiveness $d$ the weight
$f_s(d)=d(d+1)\cdots(d+s-1)$, with $f_0=1$, where $s\le4$ is the
number of its prescribed arrivals still pending. The identities
\[
 \frac{f_s(d+1)}{f_s(d)}=1+\frac{s}{d},
 \qquad d f_{s-1}(d+1)=f_s(d)
\]
show that the product of weights, killed at the first failed
prescription, has conditional multiplier $1+\sum_v s_v/T$ at
an unprescribed slot and $1/T$ at a prescribed slot; the sum is over
eligible tracked vertices. Its initial factors cost at most
$C^{|E(F)|}$, and its final value is the indicator of the prescribed
attachments. Using $1+x\le e^x$ and $s_v\ge0$, we may include all
slots in the exponent sum. Each pending edge $u\to v$ then contributes
at most $\beta\log(i_u/i_v)+O(1)$, while its prescribed slot contributes
at most $C/i_u$. Partial birth blocks and seed times change only
the constant per edge. Summing over at most $m^{|E(F)|}$ slot
assignments proves \eqref{eq:embedding-bound}.

Suppose $[k]\subseteq\mathcal C_n$ and some
$B_{6R}^{\mathrm{stop}}(x)$ has excess at least two.
Take two non-tree edges and the paths in the breadth-first tree from
their endpoints to $x$. Their union has excess two and at most $24R+2$
edges. Deleting leaves gives a connected witness $F$ of minimum
degree two and excess two. Let $I$ consist of the nonhub vertices of
$F$ at depth less than $6R$ in the stopped search from $x$, and put
$K=V(F)\setminus I$. The paths in the breadth-first tree connect $I$,
and deleting leaves preserves this property. Every edge of $F$ has
an endpoint in $I$, so $K$ is independent. Thus
\[
 \sum_{v\in V(F)}(\deg_F(v)-2)=2,
 \qquad
 2=\operatorname{tx}(F[I])+\sum_{z\in K}(\deg_F(z)-1).
\]
The first identity gives either one degree-four vertex or two
degree-three vertices. They cannot all lie in $K$, by the second
identity. Hence some $w\in I$ has degree at least three and label
$i_w>k$. Also, $F$ contains no seed edge, since both seeds are hubs.

Write $q=|E(F)|\le24R+2$. By \eqref{eq:embedding-bound}, summing
over labels with $i_w>k$, and dropping all other label restrictions,
gives at most
\[
 C^q\left(\sum_{i>k}i^{-3/2}\right)
       \left(\sum_{i=1}^n i^{-1}\right)^{|V(F)|-1}
 \le C^q k^{-1/2}(1+\log n)^{|V(F)|-1}.
\]
There are at most $(Cq)^{Cq}$ choices of witness shape, orientation,
partition $I\sqcup K$, and vertex $w$. Summing over them gives
\[
 \mathbb P\left(
   [k]\subseteq\mathcal C_n,\
   \max_x\operatorname{tx}(B_{6R}^{\mathrm{stop}}(x))\ge2
 \right)
 \le k^{-1/2}\exp\{CR[\log(R+1)+\log\log n]\}
 =n^{-\beta/4+o(1)}=o(1).
\]
Together with the early-vertex degree bound, this proves
\eqref{eq:stopped-geometry}.

For the walk estimate, fix a realization satisfying
\eqref{eq:stopped-geometry},
$x\notin\mathcal C_n$, and put $H=B_{6R}^{\mathrm{stop}}(x)$.
On $\{\tau_{\mathcal C_n}>5R\}$,
$B_R^{\mathrm{stop}}(X_{5R})\subseteq H$: every nonhub vertex at
depth less than $R$ in this stopped search is reachable from $x$
by a path avoiding hubs of length less than $6R$.
If $H$ is a tree, the result follows.
Otherwise, let $\mathscr C$ be its unique cycle and, for
$0\le t\le5R\wedge\tau_{\mathcal C_n}$, put
\[
 Y_t=\operatorname{dist}_H(X_t,\mathscr C)
       -\frac14\mathbf1_{\{\deg_G(X_t)=2\}}.
\]
Every vertex visited before the first hub hit and through time $5R$
has all its incident edges in $H$. A degree-two nonhub vertex has
no arrivals, so both neighbors have degree at least three. The
distance has mean increment zero there, and the correction gives
drift $1/4$. At a vertex of degree $d\ge3$, the sum of the distance
increments over its incident edges is $d-2$, whether or not it
lies on $\mathscr C$. At most $d-m\le d-2$ incident edges lead to
degree-two vertices, since each such edge is an arrival. The
corrected drift is therefore at least
\[
 \frac{(d-2)-\tfrac14(d-2)}d\ge\frac14.
\]
These estimates include transitions to hubs. Thus
$Y_{t\wedge\tau_{\mathcal C_n}}-\tfrac14(t\wedge\tau_{\mathcal C_n})$
is a submartingale with increments bounded in absolute value by
$3/2$. Since $Y_0\ge-1/4$, Azuma--Hoeffding gives
\[
 P_x^G\bigl(\tau_{\mathcal C_n}>5R,\ Y_{5R}\le R\bigr)\le e^{-cR}.
\]
If no hub has been hit and $X_{5R}$ is not good, its stopped
neighborhood contains $\mathscr C$, so $Y_{5R}\le R$.
Every hub is good, completing the proof.
\end{proof}

For the subsequent explorations, set
\[
 b=\lfloor ne^{-AL^2/2}\rfloor,\qquad
 \Lambda=n(b/n)^\chi,\qquad \rho=L^{-2}.
\]
Retain the urn notation of Section~\ref{sec:PA-PPT}.
Let $\Omega$ be the event that \eqref{eq:stopped-geometry} and the
following bounds hold:
\begin{equation}\label{eq:common-red-event}
 \begin{gathered}
  \frac{(j/n)^\chi}{C(\log n)^3}
  \le S_j\le C(\log n)^3(j/n)^\chi\quad(1\le j\le n),\\
  \psi_j\le C\log n/j\quad(2\le j\le n),\qquad
  \varphi_v\le CD/n\quad(v<a,\ v\notin\mathcal C_n);
 \end{gathered}
\end{equation}
\begin{equation}\label{eq:common-switch-event}
  \#\{(k,i):k\ge b,\ J_{k,i}<b\}\ge c\Lambda,\quad
  \sum_{v<a}\deg_{G_n}(v)\le Cn(a/n)^\chi,\quad 
  \max_{v\ge a}\deg_{G_n}(v)\le C(n/a)^\beta\log n;
\end{equation}
\begin{equation}\label{eq:common-age-volume}
 \deg_{G_{b-1}}(h)\le\rho\deg_{G_n}(h)\quad(h\in\mathcal C_n),
 \qquad \max_x|B_R^{\mathrm{stop}}(x)|\le e^{C_{\mathrm{vol}}L^2}.
\end{equation}

\begin{lem}
For suitable fixed constants in these bounds,
\begin{equation}\label{eq:common-regularity-input}
 \mathbb P(\Omega^c)=o(1).
\end{equation}
\end{lem}

\begin{proof}
The factorial-weight calculation in the preceding proof, with a
fixed integer $s\ge1$, gives
$\mathbb E[\deg_{G_n}(i)^s]\le(Cs)^s(n/i)^{\beta s}$.
Taking $s=\lceil\log n\rceil$, Markov's inequality and
a union bound prove the last bound in \eqref{eq:common-switch-event}
with failure probability at most $n^{-2}$.

The maximal-inequality argument in Lemma~\ref{lem:S_k-estimate},
using the summable variances of $\log(1-\psi_j)$, gives
$\max_{1\le j\le n}|\log S_j-\chi\log(j/n)|=O_{\mathbb P}(1)$.
Together with the exponential moments of $j\psi_j$ and a union
bound, this proves the first two bounds in \eqref{eq:common-red-event}
with probability $1-o(1)$. The same lemma gives
$S_j\asymp(j/n)^\chi$ uniformly for $j\ge a-1$ with high probability.
Fix urn variables satisfying these bounds. Conditional on them,
the attachment slots are independent.

For $v<a$, arrivals from labels at least $n/2$ have mean at least
$cn\varphi_v$. By Chernoff's inequality, the last bound in
\eqref{eq:common-red-event} fails with probability at most $ne^{-cD}$.
For $s=a,b$, the number of edges from $[s,n]$ to $[s-1]$ has
conditional mean
$mS_{s-1}\sum_{k=s}^nS_{k-1}^{-1}\asymp n(s/n)^\chi$.
Chernoff's inequality gives the first two bounds in
\eqref{eq:common-switch-event}, since the edges within $[a-1]$
contribute only $O(a)$ to its total degree.

For $h<a$, let $\mu_h^-$ and $\mu_h^+$ be the conditional means of
$\deg_{G_{b-1}}(h)$ and $\deg_{G_n}(h)-\deg_{G_{b-1}}(h)$.
The urn bounds give
\[
 \mu_h^+=m\varphi_h\sum_{k=b}^nS_{k-1}^{-1}\asymp n\varphi_h,
 \qquad
 \mu_h^-\le m+C(\log n)^3(b/n)^\beta\mu_h^+.
\]
Since $(\log n)^3(b/n)^\beta=o(\rho)$ and $\rho D\gg\log n$,
if $\mu_h^+\ge D/4$, Chernoff's inequality gives
\[
 \deg_{G_n}(h)-\deg_{G_{b-1}}(h)\ge\mu_h^+/2,
 \qquad \deg_{G_{b-1}}(h)\le\rho\mu_h^+/2,
\]
except with probability $Ce^{-c\rho D}$. If $\mu_h^+<D/4$, then
$\mathbb P(\deg_{G_n}(h)>D\mid\psi)\le Ce^{-cD}$.
A union bound over all $h<a$, without conditioning on whether $h$ is
a hub, bounds the conditional probability that the first inequality in
\eqref{eq:common-age-volume} fails by $Cne^{-c\rho D}=o(1)$.

For the volume bound, fix urn variables
satisfying the first two bounds in \eqref{eq:common-red-event}.
The attachment probabilities then satisfy
\[
 p_{i\to j}=\frac{\varphi_j}{S_{i-1}}
       \le C(\log n)^7i^{-\chi}j^{\chi-1}\quad(j<i).
\]
Put $f(v)=\sqrt{n/(v\vee a)}$, $F=\sqrt{n/a}$,
$K=(\log n)^{10}$ and
\[
 Q=\left\lceil4(R+1)FK^R\log n\right\rceil
   =\exp\{(A/2+10c_R)L^2+O(L+\log L)\}.
\]
From a fixed root, run breadth-first search to depth $R$, stopping
when $|E|=Q$. Examine each vertex at depth less than $R$ once,
revealing its complete star if it is young; at an old vertex, stop
when its star is complete or $D+1$ incidences have been recorded.
Only vertices with complete recorded stars become explored.
If the search finishes before reaching $Q$ labels, it contains the
stopped neighborhood: every old nonhub vertex at depth less than
$R$ has its complete star recorded, and queries at hubs only add vertices.

Reveal outgoing slots by their endpoints and test incoming slots
one target at a time, retaining every answer. Conditional on this
record and the urn, the urn uniforms for slots whose endpoints have
not been revealed remain independent and uniform outside their
excluded target intervals. All excluded
targets belong to $E$. Their total relative mass at a sender
$i\ge a$ is at most
\[
 C(\log n)^7i^{-\chi}\sum_{j\le Q}j^{\chi-1}
 \le C(\log n)^7(Q/a)^\chi=o(1).
\]
For such a slot with a young sender, each conditional attachment
probability is at most twice its value given the urn alone.
Neither membership in $\mathcal C_n$
nor final degrees are included in this conditioning.

For a young vertex $v\ge a$, the expected weight of its neighbors is
bounded by
\begin{equation}\label{eq:stopped-weighted-offspring}
 m\sum_{j<v}f(j)p_{v\to j}+m\sum_{i>v}f(i)p_{i\to v}
 \le C(\log n)^7\bigl(1+\log(n/a)\bigr)f(v).
\end{equation}
Indeed, after division by $C(\log n)^7f(v)$, the contribution from
$j<a$ is at most $C(a/v)^{\chi-1/2}$, and the remaining sums are
bounded by
\[
 C\sum_{a\le j<v}\frac{(j/v)^{\chi-1/2}}j
 +C\sum_{i>v}\frac{(v/i)^{\chi-1/2}}i
 \le C\bigl(1+\log(n/a)\bigr).
\]
This includes $\chi=1/2$. At an old vertex, the total weight of
neighbor labels added to $E$ is at most $(D+1)F=(D+1)f(v)$.

Let $W_j$ be the total $f$-weight of labels added to $E$ at depth
$j$, and let $\mathscr F_j$ be the record before that generation is
examined. Each new label has weight at most $F$. Using
$e^{\theta u}-1\le(e-1)\theta u$ for $0\le\theta u\le1$,
the conditional attachment bound and \eqref{eq:stopped-weighted-offspring}
give, by successive exposure within the generation,
\[
 \mathbb E[e^{\theta W_{j+1}}\mid\mathscr F_j,\psi]
 \le e^{K\theta W_j},\qquad 0\le\theta\le F^{-1}.
\]
The deterministic contribution from old vertices is absorbed by
$K$; collisions and the $Q$-label stop only reduce the added weight.
Iteration yields
$\mathbb E[e^{W_j/(FK^j)}\mid\psi]\le e^{f(x)/F}\le e$.
If the search reaches $Q$ labels, some generation $j\le R$ has weight at least
$Q/(R+1)\ge4FK^j\log n$. Exponential Markov inequality and a union
bound over $j\le R$ give probability at most $e(R+1)n^{-4}$.
Choose $C_{\mathrm{vol}}>A/2+10c_R$, so that $Q\le e^{C_{\mathrm{vol}}L^2}$
for large $n$. A union bound over the roots bounds the conditional
probability of failure of the volume bound by $CRn^{-3}$. Integrating
over the urn and using the preceding estimates and
Lemma~\ref{lem:stopped-tangle-free} proves
\eqref{eq:common-regularity-input}.
\end{proof}

With the scales of Section~\ref{subsec:ppt-pa-coupling}, let
$\mathcal U_n$ be the event on which all the following bounds hold:
\begin{equation}
 \max_{k_n^\circ-1\le k\le n}
 \left|\log S_k-\chi\log(k/n)\right|\le\varepsilon_n,\qquad
 \max_{k_n^\circ\le j\le n}(Z_j+j\psi_j)\le\bar g_n;
\label{eq:good-urn-event}
\end{equation}
\begin{equation}\label{eq:degree-mass-regularity}
 \max_v\deg_{G_n}(v)\le n^{3/4},\qquad
 \varphi_v\le C\bigl(\deg_{G_n}(v)+D\bigr)/n\quad(v<a),
\end{equation}
where $C$ is a fixed constant. For $C_G$ and $C$ sufficiently large,
\[
 \mathbb P(\mathcal U_n^c)
 \le 2n e^{-c\varepsilon_n^2k_n^\circ}+O(n^{-8})+ne^{-cD}
 =O(n^{-8}).
\]
Indeed, the beta moment formula gives
\[
 \frac{d^2}{d\theta^2}\log\Ebb(1-\psi_j)^\theta
 =\Psi_1(d_j+\theta)-\Psi_1(c_j+\theta)
 \le Cj^{-2},\qquad |\theta|\le cj.
\]
Here $c_j-d_j=m+\delta$ is fixed and $d_j\asymp j$.
Independence therefore yields
\[
 \log\Ebb\exp\{\theta(\log S_k-\Ebb\log S_k)\}
 \le C\theta^2/k,\qquad |\theta|\le ck.
\]
Since $\Ebb\log S_k=\chi\log(k/n)+O(k^{-1})$, Chernoff's bound
and a union bound give the first term. It is smaller than every
fixed negative power of $n$, because
$\varepsilon_n^2k_n^\circ=n^{1/3-o(1)}$.
The Gamma exponential moments and
$\Ebb e^{d_j\psi_j/2}\le2^{m+\delta}$ bound failure of the second
inequality in \eqref{eq:good-urn-event} by $O(n^{-8})$ after increasing $C_G$.
The degree-moment bound above with $s=40$ gives failure probability
$O(n^{-9})$ for the first inequality in \eqref{eq:degree-mass-regularity}.
Conditional on the urn, arrivals to $v<a$ from labels at least $n/2$
have mean at least $cn\varphi_v$. If the second inequality fails,
this mean exceeds both $2\deg_{G_n}(v)$ and $cD$, after increasing
$C$. Chernoff's inequality and a union bound give failure probability
at most $ne^{-cD}$.

Graph-level error bounds below retain their intersection with
$\Omega\cap\mathcal U_n$ until the union over starting vertices in
Proposition~\ref{prop:worst-start-reduction}. The bounds
\eqref{eq:common-red-event} are determined by $(\psi,\mathcal C_n)$,
and \eqref{eq:common-switch-event} by the final degree sequence.
The remaining argument controls $\|P_G^t(x_0,\cdot)-\pi_G\|_{\mathrm{TV}}$
uniformly over good roots.

\subsection{Reduction to a successful continuation}
\label{subsec:one-round-acceptance}

Work on $\Omega$ and fix a good root $x_0$ in the realized graph $G$.
Put $T=B_R^{\mathrm{stop}}(x_0)$. We follow the walk to a terminal of $T$
and try to continue it through a young vertex using the coupling of
Section~\ref{sec:coupling}. A return allows another attempt; a defect
ends the current walk procedure. We bound the probability that an accepted
continuation ends at a prescribed vertex. Together with control of the
duration and the spectral gap, this reduces the mixing estimate to
bounding the probability of failure.

\begin{figure}[t]
\centering
\begin{tikzpicture}[x=1cm,y=1cm,>=Stealth,
    tv/.style={circle,fill=black,inner sep=1.45pt},
    ov/.style={rectangle,fill=black,inner sep=1.75pt},
    hv/.style={circle,draw=black,fill=gray!45,minimum size=9pt,inner sep=0pt,line width=0.5pt},
    ed/.style={black!80,line width=0.55pt},
    cont/.style={black!70,line width=0.55pt,dash pattern=on 0.6pt off 1.6pt,line cap=round},
    stub/.style={gray!65,line width=0.45pt},
    leader/.style={gray!55!black,line width=0.3pt},
    ann/.style={font=\footnotesize,inner sep=1.5pt}]
  \def\R{3.3}
  % ---- the ball
  \fill[blue!4] (0,0) circle (\R);
  \draw[gray!70,densely dashed,line width=0.5pt] (0,0) circle (\R);
  % ---- tree: level 1 (r=0.95), level 2 (r=1.75), level 3 (r=2.45)
  \foreach \a in {70,190,310} \draw[ed] (0,0) -- (\a:0.95);
  \foreach \p/\c in {70/40,70/95,190/155,190/215,310/270,310/330}
     \draw[ed] (\p:0.95) -- (\c:1.75);
  \foreach \p/\c in {40/22,40/50,95/85,95/110,155/140,155/168,270/252,270/282,330/315,330/345}
     \draw[ed] (\p:1.75) -- (\c:2.45);
  % ---- the remaining levels down to depth R
  \foreach \p/\c in {22/16,22/28,50/50,85/80,85/92,110/112,140/138,168/170,
                     252/250,282/276,282/290,315/316,345/340,345/352}
     \draw[cont] (\p:2.45) -- (\c:\R);
  % ---- the hub at level 2: the search stops here; forward edges lie outside T
  \foreach \a in {185,195,...,265}
     \draw[stub] (215:1.75) -- ++(\a:0.72) node[circle,fill=gray!55,inner sep=0.7pt]{};
  % ---- path to the terminal where the walk first hits the terminal set
  \draw[blue!70!black,line width=1.3pt,line cap=round,line join=round]
     (0,0) -- (70:0.95) -- (40:1.75) -- (50:2.45);
  \draw[blue!70!black,line width=1.3pt,dash pattern=on 0.9pt off 1.9pt,line cap=round]
     (50:2.45) -- (50:\R);
  \draw[blue!70!black,line width=0.8pt] (50:\R) circle (3.6pt);
  % ---- vertices
  \foreach \a in {70,190,310} \node[tv] at (\a:0.95) {};
  \foreach \a in {40,95,155,270,330} \node[tv] at (\a:1.75) {};
  \foreach \a in {22,50,85,110,140,168,252,282,315,345} \node[tv] at (\a:2.45) {};
  \foreach \a in {16,28,50,80,92,112,138,250,276,290,340,352} \node[tv] at (\a:\R) {};
  \foreach \a in {170,316} \node[ov] at (\a:\R) {};
  \node[hv] at (215:1.75) {};
  \node[circle,fill=black,inner sep=2pt] at (0,0) {};
  \node[ann,anchor=north east] at (-0.06,-0.08) {$x_0$};
  % ---- labels
  \node[ann,gray!40!black] at (64:3.62) {depth $R$};
  \draw[leader] (16:3.42) -- (3.95,1.25) node[ann,anchor=west,black]
     {$\partial_{\mathrm{reg}}T$: young, at depth $R$};
  \draw[leader] (316:3.45) -- (3.95,-2.55) node[ann,anchor=west,black]
     {$\partial_{\mathrm{bad}}T$: old, $\deg_G\le D$, at depth $R$};
  \draw[leader] (-1.63,-0.95) -- (-3.95,-0.35) node[ann,anchor=east,black,align=right]
     {$\partial_{\mathrm{hub}}T$: hubs, at any depth,\\ retained as terminals};
\end{tikzpicture}
\caption{The stopped neighbourhood $T=B_R^{\mathrm{stop}}(x_0)$ of a good
root: hubs and vertices at depth $R$ are retained as terminals.
Circles are young vertices, squares old nonhub vertices
and the shaded circle a hub; the blue path ends at the first terminal
hit by the walk.}
\label{fig:stopped-ball}
\end{figure}
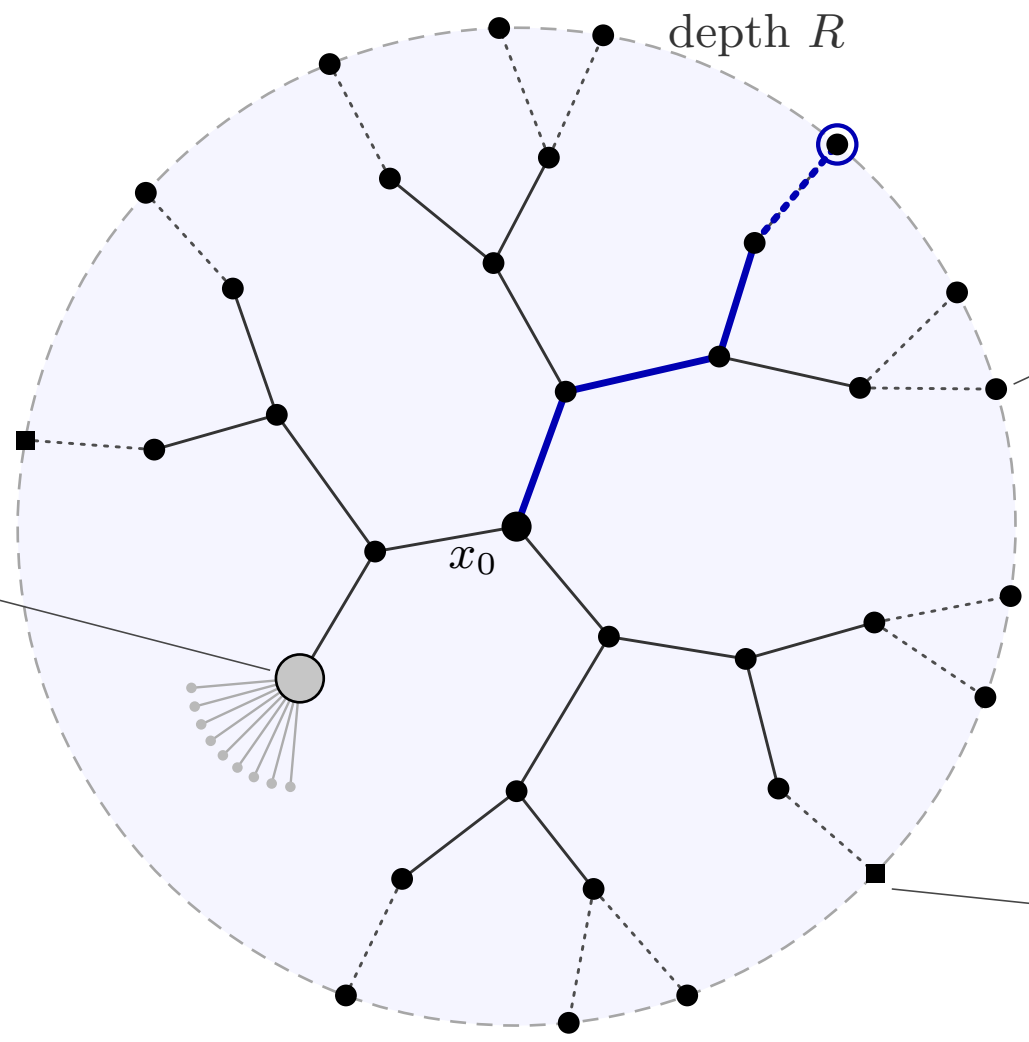

The terminal set is
\[
 \partial T=\partial_{\mathrm{reg}}T\,\dot\cup\,
             \partial_{\mathrm{bad}}T\,\dot\cup\,
             \partial_{\mathrm{hub}}T,
\]
where the first two sets contain, respectively, the young and old
nonhub vertices at depth $R$, and the third contains the hubs
of $T$, at any depth. Write $v^-$ for
the parent of a nonroot vertex and $\mathscr H_0$ for the exploration
record of $T$; see Figure~\ref{fig:stopped-ball}.
To explore $T$ and its old nonhub terminals, first reveal the hub-status
indicator of each encountered old vertex. At this stage scan incoming
slots only at nonhubs, postponing all hub scans until the old nonhub
terminals have been explored. The auxiliary search stopped at $D+1$
incidences above is used only for the volume bound.

For an oriented incident slot $g=(u,w)$ with $w\ge a$, assign
the finite matched region $\widetilde\Gamma$ constructed in
Section~\ref{sec:coupling-after-burn-in}, retaining $u$ as the
entrance parent. The \emph{forward test} is the test of
Lemma~\ref{lem:local-continuation} with testing depth $R$, started
at $w$. It stops at the first entrance \emph{return}, first other exit
from $\widetilde\Gamma$ (a \emph{defect}), or first visit to
depth $R$ in the retained region (a \emph{penetration}). An empty region gives
an immediate defect. Let $\xi_g$ be the test duration from $w$.

\begin{figure}[t]
\centering
\begin{tikzpicture}[x=1cm,y=1cm,>=Stealth,
    tv/.style={circle,fill=black,inner sep=1.45pt},
    ov/.style={rectangle,fill=black,inner sep=1.75pt},
    hv/.style={circle,draw=black,fill=gray!45,minimum size=9pt,inner sep=0pt,line width=0.5pt},
    ed/.style={black!80,line width=0.55pt},
    cont/.style={black!70,line width=0.55pt,dash pattern=on 0.6pt off 1.6pt,line cap=round},
    stub/.style={gray!65,line width=0.45pt},
    walk/.style={blue!70!black,line width=1.3pt,line cap=round,line join=round},
    walkcont/.style={blue!70!black,line width=1.3pt,dash pattern=on 0.9pt off 1.9pt,line cap=round},
    move/.style={->,blue!70!black,line width=0.8pt},
    ring/.style={blue!70!black,line width=0.8pt},
    ann/.style={font=\footnotesize,inner sep=1.5pt},
    ttl/.style={font=\small,anchor=west},
    ftest/.pic={
      \fill[blue!8] (0,0) -- (-30:0.7) arc[start angle=-30,end angle=30,radius=0.7] -- cycle;
      \draw[gray!70,densely dashed,line width=0.4pt] (-30:0.7) arc[start angle=-30,end angle=30,radius=0.7];},
    dx/.pic={
      \draw[red!75!black,line width=1pt,line cap=round]
        (-0.09,-0.09) -- (0.09,0.09) (-0.09,0.09) -- (0.09,-0.09);}]
  % =====================================================================
  % (a) young terminal: the forward test
  % =====================================================================
  \node[ttl] at (0,2.1) {\textbf{(a)}\ \ $v\in\partial_{\mathrm{reg}}T$: forward test};
  \draw[ed] (0,0) -- (0.95,0);
  \begin{scope}[shift={(0.95,0)}]
    % the forward coupling to depth R, away from the entrance
    \fill[blue!4] (0,0) -- (-35:3) arc[start angle=-35,end angle=35,radius=3] -- cycle;
    \draw[gray!70,densely dashed,line width=0.5pt] (-35:3) arc[start angle=-35,end angle=35,radius=3];
    \foreach \a in {-20,17} \draw[ed] (0,0) -- (\a:1);
    \foreach \p/\c in {-20/-27,-20/-10,17/7,17/26} \draw[ed] (\p:1) -- (\c:1.85);
    \foreach \p/\c in {-10/-15,-10/-5,7/4,7/14,26/28} \draw[cont] (\p:1.85) -- (\c:3);
    % (i) penetration: the walk reaches depth R
    \draw[walk] (0,0) -- (17:1) -- (7:1.85);
    \draw[walkcont] (7:1.85) -- (4:3);
    \foreach \a in {-20,17} \node[tv] at (\a:1) {};
    \foreach \a in {-10,7,26} \node[tv] at (\a:1.85) {};
    \foreach \a in {-15,-5,4,14,28} \node[tv] at (\a:3) {};
    \draw[ring] (4:3) circle (3.6pt);
    % defect: the coupling fails
    \pic at (-27:1.85) {dx};
  \end{scope}
  \node[tv] at (0,0) {};  \node[tv] at (0.95,0) {};
  \node[ann,anchor=north] at (0,-0.1) {$v^-$};
  \node[ann,anchor=north] at (0.95,-0.1) {$v$};
  % (ii) return through the entrance
  \draw[move] (0.85,0.13) to[bend right=45] (0.1,0.13);
  \node[ann,blue!60!black,anchor=south] at (0.47,0.36) {return};
  \node[ann,red!70!black,anchor=north] at (2.45,-1.12) {defect};
  \draw[move] (4.12,0.21) -- (5.85,0.21);
  \node[ann,blue!60!black,anchor=south west] at (4.3,0.25) {penetration};
  \node[ann,anchor=north west] at (4.3,0.15) {then to depth $\ell_1$};
  \node[ann,gray!40!black,anchor=west] at (3.55,-1.85) {depth $R$};
  % =====================================================================
  % (b) old nonhub terminal
  % =====================================================================
  \begin{scope}[xshift=6.75cm]
    \node[ttl] at (0,2.1) {\textbf{(b)}\ \ $y\in\partial_{\mathrm{bad}}T$: old nonhub};
    \draw[ed] (0,0) -- (1.1,0);
    \begin{scope}[shift={(1.1,0)}]
      \draw[walk] (0,0) -- (40:1.2);
      \draw[ed] (0,0) -- (0:1.3);
      \draw[ed] (0,0) -- (-40:1.2);
      \pic[rotate=40] at (40:1.2) {ftest};
      \node[tv] at (40:1.2) {};
      \node[hv] at (0:1.3) {};
      \node[ov] at (-40:1.2) {};
      \pic at (1.17,-0.77) {dx};
    \end{scope}
    \node[tv] at (0,0) {};  \node[ov] at (1.1,0) {};
    \node[ann,anchor=north] at (0,-0.1) {$y^-$};
    \node[ann,anchor=north] at (1.1,-0.12) {$y$};
    \draw[move] (1.0,0.13) to[bend right=45] (0.1,0.13);
    \node[ann,blue!60!black,anchor=south] at (0.55,0.38) {return};
    \node[ann,anchor=west] at (2.85,0.95) {forward test};
    \node[ann,anchor=west] at (2.62,0) {hub procedure};
    \node[ann,red!70!black,anchor=west] at (2.45,-0.77) {defect};
  \end{scope}
  % =====================================================================
  % (c) hub
  % =====================================================================
  \begin{scope}[xshift=12.75cm]
    \node[ttl] at (0,2.1) {\textbf{(c)}\ \ $h\in\partial_{\mathrm{hub}}T$: hub};
    \draw[ed] (0,0) -- (1.0,0);
    \begin{scope}[shift={(1.0,0)}]
      \foreach \a in {-120,-105,-90,-75,-60,-15,0,45,60,75,90,105,120}
         \draw[stub] (0,0) -- (\a:1.15) node[circle,fill=gray!55,inner sep=0.7pt]{};
      \draw[walk] (0,0) -- (20:1.3);
      \draw[ed] (0,0) -- (-40:1.3);
      \pic[rotate=20] at (20:1.3) {ftest};
      \node[tv] at (20:1.3) {};
      \node[ov] at (-40:1.3) {};
      \pic at (1.25,-0.84) {dx};
    \end{scope}
    \node[tv] at (0,0) {};  \node[hv] at (1.0,0) {};
    \node[ann,anchor=north east] at (0.86,-0.08) {$h$};
    \node[ann,anchor=west] at (3.0,0.45) {forward test};
    \node[ann,red!70!black,anchor=west] at (2.4,-0.84) {non-young: defect};
  \end{scope}
\end{tikzpicture}
\caption{The terminal procedures, one for each class of terminal in
Figure~\ref{fig:stopped-ball}. (a) The forward test at a young terminal;
the shaded sector indicates the forward test up to depth $R$.
(b) An old nonhub terminal chooses an incident slot and follows the
indicated rule. (c) At a hub, young slots with assigned couplings
start forward tests; all other departures give a defect. In (b)
and (c), a returning forward test is followed by another slot choice.}
\label{fig:terminal-procedures}
\end{figure}
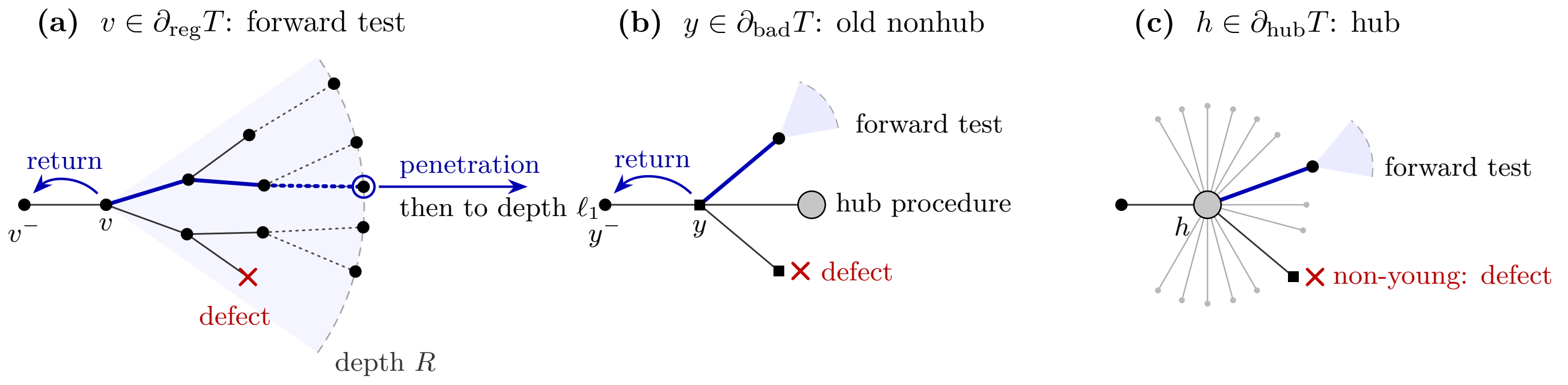

The rule at a terminal is as follows. All choices are uniform over
incident slots, with parallel edges counted separately.
\begin{enumerate}[label=\textup{(\roman*)}]
 \item At a young terminal $v$, run the forward test with entrance
 $(v^-,v)$. A return ends the procedure at $v^-$.
 \item At an old nonhub terminal $v$, take one SRW step. A step through
 the entrance ends the procedure by a return to $v^-$. A forward step
 to a young vertex starts a forward test; if that test returns to $v$,
 choose again. A forward step to a hub starts the hub procedure below.
 A forward step to another old nonhub vertex gives a defect.
 \item At a hub $h$, take one SRW step. A young slot with an assigned
 coupling starts a forward test; a return to $h$ is followed by another choice.
 Every other departure gives a defect. In particular, entrances from
 the initial stopped exploration or an old nonhub terminal are excluded.
\end{enumerate}
These rules are illustrated in Figure~\ref{fig:terminal-procedures}.
A penetration or defect is called a \emph{resolution}. Thus a nonhub
procedure either returns to its parent or resolves, whereas a hub
procedure continues until it resolves. A test declared defective in the
explorations below stops immediately. The coupled environments
and matching assigned to each slot are kept fixed on repeated visits.
The threshold $a$ governs the exploration of $T$ and its old nonhub
terminals, as well as the terminal rules.
Within a forward test or continuation, use the regularity and collision
rules of Section~\ref{sec:coupling}, with label cutoff $k_n^\circ<a$;
old vertices newly added to $E$, including hubs, require no degree
classification.
% ---------------------------------------------------------------------
% Section 5.2: the terminal procedures, one panel per terminal class.
% Same visual language as fig:stopped-ball: blue = the walk's path,
% ring = terminal reached, dashed arc = depth R, red cross = defect,
% small shaded sector = a forward test starts here.
% Needs tikz + xcolor only; no tikz library beyond what sy.sty loads.
% ---------------------------------------------------------------------

After penetration within $R^2$ steps, continue in the same
retained region as in Lemma~\ref{lem:local-continuation}, with
time and depth measured from the entrance root. Use
$C_{\rm cap}>2\gamma$ and
\[
 \kappa_n=n^{-1}e^{C_{\rm cap}\sqrt{\log n}}
\]
from that lemma. The continuation is \emph{accepted} if it reaches
matched depth $\ell_1$ before exiting and by time
$\lceil t_n^\star\rceil$, and its depth-$\ell_0$ ancestor $y$
satisfies
\[
 h_{\widetilde\Gamma}(y)\le\kappa_n/p_{\rm esc}.
\]
Here $h_{\widetilde\Gamma}$ is the killed hitting probability
defined there, using the recorded incident slots of the entire
retained region. Acceptance is determined by this finite record
and the walk at its stopping time.

Starting at $x_0$, run the appropriate procedure whenever the walk hits
$\partial T$, and resume the walk in $T$ after a return. Let $\tau_j$
be the start of the $j$th procedure, $N$ the first resolving index,
and $\sigma$ its completion time, before the continuation. Define
$\mathcal E$ to be the event that this resolution is a penetration
with an accepted continuation. Let $S$ be the time of acceptance,
with $S=\infty$ on $\mathcal E^c$. Acceptance is determined at that
time, and $S\le\sigma+\lceil t_n^\star\rceil$ on $\mathcal E$.

\begin{lem}
Uniformly over these choices and all good roots,
\begin{equation}\label{eq:resolution-time-cap}
 E_{x_0}^G\sigma\le C_mR,\qquad
 P_{x_0}^G(\mathcal E,\ X_S=z)\le C_m\kappa_n
 \quad(z\in V(G)).
\end{equation}
\end{lem}

\begin{proof}
For a procedure starting at $v$, let $\xi_v$ be its duration,
$r_v$ its resolution probability, and $q_v(z)$ its accepted endpoint
mass. We prove
\begin{equation}\label{eq:terminal-cap-time}
 r_v\ge p:=\frac{p_{\rm esc}}{1+p_{\rm esc}}=\frac14,\qquad
 E_v^G\xi_v\le C_mR,\qquad
 q_v(z)\le C_m\kappa_n r_v.
\end{equation}

For a forward test at $g=(u,w)$, let $r_g$ be its resolution
probability and $q_g(z)$ its accepted endpoint mass.
Lemma~\ref{lem:local-continuation} gives $r_g\ge p_{\rm esc}=1/3$,
$E_w^G\xi_g\le C_mR$, and
\begin{equation}\label{eq:elementary-endpoint-cap}
 q_g(z)\le\frac{\kappa_n}{p_{\rm esc}^2}r_g.
\end{equation}
These bounds also hold for a declared defect, with $r_g=1$,
$\xi_g=0$, and $q_g=0$.

At a hub $h$, let $r_e$ and $q_e(z)$ be the resolution probability and
accepted endpoint mass after departure through slot $e$. Set $r_e=1$
and $q_e(z)=0$ for a defective departure. Summing over successive
excursions that return to $h$ gives
\[
 q_h(z)
   =\frac{\sum_e q_e(z)}{\sum_e r_e}
   \le\frac{\kappa_n}{p_{\rm esc}^2}.
\]
The expected number of departures is
$\deg_G(h)/\sum_e r_e\le p_{\rm esc}^{-1}$. Each departure takes
expected time at most $1+C_mR$, so $E_h^G\xi_h\le C_mR$ and $r_h=1$.

At an old nonhub terminal $v$, put $S_v=\sum_{e\ne e_v}r_e$, where
$e_v$ is its entrance and $r_e,q_e(z)$ refer to the corresponding
forward departure. A departure into a hub includes the entire hub
procedure, so it has $r_e=1$ and satisfies the same endpoint bound.
Every forward slot has $r_e\ge p_{\rm esc}$, and therefore
\[
 r_v=\frac{S_v}{1+S_v}\ge p,\qquad
 q_v(z)
   =\frac{\sum_{e\ne e_v}q_e(z)}{1+S_v}
   \le\frac{\kappa_n}{p_{\rm esc}^2}r_v.
\]
The expected number of choices is
$\deg_G(v)/(1+S_v)\le p_{\rm esc}^{-1}$, proving the duration bound
and hence \eqref{eq:terminal-cap-time}.

It follows that $E_{x_0}^GN\le p^{-1}$. Between procedures, the
potential
$\operatorname{dist}_T(x_0,X_t)-\tfrac14\mathbf1_{\{\deg_G(X_t)=2\}}$
has drift at least $1/4$ until $\partial T$ is hit, by the calculation
in Lemma~\ref{lem:cycle-clearing}. Here $\Omega$ ensures that the seed
vertices are hubs. Thus the expected time to reach $\partial T$ is
at most $4(R+1)$ from any point of $T$. Together with
\eqref{eq:terminal-cap-time}, this gives $E_{x_0}^G\sigma\le C_mR$.

Finally, the strong Markov property at the procedure starts gives
\[
\begin{aligned}
 P_{x_0}^G(\mathcal E,X_S=z)
 &=\sum_{j\ge1}\sum_{v\in\partial T}
   P_{x_0}^G(N\ge j,X_{\tau_j}=v)
   q_v(z)\\
 &\le C_m\kappa_n\sum_{j\ge1}\sum_{v\in\partial T}
   P_{x_0}^G(N\ge j,X_{\tau_j}=v)r_v
 =C_m\kappa_n,
\end{aligned}
\]
since the last sum is $P_{x_0}^G(N<\infty)=1$.
\end{proof}

\begin{prop}\label{prop:worst-start-reduction}
Fix $\varepsilon>0$. Suppose the couplings can be chosen so that,
uniformly in $x_0$,
\begin{equation}\label{eq:boundary-input}
 \mathbb P\bigl(\Omega\cap\mathcal U_n\cap
   \{x_0\text{ is good},\ P_{x_0}^G(\mathcal E^c)>\varepsilon\}\bigr)
   =o(n^{-1}).
\end{equation}
Then, for a sufficiently large fixed $C$, with high probability,
\[
 d_n(\lceil t_n^\star\rceil+R^2+5R+s_n)\le\varepsilon+o(1),\qquad
 s_n=\lceil C\sqrt{\log n}\rceil.
\]
\end{prop}

\begin{proof}
Construct all root experiments on one common realization of the PA
graph, urn, and marks. Each experiment uses only its own retained
record in the conditional estimates. By \eqref{eq:boundary-input},
\[
\begin{split}
 &\mathbb P\bigl(\exists x_0\text{ good}:
                  P_{x_0}^G(\mathcal E^c)>\varepsilon\bigr)\\
 &\quad\le\mathbb P(\Omega^c)+\mathbb P(\mathcal U_n^c)\\
 &\qquad+\sum_{x_0=1}^n\mathbb P\bigl(\Omega\cap\mathcal U_n\cap
   \{x_0\text{ is good},\ P_{x_0}^G(\mathcal E^c)>\varepsilon\}\bigr)
 =o(1).
\end{split}
\]
Thus the common event $\mathcal U_n^c$ is excluded once, outside the
union over roots; independence of the root experiments is not needed.
Fix such a realization and a good root, and set
\[
 \nu(z)=P_{x_0}^G(\mathcal E,\ \sigma\le R^2,\
                         X_{\lceil t_n^\star\rceil+R^2}=z).
\]
By \eqref{eq:resolution-time-cap},
$\nu\le P_G^{\lceil t_n^\star\rceil+R^2}(x_0,\cdot)$ and
$\nu(V(G))\ge1-\varepsilon-C_m/R$.
For each $j\le\lceil t_n^\star\rceil+R^2$, the event
$\mathcal E\cap\{\sigma\le R^2,\ S=j\}$ is measurable at time $j$,
and its endpoint measure has atoms at most $C_m\kappa_n$.
Since $\pi_{\min}(G)\ge1/[2(n-1)]$, this measure is bounded by
$C_me^{C_{\rm cap}\sqrt{\log n}}\pi_G$.
Applying the Markov property for the remaining
$\lceil t_n^\star\rceil+R^2-j$ steps,
using stationarity, and summing over $j$ gives
\[
 \nu\le C_m(\lceil t_n^\star\rceil+R^2+1)
       e^{C_{\rm cap}\sqrt{\log n}}\pi_G.
\]
By Proposition~\ref{prop:PA-uniform-two-sided-gap}, the absolute spectral
gap is at least some fixed $g>0$ with high probability. Consequently,
\[
 \|\nu P_G^{s_n}-\nu(V(G))\pi_G\|_{\mathrm{TV}}
 \le C_m\sqrt{\lceil t_n^\star\rceil+R^2+1}\,
      e^{C_{\rm cap}\sqrt{\log n}/2-gs_n}=o(1)
\]
when $C>C_{\rm cap}/(2g)$, since $t_n^\star+R^2=O(\log n)$.
Including the discarded mass gives the
required bound from every good root.

From an arbitrary vertex, stop at the first hub hit or at time $5R$.
Lemma~\ref{lem:cycle-clearing} places the walk at a good root except
with probability $e^{-cR}$. The strong Markov property and contraction
of total variation give the additional $5R$ in the statement.
The resulting bound depends only on $G$, so it also holds with high
probability under the graph law.
\end{proof}

Sections~\ref{subsec:boundary-losses} and~\ref{subsec:hub-vertices}
estimate the failures in \eqref{eq:boundary-input}, at nonhub terminals
and hubs, respectively.

\subsection{Nonhub terminals}
\label{subsec:boundary-losses}

We estimate failures at nonhub terminals by their first hitting masses.
Let $\eta$ be the first hitting distribution on $\partial T$ from $x_0$.

\begin{lem}
For every nonhub terminal $v$,
\begin{equation}\label{eq:boundary-hitting-comparison}
 P_{x_0}^G(X_{\tau_N}=v)\le C\eta(v),\qquad
 \eta(v)\le w_*:=e^{-c_{\rm bd}R},\qquad
 c_{\rm bd}=\tfrac12\log\tfrac54.
\end{equation}
\end{lem}

\begin{proof}
If $x_0$ is a hub, there are no nonhub terminals. Otherwise, fix the
graph and all assigned coupled environments. Each terminal procedure
can be replaced by absorption with probability $r_v$ and a return
to its parent otherwise. For a nonroot vertex $u$, let $h(u)$ be the
probability of absorption below $u$ before reaching its parent, and
let $h_1(u)$ be the same probability when every terminal absorbs.
At terminals, $h(v)=r_v$ and $h_1(v)=1$. At internal vertices,
summing excursions that return to $u$ gives
\[
 h(u)=\frac{\sum_{z:z^-=u}h(z)}{1+\sum_{z:z^-=u}h(z)},
\]
with the same recursion for $h_1$.

By \eqref{eq:terminal-cap-time} and its proof, $r_v\ge1/4$, with
$r_v\ge1/3$ when $\deg_G(v)\ge3$. Since degree-two vertices are
not adjacent, the recursion extends these bounds to every nonroot
vertex: a degree-two internal vertex has
one child of degree at least three, and every other internal
nonroot vertex has at least two children. Thus $h(u)\ge1/4$, with
$h(u)\ge1/3$ when $\deg_G(u)\ge3$.

We compare the two absorption rules as in the proof of
Lemma~\ref{lem:block-approx-error}. A contribution to $h_1(u)-h(u)$
requires hitting a terminal before $u^-$ and then reaching $u^-$
without absorption. This terminal must be a nonhub, since hubs
absorb surely. It has depth $R$, so the return to $u^-$ crosses
$R-|u|+1$ successive parent edges. At each of these vertices, the
probability of reaching its parent before absorption is at most
$3/4$. The strong Markov property gives
\[
 0\le h_1(u)-h(u)\le(3/4)^{R-|u|+1},\qquad
 \frac{h_1(u)}{h(u)}\le1+3(3/4)^{R-|u|}.
\]

The excursion calculation in Lemma~\ref{lem:harmonic-kernel} gives
the conditional probability $h(z)/\sum_{y:y^-=u}h(y)$ that absorption
occurs below child $z$, given absorption before reaching $u^-$.
The same formula holds at the root and with $h_1$ in place of $h$.
At depth $j$, the ratio of these child probabilities is therefore
at most $1+3(3/4)^{R-j-1}$. Multiplying along the path to $v$ gives
\[
 \frac{P_{x_0}^G(X_{\tau_N}=v)}{\eta(v)}
 \le\prod_{j=0}^{R-1}\bigl(1+3(3/4)^{R-j-1}\bigr)
 \le\exp\left\{3\sum_{j\ge0}(3/4)^j\right\}=e^{12}.
\]

For the first hitting distribution, each choice of a prescribed
child at a vertex with at least two children has probability at
most $1/(1+1/4)=4/5$. The root has at least two children, and
degree-two vertices are not adjacent, so the path to $v$ contains
at least $\lceil R/2\rceil$ such choices. Hence
$\eta(v)\le(4/5)^{\lceil R/2\rceil}\le w_*$.
\end{proof}

Explore the old nonhub terminals. Throughout this subsection, $y$ ranges over
these terminals, and $c_y$ is the number of their forward slots
leading to old nonhub vertices. Forward slots exclude the entrance. Couplings are assigned
to the young entrances at young terminals and those reached from old
nonhub terminals. Give them weights
\[
 w_g=
 \begin{cases}
  \eta(v),&g=(v^-,v),\quad v\text{ a young terminal},\\
  \eta(y)/(\deg_G(y)-1),&g\text{ a young forward slot at }y.
 \end{cases}
\]
Then
$\sum_gw_g\le1,\, \max_gw_g\le w_*,\, \sum_gw_g^2\le w_*$.

The required hubs are the hub terminals of $T$ and the hubs reached
from its old nonhub terminals. Once those terminals have been explored,
complete the stars of every old vertex classified while exploring
$T$ and its old nonhub terminals,
including all required hubs, before starting any PA--PPT coupling.
No further degree classifications are made when these additional
neighbors are disclosed. Retain every
slot answer, expose the marks of queried targets at least $k_n^\circ$,
and reveal the normalized masses $\vartheta_\ell$ of the queried targets
below $k_n^\circ$ as in Section~\ref{subsec:ppt-pa-coupling}.
Every earlier hub-status indicator is now determined by the recorded
stars and adds no information to the PA record.

Before its own coupling begins, an entrance with a young endpoint is
\emph{reserved} relative to the current record if the endpoint's complete
star has not been recorded, its mark is unexposed, and the entrance is its only
recorded incidence. We also call the endpoint reserved for that entrance.
An entrance is \emph{spoiled} if, before its own coupling begins,
another incidence at its endpoint has been recorded; this includes
incidences recorded before the entrance itself. Its coupling is omitted.
Reserved endpoints' marks are exposed only when their own couplings
begin, which ends the reservation without spoiling the entrance.
Construct the couplings for the young entrances weighted by $w_g$ in a
fixed order, starting each only if its entrance is still reserved.
Complete each coupling before starting the next.
A spoiled entrance gives an immediate defect.

For each coupling, retain the entrance, the root type and its age (its
label divided by $n$), together with the entire preceding exploration
record. Initialize the root mark and apply the collision rules of
Section~\ref{subsec:ppt-pa-coupling}.
Defer couplings at the reserved young entrances of required hubs until
the audits in Section~\ref{subsec:hub-vertices}; intervening disclosures
may still spoil them.
Write $\mathscr H_\partial$ for the completed boundary record.

With the graph and assigned coupled environments fixed, put
$b_g=r_g-q_g(V(G))$, the probability of a defect or a penetration
with a rejected continuation. The retained region, its recorded
incident slots, and the finite hitting-probability rule determine
$b_g$. Retain these data before the next coupling.

Set $b_g=1$ at spoiled entrances and define
\[
 B_T:=\sum_gw_gb_g+\sum_y\frac{\eta(y)c_y}{\deg_G(y)-1}.
\]
Let $\mathcal G_{\mathrm{hub}}$ be the event that every required hub
procedure has quenched failure probability at most $\varepsilon/2$,
uniformly over these original entrances. For roots that are not good,
construct no couplings, set $B_T=0$, and declare $\mathcal G_{\mathrm{hub}}$
to hold.

For a good root, the terminal rules and
\eqref{eq:boundary-hitting-comparison} give
\[
 P_{x_0}^G(\mathcal E^c)\le C_mB_T+\varepsilon/2
 \qquad\text{on }\mathcal G_{\mathrm{hub}}.
\]
Indeed, at a young terminal the failure probability conditional on
resolution is at most $b_g/p$. At an old nonhub terminal $y$, summing
excursions that return to $y$ shows that, conditional on resolution,
the probability of failure without starting a hub procedure is at most
\[
 \frac{c_y+\sum_{g\text{ at }y}b_g}{p(\deg_G(y)-1)},
\]
since each forward excursion resolves with probability at least
$p$. Multiplying by the resolving-terminal probabilities and using
\eqref{eq:boundary-hitting-comparison} gives $C_mB_T$, with $p^{-1}$
absorbed into $C_m$. The strong Markov property at the first hub
procedure bounds the remaining failures by $\varepsilon/2$.
Consequently, \eqref{eq:boundary-input} follows once we prove, for
every fixed $\varepsilon_B>0$ and uniformly in $x_0$,
\[
 \mathbb P(\Omega\cap\mathcal U_n,\ B_T>\varepsilon_B)=o(n^{-1}),\qquad
 \mathbb P(\Omega\cap\mathcal U_n\cap\mathcal G_{\mathrm{hub}}^c)
 =o(n^{-1}),
\]
and choose $C_m\varepsilon_B\le\varepsilon/2$. We prove the first estimate
here and the second in Section~\ref{subsec:hub-vertices}.

We impose $|E|\le M_n$, with $M_n$ as in \eqref{eq:coupling-scales},
with $c_M<1/2$, $c_r\ge c_R$, and $c_0>A$, so $k_n^\circ<a$.
This bound includes all labels added before the couplings and counts vertices
whether or not their complete stars have been recorded.
Every required hub is in $T$ or adjacent to an old
nonhub terminal, and at most $C(D+1)(|T|+1)$ old vertices have their
status queried. On $\Omega\cap\mathcal U_n$, completing their stars
adds at most $C(D+1)(|T|+1)n^{3/4}=n^{3/4+o(1)}$ labels.
Each hub samples $k=\lceil(\log n)^2\rceil$ slots.
Thus the number of boundary couplings plus the number of slots sampled at hubs is at most
\begin{equation}\label{eq:boundary-test-count}
 Ck(D+1)(|T|+1)\le e^{C_{\mathrm{test}}L^2},
\end{equation}
for any fixed $C_{\mathrm{test}}>C_{\mathrm{vol}}$ and large $n$.
Choose this constant before the regularity constants in
Section~\ref{subsec:ppt-pa-coupling}, so that
\eqref{eq:repeated-extension-error} applies. By
\eqref{eq:walk-registration-budget}, the size of $E$, including all
labels added before the couplings, is at most
\begin{equation}\label{eq:boundary-registration-count}
 n\exp\{-\tfrac12\sqrt{\log n}+CL^2\}=o(M_n).
\end{equation}
Stop the exploration if either budget is exceeded. Before initiating
an entrance coupling, also stop if $|E|>M_n/2$.
This coupling supplies the environment for the forward test and its
continuation. The count budget is not exceeded on $\Omega$;
\eqref{eq:boundary-registration-count} rules out the stops imposed
on $|E|$ on $\Omega\cap\mathcal U_n$.

On $\mathcal U_n$, the initial PA record completed before the couplings
satisfies the history conditions of Section~\ref{subsec:ppt-pa-coupling}.
Every queried target in $E_0$ has its complete star recorded.
Put $t=|E|$.
Their nonhub degrees sum to at most $tD$, and their hub degrees to
at most $2mt$: all senders have labels in $E$ and at most $m$ birth
slots each, and every edge contributes at most twice to the degree sum.
Thus \eqref{eq:degree-mass-regularity} and \eqref{eq:good-urn-event} give
\[
 \sum_{\ell\in E_0}\vartheta_\ell
 \le\frac{CtD}{n(k_n^\circ/n)^\chi}
 =\frac{Ct}{k_n^\circ}D(k_n^\circ/n)^\beta
 =o(t\bar g_n/k_n^\circ).
\]
Here $D(k_n^\circ/n)^\beta/\bar g_n
=O((\log n)^2e^{-\beta c_1L^2})=o(1)$, verifying the aggregate cap.
Before starting the couplings, check this observable bound and
$j\psi_j\le\bar g_n$ for every $j\ge k_n^\circ$ whose mark has been exposed.
Check the latter bound at every subsequent mark exposure. A failed check stops the entire
exploration; these stops cannot occur on $\mathcal U_n$. Other mark
and comparison failures remain local, with all PA answers retained. No new low target is queried during the couplings.
At any of the preceding stops, set $B_T=1$ and declare
$\mathcal G_{\mathrm{hub}}^c$; subsequent charges and loss variables
in the estimates are zero. We bound the losses on
$\Omega\cap\mathcal U_n$, where these stops do not occur. The common
$\mathcal U_n^c$ is excluded only in
Proposition~\ref{prop:worst-start-reduction}.

We first bound old nonhub departures revealed while exploring the old
nonhub terminals, before the deferred stars are completed. Fix a
record of $T$ with $\eta(v)\le w_*$ at every nonhub terminal, as
holds on $\Omega$. Let $\mathscr G_j$ be the record
after $j$ queries in this exploration, enlarged by
$(\psi,\mathcal C_n)$. The bounds
\eqref{eq:common-red-event} are measurable in $\mathscr G_0$;
throughout this calculation suppose they hold.
During this exploration, for a slot $(k,i)$ whose endpoint has not been
revealed and a nonhub target $v<k$ whose complete star has not been recorded,
\begin{equation}\label{eq:red-query-bound}
 \mathbb P(J_{k,i}=v\mid\mathscr G_{j-1})
 \le\frac{\varphi_v}{\sum_{h<k,\,h\in\mathcal C_n}\varphi_h}
 \le C(\log n)^3\varphi_v
       \left(\frac n{\min\{k,a\}}\right)^\chi.
\end{equation}
The probability is zero if $v$ has already been excluded.
Otherwise, redirecting the slot from $v$ to a hub $h<k$ preserves
the record: $v$ remains a nonhub, $h$ remains a hub, and no negative
query has a hub as its target. The likelihood ratio is
$\varphi_h/\varphi_v$, proving the first inequality. For the second,
put $r=\min\{k-1,a-1\}$ and use
\[
 \sum_{h<k,\,h\in\mathcal C_n}\varphi_h
 \ge S_r-CD r/n
 \ge c(\log n)^{-3}(\min\{k,a\}/n)^\chi,
\]
since $(\log n)^3D(a/n)^\beta=o(1)$. We also have
\begin{equation}\label{eq:red-atom-envelope}
 \varphi_v\le C(\log n)^4n^{-\chi}v^{-\beta},\qquad
 \sum_{v\in F}\varphi_v\le C(\log n)^4(|F|/n)^\chi.
\end{equation}

Every original boundary entrance is reached from $x_0$ by a nonhub
path of length at most $R+1$. An additional incidence at its endpoint
disclosed before the first PA--PPT coupling therefore lies in
$B_{5R}^{\mathrm{stop}}(x_0)$, including incidences found in a hub's
star. On $\Omega$, at most one such edge lies outside the breadth-first
tree. It spoils at most two of these entrances, with total weight at most
$2w_*$.

Before exploring the old nonhub terminals, assign
$\lambda_y=\eta(y)/(m-1)$ to each such terminal, and zero to other
labels. At the first disclosure of an edge joining two old nonhub
vertices, charge $\lambda_y$ for
each endpoint $y$ at which this edge is a forward slot. These
charges dominate $\sum_y\eta(y)c_y/(\deg_G(y)-1)$, and each is at
most $2w_*$. For $r<s<a$, \eqref{eq:red-query-bound} is at most
$C(\log n)^3Dn^{-\beta}s^{-\chi}$. Since each slot--target pair
is queried at most once, the sum of the conditional expected
charges is at most
\[
 C(\log n)^3Dn^{-\beta}
 \sum_{r<s<a}(\lambda_r+\lambda_s)s^{-\chi}
 \le C(\log n)^3D(a/n)^\beta,
\]
where we used $\sum_y\lambda_y\le1/(m-1)$.

During the subsequent couplings, condition on the urn variables and the
current full exploration record. The degree classifications and all hub
neighbor lists are determined by its slot answers; samples from these
lists use independent auxiliary randomness. For a sender $u\ge a$,
the total relative mass of
previously excluded targets is at most
\[
 C(\log n)^7(M_n/a)^\chi=o(1).
\]
The conditional query bound from Section~\ref{subsec:ppt-pa-coupling}
therefore gives, for each slot $(u,i)$ whose endpoint has not been revealed
and every target $x<u$, including hubs,
\[
 \mathbb P(J_{u,i}=x\mid\psi,\text{record})
 \le\frac{2\varphi_x}{S_{u-1}}
 \le C(\log n)^3\varphi_x(n/u)^\chi.
\]
Before the couplings begin, the unspoiled entrances weighted by $w_g$
have distinct young endpoints. Charge an entrance's weight when another
incidence at its endpoint is first revealed before its own coupling begins.
We do not explore a reserved endpoint before its own coupling begins,
so a full endpoint disclosure can spoil its entrance only by reaching
that vertex as a target. The preceding bound and
\eqref{eq:red-atom-envelope} give expected charge at most
$C(\log n)^7/a$ per disclosure; there are at most $mM_n$ disclosures.
Alternatively, a younger-neighbor
query $J_{u,i}=x$ may reveal an incidence at a reserved endpoint
$u>x$. For its entrance $g$, sum over the $m$ slots of $u$ and the
at most $M_n$ distinct targets in $E$ queried before its coupling begins. The
conditional expected charges total at most
\[
 C(\log n)^3w_g(n/a)^\chi
 \max_{|F|\le M_n}\sum_{x\in F}\varphi_x
 \le C(\log n)^7w_g(M_n/a)^\chi.
\]
The bound holds for every target set and hence also for adaptive
queries. Summing over entrances uses $\sum_gw_g\le1$.

Apply the exponential-moment bound separately to the charges from
exploring old nonhub terminals and to the subsequent spoiled-entrance
charges. In either case, write $\Delta_j$ for the charge and
$\mathscr G_{j-1}$ for the corresponding conditioning record.
The preceding estimates give
\[
 0\le\Delta_j\le2w_*,\qquad
 \sum_j\mathbb E[\Delta_j\mid\mathscr G_{j-1}]
 \le C\bigl[(\log n)^3D(a/n)^\beta
       +(\log n)^7\{M_n/a+(M_n/a)^\chi\}\bigr]=o(1).
\]
Here $M_n/a=\exp\{-c_M\sqrt{\log n}+O(L^2)\}$.
The inequality $e^x\le1+(e-1)x$ for $0\le x\le1$ gives
\[
 \mathbb E[e^{\Delta_j/(2w_*)}\mid\mathscr G_{j-1}]
 \le\exp\left\{\frac{e-1}{2w_*}
             \mathbb E[\Delta_j\mid\mathscr G_{j-1}]\right\}.
\]
For each sum, iteration and exponential Markov inequality bound the
probability of exceeding $\varepsilon_B/4-2w_*$ by $e^{-c_{\varepsilon_B}/w_*}$
for large $n$. Integrating over the respective
initial records and including the at most $2w_*$ weight of entrances
spoiled before the couplings gives
\[
 \mathbb P\left(\Omega\cap\mathcal U_n,\
   \sum_y\frac{\eta(y)c_y}{\deg_G(y)-1}
       +\sum_{g\text{ spoiled}}w_g>\varepsilon_B/2\right)
 \le 2e^{-c_{\varepsilon_B}/w_*}.
\]

For the remaining contribution, let $\mathscr F_g^-$ be the
full record immediately before the coupling at $g$, including its
entrance, all revealed PA and PPT data, and auxiliary choices already
made, including those from failed couplings. No conditioning on the full
graph, the urn, or $\mathcal U_n$ is added. Uniformly over these
records and couplings started at reserved young entrances, including eligible hub samples,
\begin{equation}\label{eq:hub-young-conditional-input}
 \mathbb E[b_g\mid\mathscr F_g^-]\le\varepsilon_Y
\end{equation}
for any fixed $\varepsilon_Y>0$, after choosing the constants and
taking $n$ sufficiently large. Since the completed coupling record
determines the quenched loss $b_g$, the tower property identifies
this conditional expectation with the annealed test-loss probability.

Indeed, the entrance root has only its recorded entrance and an
unexposed mark. Its age is its label divided by $n$, at least
$e^{-AL^2}$, and its initial discrepancy is zero. Initialize its
type-dependent Gamma mark by \eqref{eq:joint-mark-error}, declaring
an immediate defect on failure. The initialization error, mark
goodness failure, and event $\Gamma_\rho>CL^2$ have total conditional
probability $o(1)$ for a fixed $C>0$. Outside these events,
the pretrial check $|E|\le M_n/2$ and the retained record give the
hypotheses of Lemma~\ref{lem:local-continuation}, with testing
depth $R\le r_n$. Apply \eqref{eq:local-continuation-loss} with
its error parameter chosen so that the bound is at most
$\varepsilon_Y/2$. Averaging over initialization gives
$\mathbb E[b_g\mid\mathscr F_g^-]\le\varepsilon_Y/2+o(1)$,
which proves \eqref{eq:hub-young-conditional-input}.

Condition on an initial PA record $\mathscr H$ completed before the
couplings and satisfying the displayed weight bounds. Order the initiated
boundary couplings by execution time, padding unused positions with zero contributions,
and include all intermediate disclosures in the filtration.
The weights are known before each coupling. Put
$\mu_g=\mathbb E[b_g\mid\mathscr F_g^-]$.
Conditional Hoeffding gives
\[
 \mathbb E\left[e^{\lambda w_g(b_g-\mu_g)}
                  \mid\mathscr F_g^-\right]
 \le e^{\lambda^2w_g^2/8}.
\]
Since $\sum_gw_g\mu_g\le\varepsilon_Y$ and
$\sum_gw_g^2\le w_*$, the exponential martingale bound, averaged
over these initial records, yields
\[
 \mathbb P\left(\Omega\cap\mathcal U_n,\
       \sum_{g\text{ with a coupling}}w_gb_g>\varepsilon_Y+t\right)
 \le e^{-2t^2/w_*}.
\]
Choose $\varepsilon_Y\le\varepsilon_B/4$ and $t=\varepsilon_B/4$. Combining this
with the bounds for old departures and spoiled entrances gives
\[
 \mathbb P(\Omega\cap\mathcal U_n,\ B_T>\varepsilon_B)
 \le3e^{-c_{\varepsilon_B}/w_*}
 =o(n^{-1}),
\]
provided $c_Rc_{\rm bd}>1$, since
$w_*^{-1}=e^{c_{\rm bd}R}\ge(\log n)^{c_Rc_{\rm bd}}$.

\subsection{Hub continuations}
\label{subsec:hub-vertices}
We prove
\begin{equation}
 \begin{split}
    \mathbb P(\Omega\cap\mathcal U_n\cap\mathcal G_{\mathrm{hub}}^c)
       &\le e^{-c(\log n)^{5/2}}+Ce^{CL^2-c_\varepsilon k}\\
       &=o(n^{-1}).
 \end{split}
    \label{eq:hub-continuation-input}
\end{equation}
After completing the boundary couplings, we audit the required
hubs sequentially. Their complete incident-slot lists are already
recorded. At each hub, sample
$k=\lceil(\log n)^2\rceil$ distinct incident slots uniformly
without replacement, completing each eligible coupling before taking
the next sample. Write $e_j$ for the $j$th sampled slot.
A sample is eligible if its endpoint is at least
$b$ and is still reserved for that entrance immediately before sampling;
rejected samples are not replaced.

Fix a required hub $h$, write $d=d_G(h)>D$, and let $\mathscr F_0$
be the full exploration record at the start of its audit. Write $\mathbf Z=(Z_v)_{v\in[n]}$
for the PA marks. For the fixed realization $(G,\psi,\mathbf Z)$,
call an incident slot \emph{good} if its endpoint is at
least $b$, is reserved in $\mathscr F_0$, and admits a coupled forward test and
continuation compatible with this realization and record for which
$b_e\le\varepsilon r_e/4$.
Let $\theta_h$ be the fraction of nongood slots. Conditional on
$(G,\psi,\mathbf Z,\mathscr F_0)$, this population is fixed.

\medskip
\noindent\emph{Step 1: Comparison with the sample mean.}
Let $\widehat\theta_h$ be the fraction of nongood slots among
the $k$ samples. Conditional on $(G,\psi,\mathbf Z,\mathscr F_0)$,
the sample is uniform without replacement. Hoeffding's inequality gives
\[
    \mathbb P\left(
        \theta_h>\widehat\theta_h+\frac{p\varepsilon}{8}
        \,\middle|\,G,\psi,\mathbf Z,\mathscr F_0
    \right)
    \le e^{-c_\varepsilon k}.
\]
Averaging over the PA realization gives the same bound conditional on
$\mathscr F_0$. If a stop
interrupts the audit, complete the sample uniformly from the
remaining slots only for this calculation; the additional
slots are neither recorded nor used for couplings.

Thus, to bound $\theta_h$ by $p\varepsilon/4$, it remains to
bound $\widehat\theta_h$ by $p\varepsilon/8$.

\smallskip
\noindent\emph{Step 2: Control of the sample mean.}
Let $R_j$ indicate that sample $j$ is rejected because its
endpoint is below $b$ or is not reserved. Let $D_j$ indicate
an eligible sample with $b_{e_j}>\varepsilon r_{e_j}/4$, with $D_j=0$ on rejected
samples.
After a global stop, keep the full record fixed and set
$R_j=D_j=0$ for the remaining indices.
An eligible sample satisfying $b_{e_j}\le\varepsilon r_{e_j}/4$ witnesses
goodness relative to $\mathscr F_0$, since its coupling respects that record.
Thus, for a completed audit, every sampled nongood slot is counted
by $R_j+D_j$, and
\[
    \widehat\theta_h
    \le \frac1k\sum_{j=1}^k R_j
       +\frac1k\sum_{j=1}^k D_j.
\]
We first control entrances spoiled by the additional old stars.
For an old vertex $h$, let $K_h$ count slots $u\to h$, $u\ge b$,
for which another birth slot of $u$ has an old target, including $h$.
Conditional on urn variables satisfying the first two bounds in
\eqref{eq:common-red-event}, the contributions from distinct senders
are independent and bounded by $m$, and
\[
 \mathbb E[K_h\mid\psi]
 \le Cn\varphi_h(\log n)^C(a/b)^\chi
 \le n\varphi_h e^{-cL^2}.
\]
On $\mathcal U_n$, a hub satisfies
$d_G(h)\ge\max\{D,cn\varphi_h\}$. Chernoff's inequality at the
urn-measurable threshold $\max\{D,cn\varphi_h\}/(2\log n)$,
followed by a union bound over all old vertices, gives
\[
 \mathbb P\left(\Omega\cap\mathcal U_n,\
   \exists h\in\mathcal C_n:K_h>d_G(h)/(2\log n)\right)
 \le e^{-c(\log n)^2L^2}.
\]
On $\Omega$, at most two slots at a required hub lead to young
vertices already in $T$: three such slots would give two cycles in
$B_{5R}^{\mathrm{stop}}(x_0)$.
Apart from these slots and those counted by $K_h$, every endpoint
at least $b$ disclosed at that hub is reserved before the first
PA--PPT coupling.

Fix the initial PA record just before the first PA--PPT coupling. Track
the endpoints $v\ge b$ reserved for entrances at $h$, removing each
when its entrance is spoiled or its own coupling begins.
While an endpoint $v$ remains reserved,
its mark is unexposed and $s_v=0$ in the posterior calculation of
Section~\ref{subsec:ppt-pa-coupling}, since $v\to h$ is outgoing.
Thus $\mathbb E[\psi_v\mid\mathscr F]\le C/v$ for the current full
exploration record $\mathscr F$. The excluded-mass bound gives
\begin{equation}\label{eq:hub-conditional-ordinary-slot}
 \begin{split}
 \mathbb P(\text{an endpoint disclosure hits the reserved set}
               \mid\mathscr F)&\le Cd/b,\\
 \mathbb P(J_{v,i}=x\mid\mathscr F)&\le C\bar g_n/k_n^\circ
 \end{split}
\end{equation}
for a birth slot $(v,i)$ whose endpoint has not been revealed, with $v$
reserved and $x\ge k_n^\circ$ a queried target. The latter target has its capped mark exposed.
Keep the process stopped at the record caps and the bound $|E|\le M_n$.
There are at most $mM_n$ endpoint disclosures and $M_n$ distinct
queried targets, each tested against at most $md$ birth slots of
these reserved endpoints. Let $N_h$ count only removals caused by spoiling.
The sum of its conditional expected increments is at most
$CdM_n\bar g_n/k_n^\circ$. The at most $k$ removals at the start of
couplings sampled at $h$ are counted separately.
Successive slot answers give increments at most one, so
the conditional exponential-moment bound yields, uniformly over the
completed initial records,
\[
 \mathbb P\left(N_h>d/(2\log n)\mid\text{initial record}\right)
 \le\left(\frac{CM_n\bar g_n\log n}{k_n^\circ}
      \right)^{c d/\log n}
 \le e^{-c(\log n)^{5/2}}.
\]
Here $d>D$ and
$M_n\bar g_n\log n/k_n^\circ
=\exp\{-c_M\sqrt{\log n}+O(L^2)\}$.

Let $\mathcal Q$ be the event that, before every hub sample, the
fraction of slots leading below $b$ or to an endpoint not
reserved is at most $\zeta:=\rho+(2+k)/D+1/\log n=o(1)$.
On $\Omega$, the low endpoints contribute at most $\rho$.
The preceding bounds and the at most $k$ samples at a hub imply
\begin{equation}\label{eq:hub-contact-exception}
 \mathbb P(\Omega\cap\mathcal U_n\cap\mathcal Q^c)
 \le e^{-c(\log n)^{5/2}}+e^{-c(\log n)^2L^2},
\end{equation}
after summing over the at most $e^{CL^2}$ required hubs.

Let $\mathscr F_j$ be the full record after sample $j$ and any
associated coupling.
Conditional on $(G,\psi,\mathbf Z,\mathscr F_{j-1})$, the next sample is
uniform among the
$d-j+1$ remaining slots. On $\mathcal Q$,
\[
    \mathbb E[R_j\mid G,\psi,\mathbf Z,\mathscr F_{j-1}]
    \le \frac{\zeta d}{d-j+1}
    \le 2\zeta
    \le \frac{p\varepsilon}{32}
\]
for sufficiently large $n$. Therefore, on $\mathcal Q$,
the event $\sum_{j=1}^kR_j>p\varepsilon k/16$ forces the
sum of the centered rejection indicators to exceed
$p\varepsilon k/32$. Azuma--Hoeffding, conditional on $(G,\psi,\mathbf Z)$
and then averaged over this realization, gives
\[
    \mathbb P\left(
        \mathcal Q,\
        \sum_{j=1}^kR_j>\frac{p\varepsilon k}{16}
        \,\middle|\,\mathscr F_0
    \right)
    \le e^{-c_\varepsilon k}.
\]
The event $\mathcal Q$ is intersected with the deviation
event, not added to the conditioning.

For $D_j$, let $\mathscr F_j^-$ be the full record after identifying
sample $j$ and before its coupling, without further conditioning on
$(G,\psi,\mathbf Z)$ or $\mathcal U_n$. Since $r_{e_j}\ge p$,
we have $D_j\le4b_{e_j}/(p\varepsilon)$ on eligible samples.
By \eqref{eq:hub-young-conditional-input},
\[
    \mathbb E[D_j\mid\mathscr F_j^-]
    \le \frac{4\varepsilon_Y}{p\varepsilon}(1-R_j).
\]
Choose $\varepsilon_Y\le p^2\varepsilon^2/512$, so these conditional
means are at most $p\varepsilon/64$. The increments
$D_j-\mathbb E[D_j\mid\mathscr F_j^-]$ are martingale
differences by the tower property. Azuma--Hoeffding gives
\[
    \mathbb P\left(
        \sum_{j=1}^kD_j>\frac{p\varepsilon k}{16}
        \,\middle|\,\mathscr F_0
    \right)
    \le e^{-c_\varepsilon k}.
\]
Combining this bound with the rejected-sample bound, the domination
of $\widehat\theta_h$, and Step~1 gives
\[
    \mathbb P\left(
        \Omega\cap\mathcal U_n\cap\mathcal Q,\
        \theta_h>\frac{p\varepsilon}{4}
        \,\middle|\,\mathscr F_0
    \right)
    \le Ce^{-c_\varepsilon k}.
\]
Enumerate the audits in execution order, assigning zero contributions
to unused positions. The bound~\eqref{eq:boundary-test-count} gives at most
$\exp(CL^2)$ positions, and on $\Omega\cap\mathcal U_n$ every
required hub is audited.
Average over each starting record and take a union bound to obtain
\[
    \mathbb P\left(
        \Omega\cap\mathcal U_n\cap\mathcal Q,\
        \theta_h>\frac{p\varepsilon}{4}
        \text{ for some required hub }h
    \right)
    \le Ce^{CL^2-c_\varepsilon k}.
\]
Finally, suppose $\theta_h\le p\varepsilon/4$. For each good slot,
sampled or not, fix separately a witnessing coupled forward test and continuation
compatible with $(G,\psi,\mathbf Z,\mathscr F_0)$. These choices
need not arise from one simultaneous adaptive exploration and
make no further PA revelations in the audit. Each chosen finite
record defines a walk procedure with the resolution, duration,
and endpoint bounds of Lemma~\ref{lem:local-continuation}, and is
kept fixed on repeated visits.

Declare a departure through a nongood slot an immediate defect,
setting $b_e=r_e=1$ there. Summing the excursions that return to
$h$, the probability of eventual failure is
\[
    \frac{\sum_e b_e}{\sum_e r_e}
    \le \frac{\varepsilon}{4}+\frac{\theta_h}{p}
    \le \frac{\varepsilon}{2},
\]
where the sums run over all incident slots and $r_e\ge p$.
Each original entrance to $h$ leads either to an old vertex or to
a vertex already explored in the initial stopped exploration, so
it is excluded from the good slots. The bound is therefore uniform
over these entrances. Combining the preceding bounds with
\eqref{eq:hub-contact-exception} proves
\eqref{eq:hub-continuation-input}.

Together with the boundary estimate of
Section~\ref{subsec:boundary-losses}, this proves
\eqref{eq:boundary-input}.

\medskip

\medskip
%%%%%%%%%%%%%%%%%%%%%%%%%%%%%%%%%%%%%%%%%%%%%%%%%%%%%%%%%%%%%%%%%%%%%%%%%%%%%%%%%%%%%%%%%%%%

\section{Lower bound on the mixing time}\label{sec:cutoff-proof}

Let $V_0$ be uniform on $[n]$, independently of $G_n$. For $\gamma>0$, set
\[
 t_-(\gamma):=\left\lfloor
       \frac{\log n-3\gamma\sqrt{\log n}}{\nu d}\right\rfloor.
\]
We use the retained subgraph from
Proposition~\ref{prop:walk-coupling-log-time} to localize the walk on a
set of vanishing stationary mass.

\begin{prop}\label{prop:cutoff-lower-bound}
For every $\varepsilon\in(0,1)$,
\begin{equation}
 \lim_{\gamma\to\infty}\limsup_{n\to\infty}
 \Pbb_{G_n,V_0}\left(
       \left\|P_n^{t_-(\gamma)}(V_0,\cdot)-\pi_n\right\|_{\rm TV}
                      <1-\varepsilon\right)=0.
 \label{eq:lower-bound-typical-root}
\end{equation}
Consequently,
\begin{equation}
 \lim_{\gamma\to\infty}\limsup_{n\to\infty}
 \Pbb_{G_n}\left(d_n(t_-(\gamma))<1-\varepsilon\right)=0.
 \label{eq:lower-bound-worst-start}
\end{equation}
\end{prop}

\begin{proof}
Fix $\gamma>0$ and write $t=t_-(\gamma)$. Take $A_\rho$ uniform on
$(0,1)$, independently of the PA construction, and set
$V_0=\lceil nA_\rho\rceil$ and $\Gamma_\rho=Z_{V_0}$.
These have the required PA and PPT root laws. On $\{A_\rho\ge b_n\}$,
the root satisfies the age and label cutoffs, with discrepancy
$\Delta_\rho\le(nb_n)^{-1}=o(\varepsilon_n)$.
Lemma~\ref{lem:beta-gamma} and the Gamma tail bound give the mark
requirements except with probability $o(1)$. Since
$\Pbb(A_\rho<b_n)=b_n=o(1)$, initialization fails with probability $o(1)$.

Run the construction of Proposition~\ref{prop:walk-coupling-log-time}
without an entrance, setting $\widetilde\Gamma=\varnothing$ on
initialization failure, and complete both environments with their
conditional laws. Keep the notation $X$, $\Pi_s$, and $\tau$ from that
construction. Since $\Pi_t^c\subseteq\Pi_\tau^c\cup\{\tau<t\}$,
the coupling proposition and Proposition~\ref{prop:no-truncation}
bound failure by $\Pbb(\tau<t)+\eta(\gamma)+o(1)$, after
averaging over the root data and including the initialization and
urn exceptions. The early hit requires
$\max_{s\le t}(|X_s|-\nu s)\ge\ell_1-\nu t
=2\gamma\sqrt{\log n}/d+O(1)$.
Thus \eqref{eq:maximal-distance-tightness}, enlarging $\eta$ if necessary,
gives
\begin{equation}
 \Pbb(\Pi_t^c)
 \le \Pbb(\Pi_\tau^c)+\Pbb(\tau<t)
 \le \eta(\gamma)+o(1),
 \label{eq:lower-bound-coupling-loss}
\end{equation}
where $\eta(\gamma)\to0$ as $\gamma\to\infty$ and $o(1)$ is taken as
$n\to\infty$ with $\gamma$ fixed. These probabilities average over
the environments, root, and walks.

Let $F=V(\widetilde\Gamma)$. Only explored vertices whose complete stars
were successfully compared belong to $F$; truncated and mismatched boundary
vertices are excluded. The size bound
\eqref{eq:truncated-exploration-size} and the bound of $B_n$ labels
added to $E$ per comparison in
Section~\ref{sec:coupling-after-burn-in} therefore give
\[
 |F|\le(t_n+1)e^{\omega_n}V_n,
 \qquad
 \deg_{G_n}(v)\le CB_n\quad(v\in F).
\]
Consequently, using $\pi_n(v)=\deg_{G_n}(v)/(2m(n-1))$,
\begin{equation}
 \pi_n(F)
 \le\frac{CB_n(t_n+1)e^{\omega_n}V_n}{2m(n-1)}
 =\exp\left\{-\tfrac12\sqrt{\log n}
                  +O((\log\log n)^2)\right\}=o(1).
 \label{eq:lower-bound-stationary-mass}
\end{equation}
This bound is uniform over all outcomes of the construction; the
degrees are those in the full graph $G_n$.

Let $\mathcal E_n$ be the sigma-field generated by the completed
environments, the root, the retained subgraph $\widetilde\Gamma$, and
its matching $\Phi$. The set $F$ is $\mathcal E_n$-measurable.
Sample the coupled walks using transition randomness independent of
$\mathcal E_n$. Conditional on $\mathcal E_n$, the finite walk has
kernel $P_n$ from $V_0$; on $\Pi_t$, its endpoint lies in $F$. Hence
\begin{align*}
 \left\|P_n^t(V_0,\cdot)-\pi_n\right\|_{\rm TV}
 \ge P_n^t(V_0,F)-\pi_n(F)
 \ge 1-\Pbb(\Pi_t^c\mid\mathcal E_n)-\pi_n(F).
\end{align*}
For fixed $\varepsilon\in(0,1)$ and all sufficiently large $n$,
\eqref{eq:lower-bound-stationary-mass} gives $\pi_n(F)\le\varepsilon/2$.
Markov's inequality and \eqref{eq:lower-bound-coupling-loss} yield
\begin{align*}
 \Pbb_{G_n,V_0}\left(
       \left\|P_n^t(V_0,\cdot)-\pi_n\right\|_{\rm TV}
                         <1-\varepsilon\right)
 \le \frac{2}{\varepsilon}\Pbb(\Pi_t^c)
 \le \frac{2\eta(\gamma)}{\varepsilon}+o(1).
\end{align*}
Taking $n\to\infty$ and then $\gamma\to\infty$ proves
\eqref{eq:lower-bound-typical-root}.

The worst-start bound \eqref{eq:lower-bound-worst-start} follows from
$d_n(t)\ge\|P_n^t(V_0,\cdot)-\pi_n\|_{\rm TV}$.
\end{proof}

\begin{proof}[Proof of Theorem~\ref{thm:main-cutoff}]
Fix $\varepsilon\in(0,1)$. Proposition~\ref{prop:cutoff-lower-bound}
gives the first limit in \eqref{eq:cutoff-def} after rescaling $\gamma$.
For the second limit, apply
Proposition~\ref{prop:worst-start-reduction} and the boundary estimates
of Section~\ref{sec:worst-start} with error $\varepsilon/2$.
The resulting time is $\lceil t_n^\star\rceil+R^2+5R+s_n$.
As $R=O(\log\log n)$ and $s_n=O_\varepsilon(\sqrt{\log n})$,
there is a constant $C_\varepsilon<\infty$ such that
\[
 \Pbb\bigl(d_n(t_n^\star+C_\varepsilon\sqrt{\log n})
               \le\varepsilon\bigr)\longrightarrow1.
\]
Monotonicity gives the second limit in \eqref{eq:cutoff-def}.
\end{proof}

\subsection*{Declaration of generative AI use}
The authors used generative AI tools solely to assist with checking notation, identifying potential errors in proofs, and preparing figures. All proofs were developed by the authors, who take full responsibility for the final manuscript.

\bibliographystyle{plain} % 또는 alpha, unsrt, abbrv 등
\bibliography{SRW.bib}  

\end{document}